\magnification=\magstephalf
\input stpel.mac
\let\noarrow = t
\input eplain
\indexproofingtrue
\newcount\refno\refno=1
\def\incr{\advance\refno by 1}
\edef\BlaRog{\number\refno}\incr
\edef\ColFon{\number\refno}\incr
\edef\Conr{\number\refno}\incr
\edef\DelCoord{\number\refno}\incr
\edef\Falt{\number\refno}\incr
\edef\FaCh{\number\refno}\incr
\edef\Font{\number\refno}\incr
\edef\GorOo{\number\refno}\incr
\edef\Illus{\number\refno}\incr
\edef\Katz{\number\refno}\incr
\edef\KatzST{\number\refno}\incr
\edef\Knus{\number\refno}\incr
\edef\KottTOI{\number\refno}\incr
\edef\KottIsoc{\number\refno}\incr
\edef\Kott{\number\refno}\incr
\edef\Kraft{\number\refno}\incr
\edef\LST{\number\refno}\incr
\edef\Manin{\number\refno}\incr
\edef\Mess{\number\refno}\incr
\edef\Milne{\number\refno}\incr
\edef\Durh{\number\refno}\incr
\edef\GSAS{\number\refno}\incr
\edef\DFEO{\number\refno}\incr
\edef\MAV{\number\refno}\incr
\edef\Noot{\number\refno}\incr
\edef\FOTexel{\number\refno}\incr
\edef\FOTexelB{\number\refno}\incr
\edef\RapBourb{\number\refno}\incr
\edef\RR{\number\refno}\incr
\edef\Rayn{\number\refno}\incr
\edef\ReimZink{\number\refno}\incr
\edef\Saav{\number\refno}\incr
\edef\Stamm{\number\refno}\incr
\edef\WedhOrd{\number\refno}\incr
\edef\WedhCong{\number\refno}\incr
\edef\Wedh{\number\refno}\incr

%% \nopagenumbers
%% \input Abstract
%% \line{\hfill}
%% \vfill\eject
%% \footline={\hss\tenrm\folio\hss}
%% \pageno=1

\centerline{{\sectitlefont SERRE-TATE THEORY FOR MODULI SPACES OF PEL TYPE}}
\Askip

\centerline{{\it by}\quad{\namefont Ben Moonen\footnote{*}{{\eightrm Research 
made possible by a fellowship of the Royal Netherlands Academy of Arts and 
Sciences}}\phantom{\quad{\it by}}}}
\Askip

\def\leaderfill{\leaders\hbox to.8em{\hss.\hss}\hfill}

\readtocfile
\Askip

\introduction

\introsec{intr1}
Let $A$ be an abelian variety over a perfect field~$K$ of characteristic~$p>0$. The ``Serre-Tate theory'' in the title of this paper 
refers to a collection of results about the formal deformations of~$A$ in case $A$ is {\it ordinary\/}. The first main results are 
described in the Woods Hole report of Lubin, Serre and Tate~[\LST]. Related results were obtained by Dwork and were shown to agree with 
those of Serre and Tate. As references we cite Deligne~[\DelCoord], Katz~[\KatzST], and Messing~[\Mess].

One of the main points in Serre-Tate theory is the statement that, if $A$ is ordinary, its formal deformation space has a canonical 
structure of a formal group (in fact, a formal torus) over~$W(K)$. In particular, this leads to a canonical lifting $A^\can$ 
over~$W(K)$, corresponding to the identity section of the formal group. If $K$ is finite then this canonical lifting is characterized by 
the fact that all endomorphisms of~$A$ lift to endomorphisms of~$A^\can$.

The main goal of this paper is to generalize this theory to moduli spaces of PEL type. Roughly speaking these are moduli spaces for 
triples $\ul{A} = (A,\iota,\lambda)$ where $\iota\colon \cO \to \End(A)$ is an action of a given ring~$\cO$ by endomorphisms and 
$\lambda$ is a polarization. The precise formulation of the moduli problem is somewhat involved; see section~\refn{EOonPEL} for details. Let us remark that one also fixes certain discrete invariants, such as the structure of the tangent space $\Lie(A)$ as an 
$\cO$-module. This invariant, classically referred to as the CM-type, plays an important role in this paper.

\introsec{Noot}
Let $\cA_g$ be the moduli stack of principally polarized abelian varieties. In some cases the classical Serre-Tate theory ``induces'' 
results for Shimura subvarieties of~$\cA_g$. Namely, let~$E$ be a number field and suppose $S \hookrightarrow \cA_g \otimes E$ is an 
irreducible component of a Shimura subvariety. If $v$ is a prime of~$E$ above~$p$, consider the integral model~$\cS \hookrightarrow 
\cA_g \otimes O_{E,v}$ obtained by taking the Zariski closure of~$S$ inside~$\cA_g$. If $x \in \cS \otimes \kappa(v)$ is ordinary as a 
point of~$\cA_g$ we obtain, taking formal completions at~$x$, formal schemes $\gS_x \subset \gA_x$ over~$W\big(\kappa(v)\big)$. Now 
$\gA_x$ has a canonical structure of a formal torus, and it is known that $\gS_x$ is a formal subtorus of~$\gA_x$. (At finitely many 
places~$v$ this requires a slight refinement.) If $S$ is of PEL type this follows from the results of Deligne and Illusie 
in~[\DelCoord]; in the general case this was proven by Noot in~[\Noot]. Thus, $\gS_x$ ``inherits'' a formal group structure 
from~$\gA_x$.

However, in the situation just considered it may happen that the special fibre of~$\cS$ does not meet the ordinary locus of~$\cA_g$. In 
that case the previous results give us nothing. To arrive at a meaningful theory we need a new notion of ordinariness. This is where the 
work starts.

\introsec{intrOrd}
{\it Definition of ordinariness.\/} Let $A$ be an abelian variety over a field~$K$ of characteristic~$p$. For simplicity assume that $A$ 
admits a prime-to-$p$ polarization. There are several ways to define what it means for~$A$ to be ordinary. One possible definition is 
based on the classification of Barsotti-Tate groups up to isogeny over an algebraically closed field (Dieudonn\'e, Manin). So, if $K 
\subset k = \kbar$ then $A$ with $\dim(A) = g$ is ordinary if $A_k[p^\infty]$ is isogenous to $(\Qp/\Zp)^g \times \hat\mG_m^g$.

Another approach uses only the $p$-kernel of~$A$. Namely, we have $A[p](k) \cong (\mZ/p\mZ)^f$ for some $f\in \{0,\ldots,g\}$, called 
the $p$-rank of~$A$. Then another definition of ordinariness is given by the condition that the $p$-rank is maximal, i.e., $f=g$. This 
is equivalent to the statement that $A_k[p]$ is isomorphic, as a group scheme, to $(\mZ/p\mZ)^g \times \mu_p^g$.

It is well-known that the above two definitions of ordinariness are equivalent. What is more, if $A$ is ordinary then $A_k[p^\infty]$ is 
even isomorphic to $(\Qp/\Zp)^g \times \hat\mG_m^g$.

The two approaches to ordinariness are best viewed in terms of stratifications of~$\cA_g$. On the one hand, the classification of BT 
(Barsotti-Tate groups) up to isogeny gives rise to a {\it Newton Polygon stratification\/} of~$\cA_g$ in characteristic~$p$. Two moduli 
points are in the same NP-stratum iff the associated BT are isogenous. On the other hand, there is a classification of $p$-kernel group 
schemes over $k = \kbar$ up to isomorphism, due to Kraft. This gives rise to the so-called {\it Ekedahl-Oort stratification\/} of~$\cA_g 
\otimes \Fp$, in which two points are in the same stratum iff the associated $p$-kernel group schemes are isomorphic. For more details 
about these stratifications we refer to Oort, [\FOTexel] and~[\FOTexelB] and Rapoport's report~[\RapBourb]. Although the two 
stratifications are in many respects very different, they each have a unique open stratum (the ordinary locus), and part of what was 
said above can be rephrased by saying that these two open subsets of $\cA_g \otimes \Fp$ are actually the same.

Let us now consider a PEL moduli problem. We restrict our attention to primes of good reduction. Essentially this means that we work in 
characteristic~$p>2$ such that (with $\cO$ as in~\refn{intr1}) $\cO \otimes \Zp$ is a maximal order in a product of matrix algebras over 
unramified extensions of~$\Qp$. The PEL moduli problem is represented by a stack~$\cA_\cD$ which is smooth over a finite unramified 
extension~$O_{E,v}$ of~$\Zp$; here the subscript~$\cD$ refers to the precise data that have been fixed in order to define the moduli 
problem. 

To a triple $\ul{A} = (A,\iota,\lambda)$ we can associate $X := A[p^\infty]$ with its induced action of $\cO \otimes \Zp$ and 
polarization~$\lambda$. A classification theory for such triples $\ul{X} := (X,\iota,\lambda)$ up to isogeny was developed by 
Kottwitz~[\KottIsoc] and Rapoport and Richartz~[\RR]. Their results give rise to a stratification of~$\cA_\cD$; we refer to this as the 
(generalized) NP~stratification. A point $x \in \cA_\cD$ is said to be {\it $\mu$-ordinary\/} if it lies in an open 
(=~maximal-dimensional) NP~stratum.

On the other hand, we may consider $Y := A[p]$ with its induced action of $\cO \otimes \Fp$ and polarization~$\lambda$. A classification 
of such triples $\ul{Y} := (Y,\iota,\lambda)$ was given by the author in~[\GSAS]. (This extends results of Kraft~[\Kraft] and of 
Oort~[\FOTexel] in the Siegel modular case; in the Hilbert modular case our results had previously been obtained by Goren and 
Oort~[\GorOo].) Again this gives rise to a stratification of~$\cA_\cD$, referred to as the (generalized) EO~stratification. See also 
Wedhorn~[\Wedh]. We say that a point $x \in \cA_\cD$ is {\it $[p]$-ordinary\/} if it lies in an open EO~stratum. 

Our results in~[\GSAS] include a completely explicit description of the triples~$\ul{Y}$ that occur, in terms of their Dieudonn\'e 
modules. In particular, we can describe the triple~$\ul{Y}$ that corresponds to the $[p]$-ordinary stratum directly in terms of the 
data~$\cD$ used in the formulation of the moduli problem. What is more, we can also give an explicit Barsotti-Tate group 
$\ul{X}^\ord(\cD)$ with polarization and $\cO \otimes \Zp$-action such that its $p$-kernel is of the $[p]$-ordinary type. We refer 
to~$\ul{X}^\ord(\cD)$ as the {\it standard ordinary object\/} determined by the data~$\cD$; it should be thought of as taking the role 
that is played by $(\Qp/\Zp)^g \times \hat\mG_m^g$ in the classical theory. With these notations the main result obtained in 
sections~\refn{OrdBT+O} and~\refn{OrdBT+O*e} is the following.
\Bskip

\noindent
{\it Theorem. --- Let $k$ be an algebraically closed field of characteristic~$p$. Let $\ul{A}$ correspond to a $k$-valued 
point~$x$ of the PEL moduli space~$\cA_\cD$. Write $\ul{X} := \big(A[p^\infty],\iota,\lambda\big)$ and $\ul{Y} := 
\big(A[p],\iota,\lambda\big)$. Then the following are equivalent:
\item{{\rm (a)}} $\ul{A}$ is $\mu$-ordinary, which means that $\ul{X}$ is isogenous to~$\ul{X}^\ord(\cD)$;
\item{{\rm (b)}} $\ul{A}$ is $[p]$-ordinary, which means that $\ul{Y}$ is isomorphic to the $p$-kernel of~$\ul{X}^\ord(\cD)$;
\item{{\rm (c)}} $\ul{X}$ is isomorphic to~$\ul{X}^\ord(\cD)$.}
\Cskip

This theorem gives us a well-defined ordinary locus in~$\cA_\cD$ in characteristic~$p$. We remark that if $\ul{A}$ is ordinary then in 
general it is not true that the underlying abelian variety~$A$ is ordinary in the classical sense. In fact, given the PEL data~$\cD$ it 
is very easy to decide whether the underlying~$A$ is ordinary or not; see~\refn{OrdClassic}.

In~[\DFEO] we have proven a dimension formula for the Ekedahl-Oort strata on~$\cA$. In particular, we show that the $[p]$-ordinary 
stratum is the unique EO-stratum of maximal dimension. Combined with the theorem this gives a new proof of the main result of 
Wedhorn~[\WedhOrd]:
\Bskip

\noindent
{\it Corollary.\/ {\rm (Wedhorn)} --- The ordinary locus is Zariski dense in~$\cA_\cD$.}

\introsec{intrDefo}
{\it Deformation theory of ordinary objects.\/} If $x \in \cA_\cD(K)$ is an ordinary moduli point (in the ``new'' meaning of the term) 
then we have, at least over $k = \kbar$, a completely explicit description of the corresponding triple~$\ul{X}$. Note that $\ul{X}_k$ is 
independent of the ordinary point~$x$. This becomes particularly relevant when we study deformations of ordinary objects, as it implies 
that the structure of the formal deformation space $\mD := \Def(\ul{A}_x)$ is (geometrically) the same for all ordinary points~$x$. 

In order to explain our results on deformation theory, it is perhaps best to look at an example. So, let us suppose that we have a 
CM-field~$Z$ of degree~$2m$, with totally real subfield~$Z_0$. Suppose $p$ is a prime number that is totally inert in the 
extension $\mQ \subset Z_0$ and that splits in~$Z$. Consider an order $\cO \subset Z$ which is maximal at~$p$. If $q = p^m$ then 
$\cO \otimes \Zp \cong W(\mF_q) \times W(\mF_q)$, and the non-trivial automorphism of~$Z/Z_0$ acts by interchanging the two 
factors~$W(\mF_q)$. The triples $\ul{X} = \big(X=A[p^\infty],\iota,\lambda\big)$ arising in our moduli problem are of the 
following form. We have $X = X_1 \times X_1^D$, where $X_1$ is a BT with a given action of~$W(\mF_q)$, where $X_1^D$ is its 
Serre-dual with induced $W(\mF_q)$-action, and where the polarization is obtained from the duality between the two factors. Hence 
we are reduced to the study of BT with $W(\mF_q)$-action, without any polarization. We shall from now on use the letter~$\ul{X}$ 
for this somewhat simpler object.

At this point we can make the role of the ``CM-type'' more visible. Namely, suppose $\ul{X} = (X,\iota)$ is a BT with 
$W(\mF_q)$-action, over a field $k = \kbar$ with $\charact(k) = p$. Let $M$ be the Dieudonn\'e module of~$X$. Let $\cI$ be the set 
of embeddings $\mF_q \to k$. Note that this is a set of $m$~elements that comes equipped with a cyclic ordering: if $i \in \cI$ 
then we write $i+1 := \Frob_k \circ i$ for the successor of~$i$. Now $M$, being a module over $W(\mF_q) \otimes_\Zp W(k) = 
\prod_{i\in\cI} W(k)$, decomposes into character spaces: $M = \oplus_{i\in\cI} M_i$. Frobenius and Verschiebung restrict to maps 
$F\colon M_i \to M_{i+1}$ and $V\colon M_i \leftarrow M_{i+1}$. An easy lemma shows that $d := \rank_{W(k)}(M_i)$ is independent 
of~$i$; we call it the height of~$\ul{X}$. Next we define the {\it multiplication type\/} to be the function $\gf\colon \cI \to 
\{0,1,\ldots,d\}$ given by $\gf(i) := \dim_k\big(\Ker(F\colon M_i/pM_i \to M_{i+1}/pM_{i+1})\big)$.

In the formulation of the moduli problem~$\cA_\cD$, the invariants $d$ and~$\gf$ are fixed. The structure of the ordinary object 
$\ul{X}^\ord = \ul{X}^\ord(d,\gf)$ can be made fully explicit and depends only on the pair~$(d,\gf)$. We find a natural slope 
decomposition $\ul{X}^\ord = \ul{X}^{(1)} \times \cdots \times \ul{X}^{(r)}$. This is analogous to the decomposition of $(\Qp/\Zp)^g 
\times \hat\mG_m^g$ into its ind-\'etale part $(\Qp/\Zp)^g$ (slope~$0$) and its ind-multiplicative part $\hat\mG_m^g$ (slope~$1$). But, 
in contrast with the classical case, we can have any number $r \geq 1$ of slopes. (In fact, $r-1$ is the cardinality of the set 
$\gf(\cI) \cap [1,d-1]$.) The bigger~$r$, the more complex the structure that we find on the formal deformation space~$\mD$:
\Bskip

\noindent
{\it One slope.\/} If $r=1$ then we find that $\ul{X}^\ord = \ul{X}^{(1)}$ is rigid, meaning that its formal deformation functor 
is pro-represented by~$W(k)$. In particular, there is a {\it unique\/} lifting~$\ul{X}^\can$ of~$\ul{X}^\ord$ to a BT with 
$W(\mF_q)$-action over~$W(k)$.
\Bskip

\noindent
{\it Two slopes.\/} In this case we find a nice generalization of the Serre-Tate formal group structure on the formal 
deformation space. Suppose $\ul{X}^\ord = \ul{X}^{(1)} \times \ul{X}^{(2)}$. If $R$ is an artin local $W(k)$-algebra with residue 
field~$k$ then the factors~$\ul{X}^{(\nu)}$ each admit a unique lifting~$\ulcX^{(\nu)}$ over~$R$. The first thing we show is that 
every deformation of~$\ul{X}$ is an extension of~$\ulcX^{(1)}$ by~$\ulcX^{(2)}$. As in the classical case, it follows that the 
formal deformation space~$\mD$ has the structure of a formal group over~$W(k)$. The $W(\mF_q)$-action on the two 
factors~$\ul{X}^{(\nu)}$ induces an action on the formal group~$\mD$.

Now comes the best part. Each of the two factors~$\ul{X}^{(\nu)}$ is isomorphic to a product of a number of copies of an ordinary 
object of height~$1$. So, we have integers $d^1$ and~$d^2$ with $d^1 + d^2 = d$, and we have functions $\gg^1$, $\gg^2 \colon \cI \to 
\{0,1\}$, such that $\ul{X}^{(\nu)}$ is isomorphic to the product of $d^\nu$~copies of the object~$\ul{X}^\ord(1,\gg^\nu)$. 
Moreover, the slope decomposition is such that $\gg^1(i) \leq \gg^2(i)$ for all $i \in \cI$. Now we form a new multiplication type 
$\gf^\prime\colon \cI \to \{0,1\}$ by ``subtracting'' $\gg^1$ from~$\gg^2$: let $\gf^\prime(i) = 0$ if $\gg^1(i) = \gg^2(i)$ and 
$\gf^\prime(i) = 1$ if $\gg^1(i) = 0$ and $\gg^2(i)=1$. The associated ordinary object $\ul{X}^\ord(1,\gf^\prime)$ is isoclinic 
(one slope), so by the above it has a canonical lifting $\ul{X}^\can(1,\gf^\prime)$ over~$W(k)$. The result that we find is as follows.
\Bskip

\noindent
{\it Theorem. --- With notation as above, the formal deformation space~$\mD = \Def(\ul{X}^\ord)$ has a natural structure of a BT 
with $W(\mF_q)$-action over~$W(k)$, and we have
$$
\mD \cong \ul{X}^\can(1,\gf^\prime)^{d^1d^2}\, .
$$\vskip-\lastskip\medskip}

\noindent
It is nice to compare this with the classical ordinary case: If the underlying abelian variety~$A$ is ordinary in the classical 
sense then we have $\gg^1(i)=0$ and $\gg^2(i) =1$ for all~$i\in\cI$, hence $\gf^\prime(i) = 1$ for all~$i$. This is precisely the 
case where $\mD$ is a formal torus.
\Bskip

\noindent
{\it More than two slopes.\/} If $\ul{X}^\ord = \ul{X}^{(1)} \times \cdots \times \ul{X}^{(r)}$ with $r > 2$ then we no longer 
find the structure of a formal group on the formal deformation space~$\mD$. We introduce a new notion, called an {\it 
$r$-cascade\/}, and we show that $\mD$ naturally admits such a structure. For $1 \leq a < b \leq r$, let us introduce the notation 
$\ul{X}^{(a,b)} := \prod_{\nu = a}^b \ul{X}^{(\nu)}$. A key observation is that for any deformation of~$\ul{X}$ we can lift the 
slope filtration. More precisely: if $\ulcX$ is a deformation of~$\ul{X}$ then there is a unique filtration by sub-objects
$$
0 \subset \ulcX^{(r,r)} \subset \ulcX^{(r-1,r)} \subset \cdots \subset \ulcX^{(1,r)} = \ulcX
$$
such that the special fibre of~$\ulcX^{(a,r)}$ is~$\ul{X}^{(a,r)}$. As a consequence we find that we can arrange the formal 
deformation spaces $\mD^{(a,b)} := \Def\big(\ul{X}^{(a,b)}\big)$ in a ``tower'':
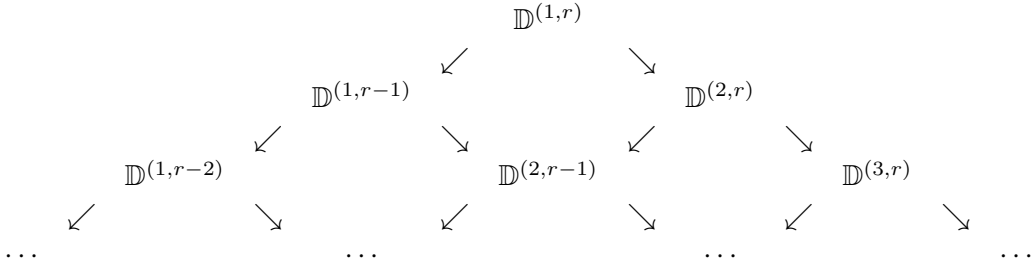
$$
\matrix{
&&&&&&\mD^{(1,r)}&&&&&&\cr
&&&&&\swarrow&&\searrow&&&&&\cr
&&&&\mD^{(1,r-1)}&&&&\mD^{(2,r)}&&&&\cr
&&&\swarrow&&\searrow&&\swarrow&&\searrow&&&\cr
&&\mD^{(1,r-2)}&&&&\mD^{(2,r-1)}&&&&\mD^{(3,r)}&&\cr
&\swarrow&&\searrow&&\swarrow&&\searrow&&\swarrow&&\searrow&\cr
\cdots&&&&\cdots&&&&\cdots&&&&\cdots\cr}
$$
Saying that $\mD^{(1,r)} = \Def(\ul{X})$ has the structure of an $r$-cascade essentially means that in each diamond
$$
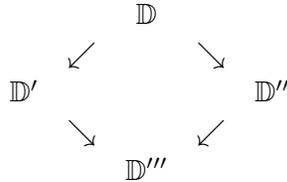
\matrix{
&&\mD&&\cr
&\swarrow&&\searrow&\cr
\mD^\prime&&&&\mD^\pprime\cr
&\searrow&&\swarrow&\cr
&&\mD^\ppprime&&\cr}
$$
in the tower, viewing $\mD^\ppprime$ as the formal base scheme, $\mD^\prime$ and $\mD^\pprime$ have the structure of a relative 
formal group, and that $\mD$ has the structure of a biextension over $\mD^\prime \times \mD^\pprime$. The structure groups 
occuring in these biextensions can all be made explicit in terms of ``subtraction of multiplication types'', as in the case of two 
slopes.

Even though an $r$-cascade is (for $r>2$) a more complicated object than a formal group, it has a number of ``group-like'' 
features that play a role in Serre-Tate theory. Thus, for instance, $\mD$~has an origin, corresponding to a canonical lifting 
of~$\ul{A}$. Also we have the notion of a torsion point, and we show, for $k=\Fpbar$, that the torsion points of~$\mD$ correspond 
precisely to the CM-liftings of~$\ul{A}$, as in the classical theory.

\introsec{intrACD} 
As is well-known, PEL moduli problems break up in three cases, labelled A, C and~D. The example that we have taken 
in~\refn{intrDefo} is of type~A, which tends to produce the most interesting new phenomena. In Case~C we essentially only find 
structures that are already covered by classical Serre-Tate theory. Case~D, finally, is technically the most problematic. Despite 
the extra work, we have included type~D throughout the paper. This requires that, in addition to the invariants $d$ and~$\gf$, we 
keep track of a further discrete invariant,~$\delta$.

\introsec{intrCong}
In the final section of this paper we discuss an application of our theory. In~\refn{EOonPEL} we introduce the moduli 
spaces~$\cA_\cD$ and we discuss how they decompose as a union of Shimura varieties. We discuss the possible values that the 
discrete invariants $(d,\gf,\delta)$ can take on these components, and we give some results on their fields of definition.

In section \refn{CongRel} we study congruence relations. The conjecture, as formulated by Blasius and Rogawski in~[\BlaRog], is 
that the Frobenius correspondence~$\Phi$ on (the Shimura components of)~$\cA_\cD$ in characteristic~$p$ satisfies a certain 
polynomial relation of which the coefficients are Hecke correspondences. This Hecke polynomial $H_{(\cG,\cX)}$ is defined in a 
purely group-theoretic way, starting from a Shimura datum~$(\cG,\cX)$.

Most of the material in section~\refn{CongRel} closely follows Wedhorn's paper~[\WedhCong]. The approach taken here is the one of Chai 
and Faltings in Chapter~\Romno 7 of their book~[\FaCh]. The main result of this section is that the desired relation $H_{(\cG,\cX)}(\Phi) = 0$ holds over the ordinary locus. We refer to the body of the text for a precise statement. 

\introsec{Ackn}
{\it Acknowledgements.\/} I thank Frans Oort and Torsten Wedhorn for their interest in my work and for stimulating discussions. The 
research for this paper was made possible by a Fellowship of the Royal Netherlands Academy of Arts and Sciences (KNAW). During my work 
on this paper I have been affiliated to the University of Utrecht (until June 2001) and the University of Amsterdam (from July 2001); I 
thank both institutions for their support.

\introsec{Notat}
{\it Notation.\/} We typically use the letter~$X$ for Barsotti-Tate groups and $Y$ for BT$_1$. The Diedonn\'e module of~$X$ (resp.~$Y$) 
is called~$M$ (resp.~$N$). In discussions about BT$_n$ for arbitrary $n \in \mN \cup \{\infty\}$ we use the letter~$X$. For deformations 
we often use $\cX$ and~$\cY$. For abelian varieties we use the letter~$A$. Underlined letters represent objects equipped with an action 
of a given ring and possibly also a polarization.     
\section{Ordinary Barsotti-Tate $\cO$-modules}{Ordi}

\ssection{BT$_1$ with $\cO$-structure}{BT1OStr}

\sssection{BTdef} 
We fix a prime number~$p$. For the definition of a Barsotti-Tate group (= $p$-divisible group) and a truncated Barsotti-Tate group we refer to Illusie~[\Illus]. We abbreviate ``Barsotti-Tate group'' to BT$_\infty$ or simply BT, and ``truncated Barsotti-Tate group of level~$n$'' to BT$_n$.

Let $\cO$ be a $\Zp$-algebra. Let $n \in \mN \cup \{\infty\}$. By a BT$_n$ with $\cO$-structure over a basis~$S$ we mean a pair $\ul{X} = (X,\iota)$ where $X$ is a BT$_n$ over~$S$ and $\iota\colon \cO \to \End_S(X)$ is a homomorphism of $\Zp$-algebras. (An alternative name would be ``BT$_n$ $\cO$-module'', which is less satisfactory for typographical reasons.)

We shall study BT$_n$ using contravariant Dieudonn\'e theory as in Fontaine~[\Font].

\sssection{pairs} 
Let $B$ be a finite dimensional semi-simple $\Fp$-algebra. Let $k$ be an algebraically closed field of characteristic~$p$. The first problem studied in~[\GSAS] is the classification of BT$_1$ with $B$-structure over~$k$. This generalizes the work of Kraft~[\Kraft], who classified group schemes killed by~$p$ without additional structure. We shall briefly review our results. 

Write $\kappa$ for the center of~$B$. Then $\kappa$ is a product of finite fields, say $\kappa = \kappa_1 \times \cdots \times \kappa_\nu$. Let~$\cI = \cI_1 \cup \cdots \cup\cI_\nu$ be the set of homomorphisms $\kappa \rightarrow k$.

Consider pairs $(N,L)$ consisting of a finitely generated $B \otimes_\Fp k$-module~$N$ and a submodule $L \subset N$. Note that the simple factors of $B \otimes_\Fp k$ are indexed by~$\cI$, so we get canonical decompositions $N = \oplus_{i\in\cI} N_i$ and $L = \oplus_{i\in\cI} L_i$. Define two functions $d$, $\gf \colon \cI \to \mZ_{\geq 0}$ by $d(i) = \len(N_i)$ and $\gf(i) = \len(L_i)$, takings lengths as $B \otimes_{\Fp} k$-modules. The pair $(d,\gf)$ determines the pair $(N,L)$ up to isomorphism.

To the pair $(N,L)$ we associate an algebraic group~$G$ over~$k$ and a conjugacy class~$\mX$ of parabolic subgroups of~$G$. First we define
$$
G := \GL_{B \otimes_\Fp k}(N) \cong \dirprod_{i\in\cI} \GL_{d(i),k}\, .
$$
Then the stabilizer $P := \Stab(L)$ is a parabolic subgroup of~$G$, and we define
$$
\mX := \hbox{conjugacy class of parabolic subgroups of~$G$ containing $P$.}
$$

\sssection{BT1pair}
Let $\ul{Y} = (Y,\iota)$ be a BT$_1$ with $B$-structure over~$k$. Write~$N$ for the Dieudonn\'e module of~$Y$ and let $L := \Ker(F) \subset N$. Let $(d,\gf)$ be the corresponding pair of functions. It can be shown (see~[\GSAS], 4.3) that the function~$d$ is constant on each of the subsets $\cI_n \subset \cI$. We refer to $(d,\gf)$ as the {\it type\/} of~$(Y,\iota)$.

\sssection{wdefGE}
Fix a type~$(d,\gf)$ with $d$ constant on each subset $\cI_n \subset \cI$. Fix a pair of $B \otimes_\Fp k$-modules $L_0 \subset N_0$ of type~$(d,\gf)$. Let $(G,\mX)$ be the associated algebraic group and conjugacy class of parabolic subgroups. Let $W_G$ be the Weyl group of~$G$, and let $W_\mX \subset W_G$ be the subgroup corresponding to~$\mX$.

To a pair $\ul{Y} = (Y,\iota)$ of type $(d,\gf)$ we associate an element $w(\ul{Y}) \in W_\mX \backslash W_G$. This is done as follows. Write~$(N,F,V)$ for the Dieudonn\'e module of~$Y$. As $Y$ is a BT$_1$ we have $\Ker(F) = \Image(V)$ and $\Image(F) = \Ker(V)$. Using this, one can show that there exists a filtration
$$
\calC_\gdot\colon
\quad (0) = \calC_0 \subset \calC_1 \subset \calC_2 \subset \cdots \subset \calC_r = N
$$
that is the coarsest filtration with the properties that
\item{(\romno1)}for every $j$ there exists an index $f(j) \in \{0,1,\ldots,r\}$ with $F(\calC_j) = 
\calC_{f(j)}$;
\item{(\romno2)}for every $j$ there exists an index $v(j) \in \{0,1,\ldots,r\}$ with 
$V^{-1}(\calC_j) = \calC_{v(j)}$.

\noindent
We refer to this filtration as the {\it canonical filtration\/} of~$N$.

Set $L := \Ker(F) \subset N$. Choose an isomorphism $\xi\colon N \isomarrow N_0$ that restricts to~$L \isomarrow L_0$. This allows us to view~$\calC_\gdot$ as a filtration of~$N_0$. Choose any refinement~$\cF_\gdot$ of~$\calC_\gdot$ to a complete flag. The relative position of~$L_0$ and $\cF_\gdot$ is given by an element $w(L_0,\cF_\gdot) \in W_\mX \backslash W_G$. It can be shown that this element is independent of the choice of~$\xi$ and the refinement~$\cF_\gdot$; see~[\GSAS], especially~4.6 for details. Now define $w(\ul{Y}) := w(L_0,\cF_\gdot)$.

With these notations, the first main result of~[\GSAS] can be stated as follows.

\sssection{GEThm}
{\it Theorem. --- Assume that $k = \kbar$. The map $\ul{Y} \mapsto w(\ul{Y})$ gives a bijection
$$
\left\{
\vcenter{
\setbox0=\hbox{{\rm isomorphism classes of}}
\copy0
\hbox to \wd0{{\rm \hfil $\ul{Y}$ of type $(d,\gf)$\hfil}}}
\right\}
\longisomarrow
W_\mX \backslash W_G\, .
$$\vskip-\lastskip\smallskip}

\sssection{modpRefl}
We retain the notation of~\refn{pairs}. Note that $\Aut(k)$ naturally acts on the set~$\mZ^\cI$.

Fix a type $(d,\gf)$. Assume that the function $d\colon \cI \to \mZ_{\geq 0}$ is constant on each of the subsets $\cI_n \subset \cI$; this is equivalent to the condition that $d$ is invariant under~$\Aut(k)$. Next consider $\Stab(\gf) := \{\alpha \in \Aut(k) \mid {}^\alpha\gf = \gf\}$. We define $E(d,\gf) \subset k$ to be the fixed field of~$\Stab(\gf)$. For instance, if $B$ is a simple algebra, let $m_0$ be the smallest positive integer with the property that $\gf(i+m_0) = \gf(i)$ for all $i\in\cI$; then $E(d,\gf) \subset k$ is the subfield with $p^{m_0}$ elements. We refer to $E(d,\gf)$ as the ``mod~$p$ reflex field''; see Remark~\refn{kappa(v)} for an explanation of this terminology.

Let $K \subset k$ be a perfect subfield. Let $\ul{Y}$ be a BT$_1$ with $\kappa$-structure over~$K$ of type $(d,\gf)$. We claim that the existence of such an object implies that $E(d,\gf) \subset K$. To see this we may replace~$k$ by the separable closure of~$K$. If $\alpha \in \Aut(k)$ then ${}^\alpha \ul{Y}$ has type $(d,{}^\alpha\gf)$. The assumption that $\ul{Y}$ is defined over~$K$ therefore implies that $\Gal(k/K) \subset \Stab(\gf)$.

\sssection{MoritaNPol}
The study of BT$_1$ with $B$-structure easily reduces to the case that $B = \kappa$ is a finite field. Indeed, as the Brauer group of a finite field is trivial we have $B \cong M_{r_1}(\kappa_1) \times \cdots \times M_{r_l}(\kappa_l)$, where the $\kappa_n$ are finite fields, $\charact(\kappa_n) = p$. Fixing such an isomorphism, every BT$_1$ with $B$-structure decomposes as $\ul{Y} = (\ul{Y}_1)^{r_1} \times \cdots \times (\ul{Y}_l)^{r_l}$, where $\ul{Y}_n$ is a BT$_1$ with $\kappa_n$-structure.

If $B = \kappa$ is a finite field, $\cI$ is simply the set of embeddings $\kappa \to k$. This set comes equipped with a natural cyclic ordering: if $i \in \cI$ we write $i+1 := \Frob_k \circ i$ for its successor. The type $(d,\gf)$ of a BT$_1$ with $\kappa$-structure consists of a non-negative integer~$d$ and a function $\gf\colon \cI \to \{0,\ldots,d\}$. The integer~$d$ is also referred to as the {\it height\/} of the truncated Barsotti-Tate $\kappa$-module. (The underlying BT$_1$ without additional structure has height $d \cdot [\kappa:\Fp]$.)

\ssection{The $[p]$-ordinary type}{[p]ordi}

\sssection{OrdNPolSit}
{\it Situation.} --- Let $k$ be an algebraically closed field of characteristic $p>0$. Let $\sigma$ be the Frobenius automorphism of~$W(k)$. Let $\kappa$ be a field of $p^m$ elements, and write $\cO = W(\kappa)$. Recall that $\cI := \Hom(\kappa,k) = \Hom\big(\cO,W(k)\big)$.

Let $\ul{X}$ be a BT with $\cO$-structure over~$k$. Write $\ul{Y} := \ul{X}[p]$, which is a BT$_1$ with $\kappa$-structure. Let $(d,\gf)$ be its type. Let $(G,\mX)$ be as in~\refn{wdefGE}. We fix a maximal torus and Borel subgroup in~$G$; this gives us a set of generators for the Weyl group~$W_G$.

\sssection{[p]ordDef}
{\it Definition.\/} --- Let $w^\ord \in W_\mX\backslash W_G$ be the class of the longest element of~$W_G$. We say that $\ul{Y}$ is {\it $[p]$-ordinary\/}, and also that $\ul{X}$ is $[p]$-ordinary, if $w(\ul{Y}) = w^\ord$.
\Cskip

By Thm.~2.1.2 of~[\DFEO], $\ul{Y}$ is $[p]$-ordinary if and only if $\Aut(\ul{Y})$ is finite.

\sssection{XordNPol}
We define a BT with $\cO$-structure $\ul{X}^\ord = \ul{X}^\ord(d,\gf)$ that is $[p]$-ordinary. We refer to $\ul{X}^\ord$ as the {\it standard ordinary\/} BT with $\cO$-structure of type~$(d,\gf)$.

Let $M$ be the free $W(k)$-module with basis $e_{i,j}$ for $i\in\cI$ and $j\in\{1,\ldots,d\}$. Write $M_i := \sum_{j=1}^d W(k)\cdot e_{i,j}$, and let $a \in \cO$ act on~$M_i$ as multiplication by~$i(a) \in W(k)$. Next define Frobenius and Verschiebung on base vectors by
$$
F(e_{i,j}) = \cases{e_{i+1,j}&if $j \leq d-\gf(i)$;\cr
p \cdot e_{i+1,j}&if $j>d-\gf(i)$;}
\qquad
V(e_{i+1,j}) = \cases{p \cdot e_{i,j}&if $j \leq d-\gf(i)$;\cr
e_{i,j}&if $j>d-\gf(i)$.}
$$
These data define a Dieudonn\'e module $\ul{M}^\ord = \ul{M}^\ord(d,\gf)$ with $\cO$-structure, and we define $\ul{X}^\ord$ to be the corresponding BT with $\cO$-structure. It is easily verified that $\ul{X}^\ord$ is indeed $[p]$-ordinary; cf.\ [\GSAS],~4.9.

We write $\ul{Y}^\ord = \ul{Y}^\ord(d,\gf) := \ul{X}^\ord[p]$. Its Dieudonn\'e module is $\ul{N}^\ord = \ul{M}^\ord/p \ul{M}^\ord$.

\sssection{EtMult}
{\it Remark.\/} --- Take $d=1$ and $\gf(i) = 0$ for all~$i$. We write $\ul{X}_\et$ for the corresponding standard ordinary BT. If $K$ is the fraction field of~$\cO$ then $\ul{X}_\et$ is none other than the ind-\'etale $p$-divisible group~$K/\cO$ with its natural $\cO$-structure. The underlying BT without $\cO$-structure is isomorphic to $(\Qp/\Zp)^m$. Let $\ul{X}_\mult$ be the standard ordinary object corresponding to $d=1$ and $\gf(i)=1$ for all~$i$; it is the Serre dual of~$\ul{X}_\et$. The underlying BT in this case is~$(\hat{\mG}_{{\rm m}})^m$.

In general, we do not know a description of~$\ul{X}^\ord$ as a functor, other than using the ``inverse Dieudonn\'e functor'' $M \mapsto X$ as in Fontaine~[\Font], \Romno 3,~1.3. 

\sssection{ajDef}
We fix a type $(d,\gf)$. Let $r-1$ be the number of values in the interval $[1,d-1]$ that occur as $\gf(i)$ for some $i\in\cI$. Define integers $0\leq a_1 \leq a_2 \leq \cdots \leq a_d$ by
$$
a_j := \# \big\{i \in \cI \bigm| \gf(i) > d-j\big\}\, .
$$ 
Let $0 \leq \lambda_1 < \lambda_2 < \cdots < \lambda_r$ be the integers occurring as $a_j$ for some~$j$; note that there are precisely~$r$ of them. Let $d^\nu := \#\{j\mid a_j =\lambda_\nu\}$ and define functions
$$
\gf^\nu\colon \cI \to \{0,d^\nu\}
\quad\hbox{by}\quad
\gf^\nu(i) = \cases{\strut 0&if $\gf(i) < \sum_{j=\nu}^r d^j$;\cr
d^\nu&if $\gf(i) \geq \sum_{j=\nu}^r d^{j\strut}$.\cr}
$$
With the obvious meaning of the notation we have $(d,\gf) = \sum_{\nu=1}^r (d^\nu,\gf^\nu)$.

Finally we define $\Ord(d,\gf)$ to be the polygon with slopes~$a_j$ ($j=1,\ldots,d$).

\sssection{MordDec}
Consider the Dieudonn\'e module $\ul{N}^\ord = \ul{N}^\ord(d,\gf)$ as in~\refn{XordNPol}. For each $j \in \{1,\ldots,d\}$ the subspace $\sum_{i\in\cI} k \cdot e_{i,j}$ of~$N$ is a Dieudonn\'e submodule, stable under the action of~$\kappa$. We find that $\ul{N}^\ord$ is a direct sum of $d$ objects of height~$1$. Grouping together isomorphic layers we get a decomposition into isotypic components
$$
\ul{N}^\ord(d,\gf) = \ul{N}^\ord(d^1,\gf^1) \oplus \cdots \oplus \ul{N}^\ord(d^r,\gf^r)\, .
$$
For $\ul{Y}^\ord$ this gives a decomposition $\ul{Y}^\ord(d,\gf) = \prod_{\nu = 1}^r \ul{Y}^\ord(d^\nu,\gf^\nu)$. Similarly we have a {\it slope decomposition\/} $\ul{X}^\ord(d,\gf) = \prod_{\nu =1}^r \ul{X}^\ord(d^\nu,\gf^\nu)$.

\sssection{EndOrdLem}
{\it Lemma. --- Notation as in\/ {\rm \refn{MordDec}}. 

{\rm (\romno1)} We have
$$
\eqalign{
\End\big(\ul{Y}^\ord\big) &= \End\big(\ul{Y}^\ord(d^1,\gf^1)\big) \times \cdots \times \End\big(\ul{Y}^\ord(d^r,\gf^r)\big)\cr
&\cong M_{d^1}(\kappa) \times \cdots \times M_{d^r}(\kappa)\, .\cr}
$$

{\rm (\romno2)} If $\ul{Y}$ is a $[p]$-ordinary BT$_1$ with $\kappa$-structure of type $(d,\gf)$ then we have a canonical decomposition $\ul{Y} = \ul{Y}^{(1)} \times \cdots \times \ul{Y}^{(r)}$ such that $\ul{Y}^{(\nu)} \cong \ul{Y}^\ord(d^\nu,\gf^\nu)$.}
\Dskip

\Proof~By definition, $\ul{Y}$ being $[p]$-ordinary means that it is isomorphic to~$\ul{Y}^\ord$. As such an isomorphism is unique up to an element of $\Aut(\ul{Y}^\ord)$ we see that (\romno2) is an immediate consequence of~(\romno1). The proof of~(\romno1) is an easy exercise, using arguments as in~[\DFEO], section~2.2. \QED
\Cskip

As mentioned in~\refn{MordDec}, $\ul{Y}$ is isomorphic to a product of $d$ objects of height~$1$; note however that this finer decomposition is {\it not\/} canonical, unless $r = d$.

\ssection{Ordinary BT with $\cO$-structure}{OrdBT+O}

\sssection{BTWaction}
Situation as in~\refn{OrdNPolSit}. Let $\ul{M} = (M,F,V,\iota)$ be the Dieudonn\'e module of~$\ul{X}$. We have a natural decomposition into character spaces $M = \dirsum_{i\in\cI} M_i$. Frobenius and Verschiebung restrict to $\sigma$-linear maps $F_i\colon M_i \to M_{i+1}$ and $\sigma^{-1}$-linear maps $V_i\colon M_i \leftarrow M_{i+1}$. 

Let $\# \kappa = p^m$, so that $\cI$ has $m$ elements. For $i\in\cI$, define the $\sigma^m$-linear endomorphism $\Phi_i$ of $M_i$ by
$$
\Phi_i = \big(M_i \mapright{F_i} M_{i+1} \mapright{F_{i+1}} \cdots \longrightarrow M_{i+m-1} \mapright{F_{i+m-1}} 
M_{i+m} = M_i\big)\, .
$$
Then $(M_i,\Phi_i)$ is a $\sigma^m$-$F$-crystal over~$k$. By construction, $F_i$ induces an isogeny from 
$\sigma^\ast(M_i,\Phi_i)$ to $(M_{i+1},\Phi_{i+1})$. In particular, the Newton polygon of $(M_i,\Phi_i)$ is independent of $i\in\cI$. We refer to this as the Newton polygon of~$\ul{X}$. (This must not be confused with the Newton polygon of~$X$.) By contrast, easy examples show that in general the Hodge polygon of $(M_i,\Phi_i)$ does depend on~$i$.

\sssection{ulHIsog}
{\it Proposition. --- Let $\ul{X}$ and $\ul{X}^\prime$ be BT with $\cO$-structure over~$k$. Then $\ul{X}$ and $\ul{X}^\prime$ are isogenous if and only if their Newton polygons are the same.\/}
\Dskip

\Proof This is an application of the theory of isocrystals with additional structure; see Kottwitz~[\KottIsoc] and Rapoport and Richartz~[\RR]. We sketch the argument.

Let $K$ be the fraction field of~$\cO$. Write $Q$ for the fraction field of~$W(k)$. Set $d := \height(\ul{X})$, and let $\Gamma := \Res_{K/\Qp} \GL_d$. The isogeny class of $\ul{X}$ is classified by $\ul{M} \otimes_{\Zp} \Qp$, which is an isocrystal with $\Gamma$-structure. This isocrystal, in turn, is classified by a $\sigma$-conjugacy class in $\Gamma(Q)$. The set $B(\Gamma)$ of such $\sigma$-conjugacy classes is studied by means of a so-called Newton map $B(\Gamma) \to \cN(\Gamma)$; in our situation this is simply the map which associates to $\ul{M} \otimes_{\Zp} \Qp$ the Newton polygon of $(M_i,\Phi_i)$. Let $\bar b \in B(\Gamma)$ be the $\sigma$-conjugacy class of $b \in \Gamma(Q)$. To the element $b$ one associates an algebraic group $J_b$, and it is shown (see [\RR], Prop.~1.17) that the fibre of the Newton map that contains $\bar b$ is in bijection with $H^1(\Qp,J_b)$. Further, $J_b$ is an inner form of a Levi subgroup of~$\Gamma$. In our situation this means that $J_b$ is an inner form of a product of general linear groups. Hence $H^1(\Qp,J_b) = 0$ and the Newton map is injective. \QED

\sssection{muordDef}
{\it Definition.\/} --- Let $\ul{X}$ be a BT with $\cO$-structure over~$k$, of type $(d,\gf)$. We say that $\ul{X}$ is {\it $\mu$-ordinary\/} if its Newton polygon equals the polygon $\Ord(d,\gf)$ defined in~\refn{ajDef}.
\Cskip

The terminology ``$\mu$-ordinary'' follows Wedhorn's paper [\WedhOrd].

\sssection{NewtHodLem}
{\it Lemma. --- Let $\ul{X}$ be as in\/ {\rm \refn{BTWaction}}. Then its Newton polygon is on or above the polygon $\Ord(d,\gf)$ defined in\/ {\rm \refn{ajDef}}, and the two polygons have the same end point.}
\Dskip

\Proof~This is [\RR], Thm.~4.2, (\romno2), taking into account loc.\ cit., Prop.~2.4,~(\romno4). \QED

\sssection{b1a1Lem}
{\it Lemma. --- Let $\ul{X}$ be as in\/ {\rm \refn{BTWaction}}. Define $a_j$ as in\/ {\rm \refn{ajDef}}. Fix $i \in \cI$, and let $0 \leq b^{(i)}_1 \leq b^{(i)}_2 \leq \cdots b^{(i)}_d$ be the Hodge slopes of $(M_i,\Phi_i)$. Then $b^{(i)}_1 \geq a_1$.}
\Dskip

\Proof~By [\Katz], 1.2.1 we have to show that $\Phi_i(M_i) \subset p^{a_1} \cdot M_i$. By definition, $a_1$ is the number of indices $\nu \in \cI$ such that $\gf(\nu) = d$. But $\gf(\nu) = d$ just means, writing $N := M/pM$, that the Frobenius $F_\nu \colon N_\nu \to N_{\nu+1}$ is zero, which is equivalent to saying that $F_\nu\colon M_\nu \to M_{\nu+1}$ lands inside $p \cdot M_{\nu+1}$. As $\Phi_i$ is obtained as a composition of all $F_{\nu}$, the lemma follows. \QED

\sssection{eqLem}
{\it Lemma. --- Let $R$ be a complete local ring with $\charact(R/\gm) = p$. Let $\rho$ be an automorphism of~$R$. Let $C_1$, $C_2$, $C_3$ and $C_4$ be matrices with coefficients in~$R$, of sizes $r \times s$, $s \times s$, $s \times r$ and $r\times r$, respectively. Assume that $C_2$ is invertible. Then the matrix equation
$$
C_1 + XC_2 + pC_3\cdot {}^\rho X + pXC_4\cdot {}^\rho X= 0
$$
has a solution for $X \in M_{r \times s}(R)$.}
\Dskip

\Proof~The idea is simply that a solution has the form $X = -C_1C_2^{-1} + p\cdot X^\prime$, where $X^\prime$ satisfies an equation
$$
C_1^\prime + X^\prime C_2^\prime + pC_3^\prime\cdot {}^\rho X^\prime + pX^\prime C_4^\prime\cdot {}^\rho X^\prime = 0\, .
$$
Writing $\Gamma := C_1C_2^{-1}$, the coefficients of the new equation are given by 
$$
C_1^\prime = C_3{}^\rho\Gamma + \Gamma C_4{}^\rho \Gamma\, ,\quad
C_2^\prime = C_2-pC_4{}^\rho\Gamma\, ,\quad
C_3^\prime = C_3-\Gamma C_4\, ,\quad\hbox{and}\quad
C_4^\prime = pC_4\, .
$$
Iterating this we obtain a power series development for $X$, which converges because $R$ is $p$-adically complete. \QED

\sssection{ord=ord}
{\it Theorem. --- Let $\ul{X}$ be a BT with $\cO$-structure over~$k$. Let $(d,\gf)$ be its type. Then the following are equivalent:
\item{\rm (a)} $\ul{X}$ is $\mu$-ordinary;
\item{\rm (b)} $\ul{X}$ is $[p]$-ordinary;
\item{\rm (c)} $\ul{X} \cong \ul{X}^\ord(d,\gf)$.}
\Dskip

\Proof~Recall the definition of~$r$ (=~the number of slopes) in~\refn{ajDef}. We first prove the theorem under the assumption that it is true for~$r=1$.

We have (c) $\Rightarrow$ (b) by definition of~$\ul{X}^\ord(d,\gf)$. Now assume that $\ul{X}$ is $\mu$-ordinary and that $r>1$. The first assumption means, by definition, that the Newton polygon of~$\ul{X}$ equals $\Ord(d,\gf)$. Let $q$ be index such that $a_1 = a_2 = \cdots = a_q < a_{q+1}$. (Note that not all $a_j$ are equal, as $r>1$.) By Mazur's basic slope estimate (see [\Katz], Thm.~1.4.1), $\Ord(d,\gf)$ lies on/above the Hodge polygon of $(M_i,\Phi_i)$, for every $i \in \cI$. Combining this with Lemma~\refn{b1a1Lem} we find that the first~$q$ slopes of these polygons are equal, i.e., $a_1 = b^{(i)}_1 = \cdots = b^{(i)}_q$. Katz~[\Katz], Thm.~1.6.1, then tells us that $(M_i,\Phi_i)$ decomposes:
$$
(M_i,\Phi_i) =  (M_i^\prime,\Phi_i^\prime) \oplus (M_i^\pprime,\Phi_i^\pprime)
$$
in such a way that 
$$
\displaylines{
\hbox{Newton slopes}\ (M_i^\prime,\Phi_i^\prime) = \hbox{Hodge slopes}\ (M_i^\prime,\Phi_i^\prime) = (a_1,a_2,\ldots,a_q)\, ;\cr
\hbox{Newton slopes}\ (M_i^\pprime,\Phi_i^\pprime) = (a_{q+1}, \ldots,a_d)\, ;\quad \hbox{Hodge slopes}\ 
(M_i^\pprime,\Phi_i^\pprime) = (b^{(i)}_{q+1},\ldots,b^{(i)}_d)\, .\cr}
$$

By construction, $M_i^\prime$ and $M_{i+1}^\pprime$ have no Newton slopes in common. Therefore, $F_i$ maps $M_i^\prime$ into $M_{i+1}^\prime$ and $M_i^\pprime$ into $M_{i+1}^\pprime$. This means that $\ul{X}$ decomposes as $\ul{X} = \ul{X}^\prime \times \ul{X}^\pprime$, such that the Dieudonn\'e module of $\ul{X}^\prime$ (resp.\ $\ul{X}^\pprime$) is $\oplus M_i^\prime$ (resp.\ $\oplus M_i^\pprime$). It easily follows from Lemma~\refn{NewtHodLem} that $\ul{X}^\prime$ and $\ul{X}^\pprime$ are again $\mu$-ordinary. By induction on the number~$r$ we may then assume that each is isomorphic to a standard ordinary~BT with $\cO$-structure. (Here we use the assumption that the theorem is true for $r=1$.)  This readily implies that $\ul{X}$ is isomorphic to the standard ordinary~BT of type $(d,\gf)$ and that $\ul{X}$ is $[p]$-ordinary.

Next assume that $\ul{X}$ is $[p]$-ordinary. As usual we write $\ul{N} = \ul{M}/p\ul{M}$ for the Dieudonn\'e module of~$\ul{Y}$. Our assumption means that $\ul{N} \cong \ul{N}^\ord(d,\gf)$. Define~$q$ as above. By~\refn{EndOrdLem} we have a canonical decomposition
$$
\ul{N} = \ul{N}^{(1)} \oplus \cdots \oplus \ul{N}^{(r)}\, .
$$
Let $\ul{N}^\prime := \ul{N}^{(1)}$ and $\ul{N}^\pprime := \oplus_{\nu=2}^r \ul{N}^{(\nu)}$. The decompositions $N_i = N_i^\prime \oplus N_i^\pprime$ are such that, for all $i \in \cI$, 
$$
\hbox{{\rm either $F_i\colon N_i \to N_{i+1}$ is injective on $N_i^\prime$,\quad or $F_i = 0$ on all of~$N_i$.}}\leqno(\refn{ord=ord}.1)
$$

We want to show that there is a sub-crystal $M_i^\prime \subset M_i$ that reduces to $N_i^\prime$ modulo~$p$. As we have seen in Lemma~\refn{b1a1Lem}, $\Phi_i$ is divisible by~$p^{a_1}$; hence we can define a $\sigma^m$-linear endomorphism~$\Psi_i$ of~$M_i$ by $\Psi_i := p^{-a_1}\cdot\Phi_i$. (Here $m = \#\cI$.) If $y \in M_i$, write $\overline y$ for its class in~$N_i$. We need two properties of~$\Psi_i$. Firstly,
$$
\hbox{{\rm if $\overline{y} \in N_i^\pprime$ then $\Psi_i(y) \in p\cdot M_i$.}}\leqno(\refn{ord=ord}.2)
$$
This is shown by the same argument as in~\refn{b1a1Lem}; note that the number of elements $i\in\cI$ for which $F\colon N_i^\pprime \to N_{i+1}^\pprime$ is zero is $> a_1$, by construction of~$N^\pprime$. The other property we need is that 
$$
\hbox{{\rm if $0 \neq \overline{y} \in N_i^\prime$ then $0 \neq \overline{\Psi_i(y)} \in N_i^\prime$.}}\leqno(\refn{ord=ord}.3)
$$
To see this we have to look back at the definition of~$\Psi_i$. For each $F_\nu\colon M_\nu \to M_{\nu+1}$ (with $\nu\in\cI$) there are two possibilities. Either $(F_\nu \bmod p) \neq 0$. If this happens and $y \in M_\nu$ with $0 \neq \overline{y} \in N_\nu^\prime$ then it follows from (\refn{ord=ord}.1) that also $0 \neq \overline{F_\nu(y)} \in N_{\nu+1}^\prime$. The other possibility, which occurs precisely $a_1$ times, is that $(F_\nu \bmod p) =0$. In this case, let $y^\prime = p^{-1} \cdot F_\nu(y)$, which is the unique element with $V_\nu(y^\prime) = y$. Looking at the structure of the Dieudonn\'e module~$\ul{N} = \ul{N}^\ord(d,\gf)$, as made explicit in~\refn{XordNPol}, we see that $0 \neq \overline{y^\prime} \in N_{\nu+1}^\prime$. Completing one full loop through the index set~$\cI$ we arrive at~(\refn{ord=ord}.3). 

Let $W=W(k)$. Choose a $W$-basis $f_1,\ldots,f_q,f_{q+1},\ldots,f_d$ for $M_i$ such that 
$$
(W\cdot f_1 + \cdots + W\cdot f_q) \otimes_W k \isomarrow N_i^\prime
\qquad\hbox{and}\qquad
(W\cdot f_{q+1} + \cdots + W\cdot f_d) \otimes_W k \isomarrow N_i^\pprime\, .
$$
Properties (\refn{ord=ord}.2) and (\refn{ord=ord}.3) imply that the matrix of $\Psi_i$ with respect to this basis has the form
$$
\pmatrix{B_1&pB_2\cr pB_3&pB_4\cr}
$$
with $B_1$ an invertible matrix (of size $q\times q$). Further, the chosen basis gives rise to a bijection
$$
M_{(d-q) \times q}(W) \longisomarrow 
\left\{\vcenter{
\setbox0=\hbox{$W$-submodules $U \subseteq M_i$ of rank~$q$,}
\copy0
\hbox to\wd0{\hfil with $U/(pM_i\cap U) \isomarrow N_i^\prime$\hfil}
}\right\}
$$
by sending a matrix $A = (a_{t,u})_{q+1\leq t\leq d,1\leq u \leq q}$ to
$$
U_A := \Span\big(f_1 + p\cdot \sum_{j=q+1}^d a_{j1} f_j, \ldots, f_q + p\cdot \sum_{j=q+1}^d a_{jq} f_j\big)\, .
$$
A straightforward computation shows that $\Psi_i(U_A) = U_{A^\prime}$ with
$$
A^\prime = (B_3 + pB_4\cdot {}^{\sigma^m}\!\!A)(B_1 + p^2B_2\cdot {}^{\sigma^m}\!\!A)^{-1}\, .
$$
Hence $U_A$ is stable under $\Psi_i$ if and only if $A \cdot (B_1 + p^2B_2\cdot {}^{\sigma^m}\!\!A) = (B_3 + pB_4\cdot {}^{\sigma^m}\!\!A)$. Lemma~\refn{eqLem} tells us that there exists an~$A$ with this property.

We proceed as follows. Choose a matrix $A$ such that $M_i^\prime := U_A$ is stable under~$\Psi_i$. Define submodules $M_{i+\nu}^\prime \subseteq M_{i+\nu}$ by an inductive procedure: if $F = 0$ on $M_{i+\nu}$, let $M_{i+\nu+1}^\prime := p^{-1} \cdot F(M_{i+\nu})$; otherwise let $M_{i+\nu+1}^\prime := F(M_{i+\nu})$. By construction, $M_{i+m}^\prime = M_i$, so that we have well-defined submodules $M_i^\prime \subseteq M_i$ for all $i \in \cI$. The direct sum $M^\prime := \oplus_{i\in\cI} M_i^\prime$ is a sub-crystal of~$M$. On the level of BT with $\cO$-structure we find that $\ul{X}$ sits 
in an exact sequence
$$
0 \tto \ul{X}^\pprime \tto \ul{X} \tto \ul{X}^\prime \tto 0\, ,\leqno(\refn{ord=ord}.4)
$$
such that the Dieudonn\'e module of $\ul{X}^\prime$ is~$\ul{M}^\prime$. The $p$-kernel of~$\ul{X}^\prime$ (resp.\ of~$\ul{X}^\pprime$) is the one given by the Dieudonn\'e module $N^\prime$ (resp.\ the Dieudonn\'e module $N/N^\prime \cong N^\pprime$). Hence $\ul{X}^\prime$ and $\ul{X}^\pprime$ are both $[p]$-ordinary. By induction on the number~$r$ we may assume that $\ul{X}^\prime$ and $\ul{X}^\pprime$ are $\mu$-ordinary, and as $\ul{X}$ is isogenous to $\ul{X}^\prime \times \ul{X}^\pprime$ it follows that also $\ul{X}$ is $\mu$-ordinary. (Note that (\refn{ord=ord}.4) splits, as now follows using (a) $\Rightarrow$~(c).)

To finish the proof of the theorem we have to consider the case~$r=1$. If $r=1$ then for all $i\in\cI$, either $\gf(i)=0$ or $\gf(i)=d$. But for such a type there is only one BT$_1$ with $\kappa$-structure, up to isomorphism. So if $\ul{X}$ is any BT with $\cO$-structure with~$r=1$ then $\ul{X}$ is $[p]$-ordinary. On the other hand, in the first part of the proof we have seen that $\ul{X}$ is also $\mu$-ordinary. Finally, in the last part of the proof we have seen that $\Phi_i = p^{a_1} \cdot \Psi_i$ with $\Psi_i\colon M_i \to M_i$ a $\sigma^m$-linear bijection. (If $r=1$ then $M = M^\prime$.) As $k=\kbar$ we can choose a 
$W(k)$-basis $e_{i,1},\ldots,e_{i,d}$ for $M_i$ such that $\Psi_i(e_{i,j}) = e_{i,j}$. Now we choose bases for the character spaces $M_{i+n}$, for $n=1,\ldots,m$, by induction:
$$
e_{i+n+1,j} := \cases{F_{i+n}(e_{i+n,j})&if $F_{i+n} \neq 0 \bmod p$;\cr
p^{-1} F_{i+n}(e_{i+n,j})&if $F_{i+n} = 0 \bmod p$.}
$$
Our choice of the $e_{i,j}$ is such that $e_{i+m,j} = e_{i,j}$. The conclusion is that $\ul{M}$ is isomorphic to the standard ordinary Dieudonn\'e module of the given type. This completes the proof. \QED

\sssection{BTordDef}
{\it Definition.\/} --- Let $K$ be a field of characteristic~$p$. Let $\cO$ be a finite unramified extension of~$\Zp$. If $\ul{X}$ is a BT with $\cO$-structure over~$K$ then we say that $\ul{X}$ is {\it ordinary\/} if $\ul{X} \otimes_K k$ satisfies the equivalent conditions of~\refn{ord=ord} for some (equivalently: every) field~$k = \kbar \supset K$.

\sssection{OrdGenRem}
{\it Remark.\/} --- Thus far we have defined ordinariness only for BT with $\cO$-structure, where $\cO$ is a finite unramified extension of~$\Zp$. We can extend this, in an obvious way, to the situation where $\cO$ is a maximal order in a product of matrix algebras over finite unramified extensions of~$\Qp$. We leave this to the reader.

\sssection{OrdClassic}
{\it Remark.\/} --- Let $\ul{X} = (X,\iota)$ be an ordinary BT with $\cO$-structure over~$k =\kbar$, of type $(d,\gf)$. As usual we let $m = [\cO:\Zp]$. Then $X$, the underlying BT (without $\cO$-structure), is {\it not\/}, in general, ordinary in the classical sense. In fact, using our explicit description of $\ul{X}^\ord(d,\gf)$ and the notation of~\refn{EtMult}, we find
\item{(a)} $X \cong (\Qp/\Zp)^{md} \;\Leftrightarrow\; \ul{X} \cong \ul{X}_\et^d  \;\Leftrightarrow\; \gf(i) = 0\quad \hbox{for all $i\in\cI$}$;
\item{(b)} $X \cong (\hat\mG_{{\rm m}})^{md} \;\Leftrightarrow\; \ul{X} \cong \ul{X}_\mult^d \;\Leftrightarrow\;  \gf(i) = d\quad \hbox{for all $i\in\cI$}$;
\item{(c)} $X$ is ordinary if and only if $\gf$ is a constant function on~$\cI$.

\noindent
In terms of the ``mod $p$ reflex field'' $E(d,\gf)$ defined in~\refn{modpRefl} we find that $X$ is ordinary if and only if $E(d,\gf) = \Fp$. This result is to be compared with [\WedhOrd],~Thm.~1.6.3. 

\sssection{EndoLem2}
{\it Lemma. --- Let $\cO$ be a finite unramified extension of~$\Zp$ with residue field $\kappa \cong \mF_{p^m}$. Fix a type $(d,\gf)$, and let
$$
\ul{X}^\ord(d,\gf) = \ul{X}^\ord(d^1,\gf^1) \times \cdots \times \ul{X}^\ord(d^r,\gf^r)\, .
$$
be the decomposition of the corresponding standard ordinary BT over $k=\kbar$, as in\/ {\rm \refn{MordDec}}. Then
$$
\eqalign{
\End\big(\ul{X}^\ord(d,\gf)\big) &= \End\big(\ul{X}^\ord(d^1,\gf^1)\big) \times \cdots \times 
\End\big(\ul{X}^\ord(d^r,\gf^r)\big)\cr
&\cong M_{d^1}\big(W(\kappa)\big) \times \cdots \times M_{d^r}\big(W(\kappa)\big)\, .}
$$\vskip-\lastskip\smallskip}
\Dskip

\Proof~Easy exercise. \QED

\sssection{HDecCor}
{\it Corollary. --- Let $K$ be a perfect field, $\charact(K) = p$. Let $\cO$ be a finite unramified extension of~$\Zp$. Let $\ul{X}$ be an ordinary BT with $\cO$-structure over~$K$. Then there is a canonical decomposition
$$
\ul{X} = \ul{X}^{(1)} \times \cdots \times \ul{X}^{(r)}
$$
over~$K$ such that, for $k \supset K$ algebraically closed, $\ul{X}^{(\nu)} \otimes_K k \cong \ul{X}^\ord(d^\nu,\gf^\nu)$.}
\Dskip

\Proof~This follows from \refn{EndoLem2} using descent. \QED
\section{Deformation theory of ordinary objects}{DefoTh}

\ssection{Deformation theory of BT with endomorphisms}{DefoBTEnd}

\sssection{DefoFun}
Let $K$ be a perfect field of characteristic $p>0$. Write $\catC_{W(K)}$ for the category of pairs $(R,\beta)$, with $R$ an artinian local $W(K)$-algebra and $\beta\colon R/\gm_R \isomarrow K$ an isomorphism. Morphisms in $\catC_{W(K)}$ are local homomorphisms of $W(K)$-algebras which are compatible with the given isomorphisms~$\phi$.

Let $(R,\beta) \in \catC_{W(K)}$. If $\ul{X}$ is a BT$_n$ with $\cO$-structure over~$K$ then by a {\it deformation\/}, or {\it lifting\/}, of~$\ul{X}$ over~$(R,\beta)$ we mean a pair $(\ulcX,\alpha)$ where $\ulcX$ is a BT$_n$ with $\cO$-structure over~$R$ and $\alpha\colon \ulcX \otimes_{R,\beta} K \isomarrow \ul{X}$. In practice we often omit~$\beta$ from the notation.

The formal deformation functor of~$\ul{X}$ is the covariant functor $\Def(\ul{X})\colon \catC_{W(K)} \to \Sets$ whose value 
on a pair $(R,\beta)$ is the set of isomorphism classes of deformations of~$\ul{X}$ over~$R$. 

\sssection{DefoSit}
{\it Situation.} --- Let $K$ be a perfect field of characteristic $p>0$. Let $\cO$ be an unramified extension of~$\Zp$ of degree~$m$ with residue field~$\kappa$; in other words: $\kappa \cong \mF_{p^m}$ and $\cO \cong W(\kappa)$. Let $\ul{X}$ be an {\it ordinary\/} BT with $\cO$-structure over~$K$, of type~$(d,\gf)$. Let $\ul{Y} := \ul{X}[p]$. We write $r$ for the number of slopes, as defined in~\refn{ajDef}.

\sssection{ProRepProp}
{\it Proposition\/ {\rm (Wedhorn, [\Wedh])}. --- Situation as above. Write $T$ (resp.\ $T^D$) for the tangent space of~$X$ (resp.~$X^D$) at the origin.

{\rm (\romno1)} The functor\/ $\Def(\ul{X})$ is pro-representable and formally smooth over~$W(K)$.

{\rm (\romno2)} The functor\/ $\Def(\ul{Y})$ is formally smooth over~$W(K)$. The canonical map $\gamma\colon \Def(\ul{X}) \to \Def(\ul{Y})$ is a hull.

{\rm (\romno3)} The tangent spaces of\/ $\Def(\ul{X})$ and\/ $\Def(\ul{Y})$ are both canonically isomorphic to 
the $K$-vector space\/ $T^D \otimes_{\kappa \otimes_{\Fp} K} T$, and via these identifications $\gamma$ induces 
the identity map on tangent spaces.}

\sssection{TangSp}
In the situation of the proposition, the tangent space of $\Def(\ul{X})$ can also be described in terms of the 
Dieudonn\'e module $\ul{N}$ of~$\ul{Y}$. Namely, there are natural isomorphisms $T \cong N[F] := \Ker(F)$ and $T^D \cong 
\big(N/N[F]\big)^\vee$. This gives that over $k = \kbar$ the tangent space of $\Def(\ul{X})$ is isomorphic to 
$\dirsum_{i\in\cI} \Hom\big(N_i/N_i[F],N_i[F]\big)$. In particular, the relative dimension of $\Def(\ul{X})$ over~$W(K)$ is equal to $\sum_{i\in\cI} \gf(i)\cdot \big(d-\gf(i)\big)$.

\sssection{1slopeRig}
{\it Corollary. --- Suppose $\ul{X}$ is isoclinic, meaning that $r=1$. Then $\ul{X}$ is rigid, i.e., $\Def(\ul{X})$ is 
pro-represented by~$W(K)$.}
\Dskip

\Proof~The assumption means that for all $i\in\cI$ either $\gf(i)=0$ or $\gf(i)=d$. 
\QED

\sssection{4tups}
Let $k$ be an algebraically closed field, $\charact(k) = p$. Write $\sigma$ for the Frobenius automorphism 
of~$W(k)$. Let $A$ be a formally smooth $W(k)$-algebra. Fix an endomorphism $\phi_A$ lifting the Frobenius 
endomorphism of $A_0 := A/pA$ and with $(\phi_A)_{|W(k)} = \sigma$. 

Consider $4$-tuples $\big(\cM,\Fil^1(\cM),\nabla,F_\cM\big)$ with
\item{---}$\cM$ a free $A$-module of finite rank;
\item{---}$\Fil^1(\cM) \subset \cM$ a direct summand;
\item{---}$\nabla\colon \cM \to \cM \otimes \hat\Omega_{A/W(k)}$ an integrable, topologically quasi-nilpotent 
connection;
\item{---}$F_\cM\colon \cM \to \cM$ a $\phi_A$-linear endomorphism,

\noindent
such that, writing $\cMtilde := \cM + p^{-1} \Fil^1(\cM)$ and $\cM_0 := \cM \otimes_A A_0$,
\item{---}$F_\cM$ induces an isomorphism $\cMtilde \otimes_{A,\phi_A} A \isomarrow \cM$, and
\item{---}$\Fil^1(\cM) \otimes_A A_0 = \Ker\big(F_\cM \otimes \Frob_{A_0} \colon \cM_0 \to \cM_0\big)$.

\noindent
With the obvious notion of a morphism, such $4$-tuples form a category $\catMF^\nabla_{[0,1]}(A)$. Crystalline 
Dieudonn\'e theory establishes an anti-equivalence
$$
\big(\hbox{BT over~$A$}\big)
\mapright{\hbox{{\sevenrm anti-eq}}}
\catMF^\nabla_{[0,1]}(A)\, 
.\leqno(\refn{4tups}.1)
$$
See [\Durh], \S~4 for further discussion.

\sssection{UnivDefo}
We need a description of the $4$-tuple $\big(\cM,\Fil^1(\cM),\nabla,F_\cM\big)$ corresponding to the universal deformation of an ordinary BT with $\cO$-structure, $\cO$ a finite unramified extension of~$\Zp$. We use a construction due to Faltings~[\Falt]. We start with the Dieudonn\'e module $\ul{M}$ of~$\ul{X}$. If $(d,\gf)$ is the type of~$\ul{X}$ then on the standard basis $\{e_{i,j}\}$ for~$M$ (as in~\refn{XordNPol}) the Hodge filtration is given by
$$
\Fil^1(M) = \dirsum_{i\in\cI, j > d-\gf(i)} W(k) \cdot e_{i,j}\, .
$$
The submodule $M^0 := \dirsum_{i\in\cI, j \leq d-\gf(i)} W(k) \cdot e_{i,j}$ is a complement for~$\Fil^1(M)$ in~$M$.

Inside the algebraic group $\GL_{\cO \otimes_{\Zp} W(k)}(M) = \prod_{i\in\cI} \GL_{d,W(k)}$, the stabilizer of 
$M^0$ is a parabolic subgroup. Let $U$ be its unipotent radical. In down-to-earth terms, $U = \prod_{i\in\cI} 
U_i$ with
$$
U_i = \hbox{{\rm group of matrices}}\quad \pmatrix{1&\ast\cr 0&1\cr}\, ,
$$
where the upper right-hand block has size $\big(d-\gf(i)\big) \times \gf(i)$.

Let $A$ be the completed local ring of $U$ at the identity element. Thus, $A$ is a formal power series ring
$$
A = W(k)\Big[\!\Big[ u^{(i)}_{r,s}\Big]\!\Big]\, ,
$$
where the indices range over $i\in\cI$, over $r \in\{d+1-\gf(i),\ldots,d\}$ and $s\in\{1,\ldots,d-\gf(i)\}$. It 
will be convenient to formally put $u^{(i)}_{r,s} = 0$ if $r$ or~$s$ is not in the specified range. Define $\phi_A\colon A \to A$ to be the lifting of Frobenius with $(\phi_A)_{|W(k)} = \sigma$ and $\phi_A(u^{(i)}_{r,s}) = \big(u^{(i)}_{r,s}\big)^p$.

Note that we have a tautological element $g^\univ \in \GL_{\cO \otimes_{\Zp} W(k)}(M)(A)$.

\sssection{DefoFalt}
{\it Proposition\/ {\rm (Faltings, [\Falt], \S~7)}. --- Define
$$
\cM := M \otimes_{W(k)} A\, ,\quad 
\Fil^1(\cM) := \Fil^1(M) \otimes_{W(k)} A\, ,\quad\hbox{and}\quad
F_\cM := g^\univ \circ \big(F_M \otimes \phi_A\big)\, .
$$
Then there is a unique topologically quasi-nilpotent connection $\nabla \colon \cM \to \cM \otimes 
\hat\Omega_{A/W(k)}$ that is compatible with~$F_\cM$, and this connection is integrable. The ring $\cO$ acts on 
the $4$-tuple $\big(\cM,\Fil^1(\cM),\nabla,F_\cM\big)$ by endomorphisms. Via the anti-equivalence of 
categories\/ {\rm (\refn{4tups}.1)} this $4$-tuple corresponds to the universal deformation of~$\ul{X}$ as a 
BT with $\cO$-action.}
\Cskip

This result is also discussed in [\Durh], \S~4.

\sssection{FiltLift}
{\it Proposition. (Lifting of the slope filtration.) --- Let $K$ be a perfect field, $\charact(K) = p>0$. Let $\cO$ be a 
finite unramified extension of~$\Zp$ with residue field $\kappa \cong \mF_{p^m}$. Let $\ul{X}$ be an ordinary BT with 
$\cO$-structure over~$K$. Consider the decomposition of~$\ul{X}$ as in\/ {\rm \refn{HDecCor}}, and 
define (for $1\leq a \leq b\leq r$)
$$
\ul{X}^{(a,b)} := \prod\nolimits_{j=a}^b \ul{X}^{(j)}\, ,
$$
so that we have a slope filtration
$$
0 \subset \ul{X}^{(r,r)} \subset \ul{X}^{(r-1,r)} \subset \cdots \subset \ul{X}^{(1,r)} = \ul{X}\, 
.
$$
Let $R$ be in $\catC_{W(K)}$. If $(\ulcX,\alpha)$ is a deformation of $\ul{X}$ over~$R$ then there is a unique 
filtration of $\ulcX$ by sub-BT with $\cO$-structure,
$$
0 \subset \ulcX^{(r,r)} \subset \ulcX^{(r-1,r)} \subset \cdots \subset \ulcX^{(1,r)} = \ulcX
$$
such that $\alpha\colon \ulcX \otimes_R K \isomarrow \ul{X}$ restricts to isomorphisms $\ulcX^{(a,r)} \otimes_R 
K \isomarrow \ul{X}^{(a,r)}$.}
\Dskip

\Proof~The unicity is an easy consequence of Grothendieck-Messing deformation theory. It suffices to prove 
the existence in the case that $K = k$ is algebraically closed; the general case then follows by descent, using 
the unicity of the lifted filtration. Further it clearly suffices to show that the desired filtration exists in 
case $\ulcX$ is the universal deformation of~$\ul{X}$. Finally it suffices to show that $\ul{X}^{(2,r)}$ lifts 
to a sub-object of~$\ulcX$.

We use the description of the universal deformation of~$\ul{X}$ given in~\refn{UnivDefo} and~\refn{DefoFalt}. 
Let $\ul{M}^\prime \subset \ul{M}$ be the submodule corresponding to the quotient $\ul{X} \twoheadrightarrow 
\ul{X}/\ul{X}^{(2,r)} = \ul{X}^{(1)}$. Let $d^\prime = \min\{d-\gf(i)\mid i\in\cI\}$ be the rank of~$M^\prime$ 
over $\cO \otimes_{\Zp} W(k)$, so that
$$
M^\prime = \dirsum_{i\in\cI} \Span\Big(e_{i,j};\; j=1,\ldots,d^\prime\Big)\, .
$$
Set $\cM^\prime := M^\prime \otimes_{W(k)} A$ and $\Fil^1(\cM^\prime) := \Fil^1(M^\prime) \otimes_{W(k)} A = 
\cM^\prime \cap \Fil^1(\cM)$. Note that the tautological element $g^\univ$ acts trivially on~$\cM^\prime$, so 
$F_\cM$ restricts to $F_{\cM^\prime} := F_{M^\prime} \otimes \phi_A$ on~$\cM^\prime$.

We claim that $\nabla(\cM^\prime) \subseteq \cM^\prime \otimes \hat\Omega_{A/W(k)}$, i.e., $\nabla$ restricts to 
a connection~$\nabla^\prime$ on~$\cM^\prime$. If this is true then 
$\big(\cM^\prime,\Fil^1(\cM^\prime),\nabla^\prime,F_{\cM^\prime}\big)$ is a sub-object of 
$\big(\cM,\Fil^1(\cM),\nabla,F_\cM\big)$ and the desired lifting of $\ul{X}^{(2,r)}$ is the one corresponding 
to the quotient-crystal.

In the proof of the claim we use the standard basis $\{e_{i,j}\}$ for the module~$\cM$. Let $\{f_{i,j}\}$ be the 
corresponding basis for $\cMtilde \otimes_{A,\psi_A} A$; concretely,
$$
f_{i,j} = \cases{e_{i,j} \otimes 1&if $j> d-\gf(i)$;\cr
p^{-1}e_{i,j} \otimes 1&if $j \leq d -\gf(i)$.}
$$
With respect to these bases Frobenius is given by
$$
F_\cM(f_{i,j}) = g^\univ_{i+1}(e_{i+1,j})\, ,\qquad g^\univ_{i+1} = \pmatrix{1&u^{(i+1)}_{r,s}\cr 0&1\cr}\, .
$$
The connection $\nabla$ decomposes into factors $\nabla_i \colon \cM_i \to \cM_i \otimes \hat\Omega_{A/W(k)}$; 
let $D^{(i)}$ be the $d\times d$ matrix of $1$-forms of~$\nabla_i$ on the basis $e_{i,1},\ldots,e_{i,d}$. 
Similarly, let $\tilde D^{(i)}$ be the matrix of the induced connection $\tilde\nabla_i$ on $\cMtilde_i 
\otimes_{A,\phi_A} A$ with respect to the basis $f_{i,1},\ldots,f_{i,d}$. Then $\tilde D^{(i)}$ is obtained from 
$D^{(i)}$ by applying ${\rm d}\phi_A$ to all coefficients.

By definition, $\nabla$ is compatible with Frobenius, meaning that $(F \otimes \id) \circ \tilde\nabla = \nabla 
\circ F$. But
$$
\nabla \circ F(f_{i,j}) = \sum_{\mu=1}^d e_{i+1,\mu} \otimes D^{(i+1)}_{\mu,j} + \sum_{r,\mu=1}^d 
u^{(i+1)}_{r,j} e_{i+1,\mu} \otimes D^{(i+1)}_{\mu,r} + \sum_{r=1}^d e_{i+1,r} \otimes {\rm d}u^{(i+1)}_{r,j}\, 
,
$$
whereas
$$
(F\otimes\id) \circ \tilde\nabla (f_{i,j}) = \sum_{\mu=1}^d \Big\{e_{i+1,\mu} + \sum_{s=1}^d u^{(i+1)}_{s,\mu} 
e_{i+1,s}\Big\} \otimes {\rm d}\phi_A\Big(D^{(i)}_{\mu,j}\Big)\, .
$$
Here recall that we formally put $u^{(i)}_{r,s} = 0$ if either $r \leq d-\gf(i)$ or $s>d-\gf(i)$. Comparing 
coefficients of $e_{i+1,\nu}$ (for fixed $\nu \in \{1,\ldots,d\}$) we find
$$
D^{(i+1)}_{\nu,j} + \sum_{r=1}^d u^{(i+1)}_{r,j} D^{(i+1)}_{\nu,r} + {\rm d}u^{(i+1)}_{\nu,j} = {\rm d}\phi_A 
\Big(D^{(i)}_{\nu,j}\Big) + \sum_{\mu=1}^d u^{(i+1)}_{\nu,\mu} \cdot {\rm d}\phi_A \Big(D^{(i)}_{\mu,j}\Big)\, 
.\leqno(\refn{FiltLift}.1)
$$

Now we specialize to the case that $\nu \leq d^\prime$ and $j<d^\prime$. As $u^{(i+1)}_{\nu,j} = 0$ and 
$u^{(i+1)}_{r,j} \neq 0$ only for $r>d^\prime$, (\refn{FiltLift}.1) becomes
$$
D^{(i+1)}_{\nu,j} = {\rm d}\phi_A \Big(D^{(i)}_{\nu,j}\Big) + \sum_{\mu=1}^d u^{(i+1)}_{\nu,\mu} \cdot {\rm 
d}\phi_A \Big(D^{(i)}_{\mu,j}\Big) - \sum_{r=d^\prime+1}^d u^{(i+1)}_{r,j} D^{(i+1)}_{\nu,r}\, .
$$
By induction on~$n$ this readily implies that for all $i\in\cI$, all $\nu \leq d^\prime$ and $j>d^\prime$ we 
have $D^{(i)}_{\nu,j} \in \gm_A^n \cdot \hat\Omega_{A/W(k)}$. But $A$ is noetherian, so $\cap_{n\geq 0}\; 
\gm_A^n \cdot \hat\Omega_{A/W(k)} = 0$, so $D^{(i)}_{\nu,j} = 0$ whenever $\nu\leq d^\prime$ and $j>d^\prime$. 
This is what we wanted to prove. \QED

\ssection{Cascades}{Casc}

\sssection{CascDef}
{\it Definition.\/} --- Let $\cT$ be a topos with final object~$S$. Let~$r$ be a positive integer. An {\it 
$r$-cascade\/} in~$\cT$ consists of the following data:

(1) commutative $\cT$-groups $G^{(i,j)}$ for $1\leq i < j\leq r$;

(2) objects $\Gamma^{(i,j)}$ for $1\leq i< j \leq r$; if $i \geq j$ then we put $\Gamma^{(i,j)} := S$;

(3) morphisms $\lambda^{(i,j)}\colon \Gamma^{(i,j)} \to \Gamma^{(i,j-1)}$ and $\rho^{(i,j)}\colon \Gamma^{(i,j)} 
\to \Gamma^{(i+1,j)}$ satisfying the commutativity relation $\rho^{(i,j-1)} \circ \lambda^{(i,j)} = 
\lambda^{(i+1,j)} \circ \rho^{(i,j)}$;

(4) the structure on $\Gamma^{(i,j)}$ of a biextension of $\big(\Gamma^{(i,j-1)},\Gamma^{(i+1,j)}\big)$ by 
$G^{(i,j)} \times \Gamma^{(i+1,j-1)}$ in the category $\cT_{/\Gamma^{(i+1,j-1)}}$.
\Cskip

Part (4) of the data is meaningful by induction on $j-i$. If $j=i+1$ then, by convention, $\Gamma^{(i,j-1)} = 
\Gamma^{(i+1,j)} = S$ and (4) means that $\Gamma^{(i,j)}$ is to be equipped with the structure of a commutative 
$\cT$-group isomorphic to~$G^{(i,j)}$. If $j=i+m$ and data as in (4) are available on all 
$\Gamma^{(i^\prime,j^\prime)}$ with $j^\prime-i^\prime < m$ then $\Gamma^{(i,j-1)}$ and $\Gamma^{(i+1,j)}$ both 
have the structure of a commutative group over $\Gamma^{(i+1,j-1)}$ (as part of their structure of a 
biextension), so that (4) is meaningful for~$\Gamma^{(i,j)}$.

\sssection{CascExa}
{\it Example.\/} --- A $1$-cascade only consists of the final object~$S$. A $2$-cascade is just a commutative $\cT$-group. A 
$3$-cascade is a biextension. A $4$-cascade is a commutative diagram
$$
\matrix{&&&&\Gamma^{(1,4)}&&&&\cr
&&&\swarrow&&\searrow&&&\cr
&&\Gamma^{(1,3)}&&&&\Gamma^{(2,4)}&&\cr
&\swarrow&&\searrow&&\swarrow&&\searrow&\cr
\Gamma^{(1,2)}&&&&\Gamma^{(2,3)}&&&&\Gamma^{(3,4)}\cr}
$$
with
\item{---}$\Gamma^{(i,i+1)}$ commutative $\cT$-groups isomorphic to $G^{(i,i+1)}$ (for $i=1$, $2$,~$3$);
\item{---}$\Gamma^{(1,3)}$ a biextension of $\big(\Gamma^{(1,2)},\Gamma^{(2,3)}\big)$ by $G^{(1,3)}$;
\item{---}$\Gamma^{(2,4)}$ a biextension of $\big(\Gamma^{(2,3)},\Gamma^{(3,4)}\big)$ by $G^{(2,4)}$;
\item{---}$\Gamma^{(1,4)}$ a biextension of $\big(\Gamma^{(1,3)},\Gamma^{(2,4)}\big)$ by $G^{(1,4)} \times 
\Gamma^{(2,3)}$ in the category $\cT_{/\Gamma^{(2,3)}}$.

We often refer to a cascade by the single letter~$\Gamma$. We call the $G^{(i,j)}$ the {\it 
group-constituents\/} of the cascade and $\Gamma^{(i,j)}$ the {\it $(i,j)$-truncation\/}.

If $\Gamma$ is an $r$-cascade, $r\geq 2$, then $\Gamma^{(1,r-1)}$ and $\Gamma^{(2,r)}$ both inherit a natural 
structure of an $(r-1)$-cascade.

\sssection{HomCascDef}
{\it Definition.\/} --- Let $\Gamma$ and $\Delta$ be $r$-cascades, with group constituents $G^{(i,j)}$ and 
$H^{(i,j)}$, respectively. A {\it homomorphism $f\colon \Gamma \to \Delta$ of cascades\/} is a collection of 
maps $f^{(i,j)}\colon \Gamma^{(i,j)} \to \Delta^{(i,j)}$ and homomorphisms of groups $h^{(i,j)}\colon G^{(i,j)} 
\to H^{(i,j)}$ satisfying the following two conditions. 

(a)~The maps $f^{(i,j)}$ are compatible, in the obvious sense, with the given morphisms $\lambda$ and~$\rho$; 
symbolically we may write this as the conditions that $f \circ \lambda_\Gamma = \lambda_\Delta \circ f$ and $f 
\circ \rho_\Gamma = \rho_\Delta \circ f$.

(b)~For all $1\leq i< j \leq r$, let a ``$\; \tilde{ }\; $'' denote a base change via $f^{(i+1,j-1)}$, and write 
$\tilde f^{(\mu,\nu)}\colon \Gamma^{(\mu,\nu)} \to \widetilde\Delta^{(\mu,\nu)}$ for the morphism over 
$\Gamma^{(i+1,j-1)}$ induced by $f^{(\mu,\nu)}$. Then the quadruple $\big(\tilde f^{(i,j-1)}, \tilde 
f^{(i+1,j)}, h^{(i,j)}, \tilde f^{(i,j)}\big)$ is a homomorphism of biextensions over~$\Gamma^{(i+1,j-1)}$.   

\sssection{TorsPtDef}
{\it Definition.\/} --- Let $\Gamma$ be an $r$-cascade. Let $x \in \Gamma^{(1,r)}(R)$ be an $R$-valued point, 
for some $R\in\cT$. By induction on~$r$ we define what it means for~$x$ to be a {\it torsion point\/}: If $r=1$ 
then every~$x$ is torsion. If $r\geq 2$ then we say that $x$ is a torsion point if 

(a)~$\lambda(x) \in \Gamma^{(1,r-1)}(R)$ and $\rho(x) \in \Gamma^{(2,r)}(R)$ are torsion points;

(b)~$x$ is a torsion point of $\Gamma^{(1,r)}$ viewed as a group over $\Gamma^{(1,r-1)}$, and also a torsion 
point of~$\Gamma^{(1,r)}$ viewed as a group over~$\Gamma^{(2,r)}$.
\Cskip

Note that if (a) holds then in (b) it suffices to require that $x$ is a torsion point for one of the two group 
laws. To see this one uses that a biextension of a pair of groups $(\Gamma_1,\Gamma_2)$ by a third group has a 
canonical trivialization over~$\{0\} \times \Gamma_2$ and over $\Gamma_1 \times \{0\}$.

\sssection{0ofCasc}
An $r$-cascade $\Gamma$ has a natural zero section $0 \in \Gamma(S)$. As in the above we define it by induction 
on~$r$. We leave the details to the reader.  

\sssection{DualCasc}
{\it Definition.\/} --- Let $\Gamma$ be an $r$-cascade. Then we define the {\it dual cascade\/}, notation 
$\Gamma^\vee$, to be the $r$-cascade obtained after replacing all index pairs $(i,j)$ by $(r+1-j,r+1-i)$. Thus, 
the group constituents are $G^{\vee,(i,j)} := G^{(r+1-j,r+1-i)}$, the truncations are $\Gamma^{\vee,(i,j)} := 
\Gamma^{(r+1-j,r+1-i)}$, and the biextension structures are the same (after re-indexing) as those occurring 
in~$\Gamma$.

\ssection{The cascade structure on the deformation space}{CascDefo}

\sssection{Ext/R}
Situation as in~\refn{DefoSit}. We first study the case that the Newton polygon of~$\ul{X}$ has precisely two slopes, i.e., $r=2$. By \refn{HDecCor} we have a decomposition $\ul{X} = \ul{X}^{(1)} \times \ul{X}^{(2)}$ where the factors are both isoclinic. As shown in~\refn{1slopeRig}, if $R \in \catC_{W(K)}$ then there is a unique lifting $\ulcX^{(\nu)}$ of $\ul{X}^{(\nu)}$ over~$R$, for $\nu \in\{1,2\}$. If we want to indicate over which ring~$R$ we are working we use the notation~$\ulcX^{(\nu)}_R$.

Consider the category $\catEXT\big(\ulcX^{(1)},\ulcX^{(2)}\big)$ of extensions of $\ulcX^{(1)}$ by $\ulcX^{(2)}$ as fppf sheaves of $\cO$-modules over~$R$. If there is no risk of confusion we simply write $\catEXT_R$ for this category, and $\Ext_R$ denotes the set of isomorphism classes in~$\catEXT_R$. 

Let $\ulcX$ be an object of $\catEXT_R$. It follows from [\Mess], \Romno 1, (2.4.3) that, forgetting the structure of an extension, $\ulcX$ is again a BT with $\cO$-structure. By looking at the Newton polygon, using [\Katz], Lemma~1.3.4, we see that $\ulcX \otimes_R K$ is ordinary. Applying Thm.~\refn{ord=ord}, and using that here are no non-trivial homomorphisms from $\ul{X}^{(1)}$ to $\ul{X}^{(2)}$, it follows that there is a unique trivialization $\alpha\colon \ulcX \otimes_R K \isomarrow \ul{X}^{(1)} \times \ul{X}^{(2)}$ as extensions. 

\sssection{(1,fprime)}
Let $\ul{X}$ be as in\/ {\rm \refn{DefoSit}}, and assume $r=2$. Let $(d,\gf) = (d^1,\gf^1) + (d^2,\gf^2)$ be the decomposition of the type of $\ul{X}$ as in section~\refn{ajDef}. Define a new type $(d^\prime,\gf^\prime)$ by $d^\prime =1$ and
$$
\gf^\prime(i) = \cases{0&{\rm if $\gf^1(i) = \gf^2(i) = 0$;}\cr
0&{\rm if $\gf^1(i) = d^1$ and $\gf^2(i) = d^2$;}\cr
1&{\rm if $\gf^1(i) = 0$ and $\gf^2(i) = d^2$.}\cr}
$$
As functions on~$\cI$ we have $\gf^2 = \gf^1 + (d^2-d^1) \cdot \gf^\prime$. The corresponding ordinary object $\ul{X}(1,\gf^\prime)$ is isoclinic, so by~\refn{1slopeRig} it has a unique lifting $\ul{X}^\can(1,\gf^\prime)$ to a BT with 
$\cO$-structure over~$W(k)$.

\sssection{2SlThm}
{\it Theorem. --- Let $\ul{X}$ be as in\/ {\rm \refn{DefoSit}} with $r=2$. With notations as above, the category $\catEXT_R$ is equivalent to the category $\catDEF_R(\ul{X})$ of deformations of $\ul{X}$ over~$R$. The functor $\Def(\ul{X})$ has a natural structure of a BT with $\cO$-structure over~$W(K)$. If $k$ is an algebraically closed field containing~$K$ then 
$$
\Def(\ul{X}) \otimes_{W(K)} W(k) \cong \ul{X}^\can(1,\gf^\prime)^{d^1d^2}\, ,\leqno{\rm 
(\refn{2SlThm}.1)}
$$
as BT with $\cO$-structure.}
\Dskip

\Proof~If $\ulcX$ is an object of $\catEXT_R$, let $\alpha\colon \ulcX \otimes_R K \isomarrow \ul{X}^{(1)} 
\times \ul{X}^{(2)}$ be the unique trivialization of $\ulcX \otimes K$ as extension. Forgetting the structure of 
an extension on~$\ulcX$, the pair $(\ulcX,\alpha)$ is a deformation of~$\ul{X}$ over~$R$. One easily checks that 
this defines a functor $h\colon \catEXT_R \to \catDEF_R(\ul{X})$. In the opposite direction, suppose 
$(\ulcX,\beta)$ is a deformation of $\ul{X}$ over~$R$. By~\refn{FiltLift} and rigidity of~$\ul{X}^{(2)}$ we have 
$\ulcX^{(2)} \hookrightarrow \ulcX$. By [\Mess], \Romno 1,~(2.4.3) the quotient $\ulcX/\ulcX^{(2)}$ is again a 
BT with $\cO$-structure. By rigidity of~$\ul{X}^{(1)}$ it then follows that $\ulcX$ is an extension of 
$\ulcX^{(1)}$ by~$\ulcX^{(2)}$, in such a way that the given identification $\beta\colon \ulcX \otimes_R K 
\isomarrow \ul{X}$ becomes a trivialization of the extension over~$K$. This gives a quasi-inverse to the 
functor~$h$.

We find that $\Def_R(\ul{X}) \cong \Ext_R$, which has a natural group structure. Hence $\Def(\ul{X})$ has the 
structure of a smooth formal group over~$W(K)$. Further, $\cO$ acts on it through its action on~$\ul{X}^{(1)}$.

Let us now show that $\Def(\ul{X})$ is a BT. By [\Mess], \Romno 2, (4.3) and~(4.5) it suffices to show that for 
every $R \in \catC_{W(K)}$ multiplication by~$p$ is an epimorphism of $\Def(\ul{X}) \otimes R$ to itself. For 
this it is enough to show that $D := \Def(\ul{X}) \otimes_{W(K)} K$ is a BT over~$K$. Further we may assume that 
$K = k = \kbar$. (Another proof of these reduction steps can be found in Conrad's notes [\Conr], \S~3.1) By the 
classification theory of formal groups (see Manin [\Manin], \Romno 2.4), in order to conclude that $D$ is a BT 
it suffices to show that $D[p]$ is a finite group scheme.

As usual we write $\ul{Y}^{(\nu)}$ for the $p$-kernel of $\ul{X}^{(\nu)}$; similarly, let $\ulcY^{(\nu)}_R := 
\ulcX^{(\nu)}_R[p]$. If $R$ is a $k$-algebra then $\ulcY^{(\nu)}_R = \ul{Y}^{(\nu)} \otimes_k R$.

Let $R$ be an artinian local $k$-algebra. We have an exact sequence
$$
\Hom\Big(\ulcY_R^{(1)},\ulcX_R^{(2)}\Big) \tto \Ext_R \mapright{\times p} \Ext_R \tto \cdots\, 
,\leqno(\refn{2SlThm}.2)
$$
where the $\Hom(\ )$ denotes homomorphisms of sheaves of $\cO$-modules over~$R$. Clearly we have 
$\Hom\big(\ulcY_R^{(1)},\ulcX_R^{(2)}\big) \cong \Hom\big(\ulcY_R^{(1)},\ulcY_R^{(2)}\big)$, and as a functor in 
$k$-algebras~$R$ the latter is representable by a group scheme $G := \bHom\big(\ul{Y}^{(1)},\ul{Y}^{(2)}\big)$. 
But $G$ is isomorphic to a closed subgroup scheme of the group scheme $\bAut\big(\ul{Y}\big)$, which is finite 
by~[\DFEO], Thm.~2.1.2 and our definition of the $[p]$-ordinary type. Hence we have a finite group scheme~$G$ 
and a homomorphism $G \to D[p]$ which is surjective on $R$-valued points for every~$R$ as above. This is easily 
seen to imply that $D[p]$ is finite.

In order to prove the last assertion of the theorem we may assume that $K = k$, and one easily reduces to the 
case that $d^1=d^2=1$. If we can show that $D := \Def(\ul{X}) \otimes_{W(k)} k$ is isomorphic to 
$\ul{X}(1,\gf^\prime)$ then (\refn{2SlThm}.1) follows from the rigidity result~\refn{1slopeRig}.

{}From now on we assume that $d^1=d^2=1$ and $K=k$. We have $\ul{Y}^{(\nu)} = \ul{Y}^\ord(1,\gf^\nu)$. Let 
$\ul{N}_D$ be the Dieudonn\'e module of~$D[p]$ and $\ul{N} = \ul{N}(1,\gf^\prime)$ that of 
$\ul{Y}(1,\gf^\prime)$. We know that $N_D[F]$ is isomorphic to the tangent space of~$D$ at the origin, 
so~\refn{TangSp} gives $N_D[F] \cong N[F]$ as modules over $\cO \otimes k$. By our classification results the 
proof is complete if we can show that~$N_D$, which is a free module over~$\cO \otimes k$, is of rank~$1$. We shall 
use a result of Raynaud~[\Rayn] to prove that the affine algebra of $G := 
\bHom\big(\ul{Y}^{(1)},\ul{Y}^{(2)}\big)$ has $k$-dimension (at most) equal to~$p^m$. Loc.\ cit., Cor.~1.5.1 
tells us that the affine algebras of~$\ul{Y}^{(1)}$ and~$\ul{Y}^{(2)}$ are of the form $A^{(1)} = k[x_i;\; 
i\in\cI]/\ga$ and $A^{(2)} =  k[y_i;\; i\in\cI]/\gb$ with
$$
\ga = \Big(x_i^p - \big(1-\gf^1(i)\big) x_{i+1};\; i\in\cI\Big)\, ,
\quad\hbox{and}\quad
\gb = \Big(y_i^p - \big(1-\gf^2(i)\big) y_{i+1};\; i\in\cI\Big)\, ,
$$
and $a\in \kappa := \cO/p\cO$ acts on $x_i$ and $y_i$ as multiplication by $i(a) \in k$. If $R$ is a $k$-algebra 
then a homomorphism $\ul{Y}^{(1)}_R \to \ul{Y}^{(2)}_R$ is given by a homomorphism $A^{(2)} \otimes_k R \to 
A^{(1)} \otimes_k R$ with $y_i \mapsto \gamma_i \cdot x_i$ for certain $\gamma_i \in R$. For each~$i$ there are 
three possibilities.
\item{(a)}If $\gf^1(i) = 0 = \gf^2(i)$ then the relations $x_i^p = x_{i+1}$ and $y_i^p = y_{i+1}$ give that 
$\gamma_{i+1} = \gamma_i^p$. 
\item{(b)}If $\gf^1(i) = 0$ and $\gf^2(i) = 1$ then we have $x_i^p = x_{i+1}$ and $y_i^p = 0$, hence $\gamma_i^p 
= 0$. 
\item{(c)}The third possibility is that $\gf^1(i) = 1 = \gf^2(i)$. We claim that in this case again 
$\gamma_{i+1} = \gamma_i^p$. To see this we can use Cartier duality, which interchanges cases (a) and~(c). 
Alternatively, we can see that $\gamma_{i+1} = \gamma_i^p$ by using the explicit formulas for the 
comultiplication given by [\Rayn], Cor.~1.5.1.

\noindent
The conclusion is that the affine algebra $A_G$ of~$G$ is a quotient of the ring $k[z_i;\; i\in\cI]/\gc$ with
$$
\gc = \Big(z_i^p - \big(1-\gf^\prime(i)\big) z_{i+1};\; i\in\cI\Big)\, .
$$
This shows that $\dim_k(A_G) \leq p^m$, as claimed. But then $N_D$, which is free, has rank $\leq 1$ over $\cO 
\otimes k$. On the other hand, $N_D[F] \neq 0$. Hence $N_D$ is free of rank~$1$. This completes the proof of the 
theorem. \QED

\sssection{2SlRem}
{\it Remarks.\/} --- (\romno1)\enspace Let $\ul{X} = (X,\iota)$ be as in the theorem. By~\refn{OrdClassic}, $X$ is ordinary 
(in the classical sense) if and only if $\gf$ is a constant function, which means that $\gf^\prime(i)=1$ for all~$i$. This is 
precisely the case that~$\Def(\ul{X})$ is a formal torus. Of course, the formal group structure on~$\Def(\ul{X})$ is in this 
case the same as the one defined by Serre and Tate.

(\romno2)\enspace As a corollary of the proof we find that $\bHom\big(\ul{Y}^{(1)},\ul{Y}^{(2)}\big)$ is a 
$K$-form of the group scheme $\ul{Y}(1,\gf^\prime)^{d^1d^2}$, and that the first map in~(\refn{2SlThm}.2) 
is injective. This last result also follows from the fact that for $R \in \catC_{W(K)}$ we have 
$\Hom\big(\ul{X}^{(1)}_R,\ul{X}^{(2)}_R\big) = 0$, as can be shown using [\Illus], Thm.~4.4.

\sssection{FStopos}
Let $K$ be a perfect field, $\charact(K) = p>0$. Let $W = W(K)$. Write $\catFS_W$ for the category of affine 
formal schemes $\gX$ over~$W$ with the property that $\Gamma(\gX,O_\gX)$ is a pro-finite $W$-algebra. By a 
theorem of Grothendieck, $\catFS_W$ is equivalent to the category of left-exact covariant functors $\catC_W \to 
\Sets$.

On $\catFS_W$ we consider the flat topology; see SGA~3, \Romno 7$_{{\rm B}}$, 1.5. This topology is coarser than 
the canonical topology. Hence if we write $\widehat\catFS_W$ for the topos of sheaves (for the flat topology) 
then we have a natural functor $i\colon \catFS_W \to \widehat\catFS_W$.

\sssection{CascConstr}
{\it Construction.\/} --- Let $\ul{X}$ be as in \refn{DefoSit}. We are going to define the structure of a 
cascade on the formal deformation space of~$\ul{X}$. More precisely, write $\ul{X} = \ul{X}^{(1)} \times \cdots 
\times \ul{X}^{(r)}$ as in~\refn{HDecCor}. For $1 \leq a \leq b \leq r$, let $\ul{X}^{(a,b)} := \prod_{j=a}^b 
\ul{X}^{(j)}$, and define
$$
\Gamma^{(a,b)} := \Def\Big(\ul{X}^{(a,b)}\Big)\, .
$$
Using Prop.~\refn{FiltLift} we obtain natural morphisms $\lambda \colon \Gamma^{(a,b)} \to \Gamma^{(a,b-1)}$ and $\rho\colon \Gamma^{(a,b)} \to \Gamma^{(a+1,b)}$. Finally, define
$$
G^{(a,b)} := \Def\Big(\ul{X}^{(a)} \times \ul{X}^{(b)}\Big)\, ,
$$
viewed as a BT with $\cO$-structure over~$W = W(K)$. We shall define the structure of an $r$-cascade on the 
collection of data $\{\Gamma^{(a,b)},\lambda,\rho\}$, with group constituents $G^{(a,b)}$. Here we work in the 
category~$\widehat\catFS_W$.

To begin with, fix indices $a<b$. We claim that we can choose coordinates such that
$$

\matrix{\Gamma^{(a,b)}\cr
\swarrow\qquad\quad\searrow\cr
\Gamma^{(a,b-1)}\qquad\qquad\Gamma^{(a+1,b)}\cr
\searrow\qquad\quad\swarrow\cr
\Gamma^{(a+1,b-1)}\cr}
\qquad\hbox{{\rm is given by}}\qquad
\matrix{A[\![u_i,v_j,w_1,\ldots,w_f]\!]\cr
\nearrow\qquad\qquad\nwarrow\cr
A[\![u_1,\ldots,u_d]\!]\qquad\qquad A[\![v_1,\ldots,v_e]\!]\cr
\nwarrow\qquad\qquad\nearrow&\cr
A\cr}
$$
where $A \cong W[\![x_1,\ldots,x_h]\!]$ is the affine algebra of $\Gamma^{(a+1,b-1)}$. This is an easy algebra 
exercise; one uses that each $\Gamma^{(a,b)}$ is formally smooth over~$W$, and computes the induced maps on 
tangent spaces in terms of Dieudonn\'e modules. In particular we find that the natural map
$$
\pi\colon \Gamma^{(a,b)} \tto \Gamma^{(a,b-1)} \times_{\Gamma^{(a+1,b-1)}} \Gamma^{(a+1,b)}
$$
is formally smooth, hence topologically flat and surjective.

We want to define on $\Gamma^{(a,b)}$ the structure of a group over $\Gamma^{(a,b-1)}$, as well as the structure 
of a group over $\Gamma^{(a+1,b)}$. Let $R\in \catC_W$, with $W=W(K)$. Let $\eta \in \Gamma^{(a,b-1)}(R)$ 
correspond to a deformation $(\ulcF,\beta)$ of $\ul{X}^{(a,b-1)}$ over~$R$. With similar arguments as in the 
proof of~\refn{2SlThm}, one shows that there is a bijection
$$
\big\{\zeta \in \Gamma^{(a,b)}(R)\bigm| \zeta \mapsto \eta\big\}
\isomarrow \Ext\big(\ulcF,\ulcX^{(b)}_R\big)\, .
$$
The desired group structure on $\Gamma^{(a,b)}$ over $\Gamma^{(a,b-1)}$ is then obtained by transporting the 
group structure on $\Ext\big(\ulcF,\ulcX^{(b)}_R\big)$. For the group structure on $\Gamma^{(a,b)}$ over 
$\Gamma^{(a+1,b)}$ the construction is similar.

The next point is to show, for $b>a+1$, that $\Gamma^{(a,b)}$ has a natural structure of a $G^{(a,b)}$-torsor 
over $\Pi := \Gamma^{(a,b-1)} \times_{\Gamma^{(a+1,b-1)}} \Gamma^{(a+1,b)}$. For this we use Grothendieck's 
notion of a blended extension (``extension panach\'ee''); see SGA~7, \Romno 9,~9.3. We write $\Gamma(?)$ for the 
fibre of $? \in \Pi(R)$ under $\pi\colon \Gamma(R) \to \Pi(R)$.

An $R$-valued point $\eta_1 \times_\xi \eta_2$ of $\Pi$ is given by:
\item{---} a deformation $(\ulcG,\alpha)$ of $\ul{X}^{(a+1,b-1)}$ over~$R$;
\item{---} the class of an extension $0 \tto \ulcG \tto \ulcF_1 \tto \ulcX^{(a)}_R \tto 0$;
\item{---} the class of an extension $0 \tto \ulcX^{(b)}_R \tto \ulcF_2 \tto \ulcG \tto 0$.

\noindent
(The extensions live in the category of sheaves of $\cO$-modules over~$R$.) With similar arguments as before we 
find that the fibre $\Gamma(\eta_1 \times_\xi \eta_2)$ is in natural bijection with the set 
$\Extpan\big(\ulcF_1,\ulcF_2\big)$ of blended extensions of $\ulcF_1$ by $\ulcF_2$. As shown in loc.\ cit., this 
set is either empty or it is is principal homogeneous under the group $\Ext\big(\ulcX^{(a)}_R,\ulcX^{(b)}_R\big) 
= G^{(a,b)}(R)$. But we have already excluded the first option. Using the obvious functoriality with respect 
to~$R$, and using that $\pi\colon \Gamma^{(a,b)} \to \Pi$ is a flat covering, we find that, indeed, 
$\Gamma^{(a,b)}$ is a $G^{(a,b)}$-torsor over~$\Pi$.

Finally we have to give $\Gamma^{(a,b)}$ the structure of a biextension. We use a pointwise notation. Let 
$\eta_1$, $\eta_1^\prime \in \Gamma^{(a,b-1)}(R)$ and $\eta_2$, $\eta_2^\prime \in \Gamma^{(a+1,b)}(R)$, all 
mapping to the same point $\xi \in \Gamma^{(a+1,b-1)}(R)$. The points $\eta_1$ and $\eta_1^\prime$ correspond to 
extensions $\ulcF_1$ and $\ulcF_1^\prime$ of $\ulcX^{(a)}_R$ by~$\ulcG$; the points $\eta_2$ and $\eta_2^\prime$ 
to extensions $\ulcF_2$ and $\ulcF_2^\prime$ of $\ulcG$ by~$\ulcX^{(b)}_R$.  

A point of $\Gamma(\eta_1 \times_\xi \eta_2)$ can be viewed as an extension~$E$ of $\ulcF_1$ by $\ulcX^{(b)}_R$; 
similarly for points of $\Gamma(\eta_1^\prime \times_\xi \eta_2)$. This interpretation gives rise to maps
$$
\phi(\eta_1,\eta_1^\prime;\eta_2) \colon \Gamma(\eta_1 \times_\xi \eta_2) \times \Gamma(\eta_1^\prime \times_\xi 
\eta_2) \tto \Gamma\big((\eta_1 \wedge \eta_1^\prime) \times_\xi \eta_2\big)
$$
by $\big([E],[E^\prime]\big) \mapsto [E \wedge E^\prime]$. A point of $\Gamma(\eta_1 \times_\xi \eta_2)$ can 
also be viewed as an extension~$F$ of $\ulcX^{(a)}_R$ by $\ulcF_2$, and similarly for $\Gamma(\eta_1 \times_\xi 
\eta_2^\prime)$. This interpretation gives rise to maps
$$
\psi(\eta_1;\eta_2,\eta_2^\prime) \colon \Gamma(\eta_1 \times_\xi \eta_2) \times \Gamma(\eta_1 \times_\xi 
\eta_2^\prime) \tto \Gamma\big(\eta_1 \times_\xi (\eta_2\wedge \eta_2^\prime)\big)
$$
by $\big([F],[F^\prime]\big) \mapsto [F \wedge F^\prime]$. These maps $\phi(\eta_1,\eta_1^\prime;\eta_2)$ and 
$\psi(\eta_1;\eta_2,\eta_2^\prime)$ give the desired structure of a biextension; we leave it to the reader to 
verify that they satisfy all required compatibilities.

\sssection{CascHDual}
{\it Remark.} --- If $\ul{X}^D$ is the Cartier dual of $\ul{X}$ then $\Def\big(\ul{X}^D\big)$ is naturally 
isomorphic, as a cascade, to the dual of $\Def(\ul{X})$.

\sssection{PhiCan}
Situation as in \refn{DefoSit}. Write~$\sigma$ for the Frobenius automorphism of~$W(K)$. For any $n\geq 1$ we have a canonical isomorphism of formal $W(K)$-schemes
$$
\Def\big(\ul{X}^{(\sigma^n)}\big) \isomarrow \Def(\ul{X})^{(\sigma^n)}\, .\leqno(\refn{PhiCan}.1)
$$
For simplicity, write $\mD := \Def(\ul{X})$ and $D := \mD \hat\otimes_{W(K)} K$. If $\ulcX$ is a deformation of~$\ul{X}$ over 
an algebra $R \in \catC_W$ of characteristic~$p$ then $(\Frob_R^n)^\ast \ulcX$ is a deformation of~$\ul{X}^{(\sigma^n)}$. 
Via~(\refn{PhiCan}.1) this gives a morphism
$$
\phi_n \colon D \tto D^{(\sigma^n)}\, .
$$
This morphism is none other than the $n$th power relative Frobenius of~$D$ over~$K$. 

Fix an algebraically closed field~$k$ containing~$K$. Define $m_0$ to be period of~$\gf$, i.e., the smallest positive integer such that $\gf(i+m_0) = \gf(i)$ for all $i\in\cI = \Hom(\kappa,k)$. (Cf.~\refn{modpRefl}.) We are 
going to define a lifting
$$
\Phi^\can\colon \mD \tto \mD^{(\sigma^{m_0})}
$$
of $\phi_{m_0}$. Recall that we have a decomposition $\ul{X} = \ul{X}^{(1)} \times \cdots \times \ul{X}^{(r)}$ 
with $\ul{X}^{(\nu)}$ isoclinic of slope~$\lambda_\nu$. Using that $\ul{X}$ is a $K$-form of the standard ordinary 
object of type $(d,\gf)$ we see that
$$
\eqalign{
\ul{X}\big[F^{m_0}\big] &= \ul{X}^{(1)}\big[F^{m_0}\big] \times \cdots \times \ul{X}^{(r)}\big[F^{m_0}\big]\cr
&= \ul{X}^{(1)}\big[p^{\lambda_1 m_0/m}\big] \times \cdots \times \ul{X}^{(r)}\big[p^{\lambda_r m_0/m}\big]\, .\cr}
$$
Note that $\lambda_\nu m_0/m \in \mZ$ for all~$\nu$. Now suppose $\ulcX$ is a lifting of~$\ul{X}$ over some $R\in 
\catC_{W(K)}$. As shown in~\refn{FiltLift} we have a slope filtration
$$
0 \subset \ulcX^{(r,r)} \subset \ulcX^{(r-1,r)} \subset \cdots \subset \ulcX^{(1,r)} = \ulcX\, .
$$
Define a finite subgroup scheme $Q \subset \ulcX$ by
$$
Q := \ulcX^{(r,r)}\big[p^{\lambda_r m_0/m}\big] + \ulcX^{(r-1,r)}\big[p^{\lambda_{r-1}m_0/m}\big] + \cdots + 
\ulcX^{(1,r)}\big[p^{\lambda_1 m_0/m}\big]\, .
$$
Note that $Q \otimes_R K = \ul{X}\big[F^{m_0}\big]$. One easily verifies that $\ulcX/Q$ is again a BT with 
$\cO$-structure over~$R$, which is a lifting of $\ul{X}^{(\sigma^{m_0})}$. This construction defines a functor 
$\Def(\ul{X}) \to \Def\big(\ul{X}^{(\sigma^{m_0})}\big)$, and by composition with (\refn{PhiCan}.1) we get a 
morphism $\Phi^\can\colon \mD \to \mD^{(\sigma^{m_0})}$ that lifts~$\phi_{m_0}$.

\sssection{PhiCanHom}
{\it Proposition. --- With respect to the cascade structure on $\mD := \Def(\ul{X})$ defined in\/ {\rm 
\refn{CascConstr}}, the morphism $\Phi^\can\colon \mD \to \mD^{(\sigma^{m_0})}$ is a homomorphism of cascades.}
\Dskip

The proof of the proposition is tedious but straightforward; we leave it to the reader. We do not know if one 
can {\it characterize\/} the cascade structure on $\Def(\ul{X})$ by its property that $\Phi^\can$ defines a 
homomorphism, as in the ``classical'' ordinary case---cf.\ the appendix by Katz to~[\DelCoord].

\sssection{HcanDef}
{\it Definition.} --- Situation as in \refn{DefoSit}. The {\it canonical lifting\/} $\ul{X}^\can$ of~$\ul{X}$ 
over~$W(K)$ is the lifting corresponding to the zero section of the cascade $\Def(\ul{X})$. Concretely, if 
$\ul{X} = \ul{X}^{(1)} \times \cdots \times \ul{X}^{(r)}$ is the slope decomposition of~$\ul{X}$ as 
in~\refn{HDecCor} then by \refn{1slopeRig} each isoclinic factor $\ul{X}^{(\nu)}$ has a unique lifting 
$\ulcX^{(\nu)} = \ulcX^{(\nu)}_{W(K)}$ over~$W(K)$, and
$$
\ul{X}^\can := \ulcX^{(1)} \times \cdots \times \ulcX^{(r)}\, .
$$

\sssection{EndoHtil}
Let $R$ be a complete local $W(K)$-algebra with residue field~$K$. If $\ulcX$ is a lifting of~$\ul{X}$ over~$R$ 
then the natural map $\End_R(\ulcX) \to \End_K(\ul{X})$ is injective; this follows from Illusie~[\Illus], d) of Thm.~4.4. 

Write $L$ for the fraction field of~$\cO$. Recall that $(d,\gf)$ is the type of~$\ul{X}$. Given a lifting 
$\ulcX$ over~$R$, write $\tilde R := R \hat\otimes_{W(K)} W(\overline K)$. We say that $\ulcX$ is {\it of 
CM-type\/} if $\End_{\tilde R}(\ulcX) \otimes_{\Zp} \Qp$ contains a commutative semi-simple $L$-subalgebra $\cE$ 
with $\dim_L(\cE) = d$.

\sssection{TorsPts}
{\it Proposition. --- {\rm (\romno1)}~The canonical lifting $\ul{X}^\can$ is the unique lifting of~$\ul{X}$ with the property that (geometrically) all endomorphisms lift. More precisely, suppose $K=k$ is algebraically closed. Let $\ulcX$ be a lifting of~$\ul{X}$ over~$R \in \catC_{W(k)}$. Then the map $\End_R(\ulcX) \tto \End_k(\ul{X})$ is an isomorphism if and only if $\ulcX \cong \ul{X}^\can \otimes_{W(k)} R$. 

{\rm (\romno2)}~Let $\mD := \Def(\ul{X})$ have the structure of an $r$-cascade defined in\/~{\rm 
\refn{CascConstr}}. Let $R$ be a complete local $W(K)$-algebra with residue field~$K$. Let $s \in \mD(R)$ 
correspond to a lifting $\ulcX$ of~$\ul{X}$ over~$R$. Then the following properties are equivalent:
\itemitem{\rm (a)}~$s$ is a torsion point;
\itemitem{\rm (b)}~$\ulcX$ is isogenous to~$\ul{X}^\can$;
\itemitem{\rm (c)}~$\ulcX$ is of CM-type.}
\Dskip

\Proof (\romno1) Write $\cO^\prime := \End(\ul{X})$. Then $\cO^\prime$ is a product of matrix algebras over 
finite unramified extensions of~$\Zp$; see~\refn{EndoLem2}. Write $\ul{X}^\prime$ for $\ul{X}$ viewed as a BT 
with $\cO^\prime$-structure. On the one hand, using the explicit description of the ordinary type over $k = 
\kbar$, it is clear that all endomorphisms of~$\ul{X}$ lift to endomorphisms of~$\ul{X}^\can$. On the other 
hand, it is not difficult to see that $\ul{X}^\prime$ is again ordinary, and that it is rigid.

(\romno2) To see that (a) $\Rightarrow$ (b) we argue by induction on~$r$, the number of slopes. We use the notation of~\refn{FiltLift}. If $r=1$ then $\ul{X}$ is rigid and there is nothing to prove. For $r\geq 2$ we have an extension
$$
0 \tto \ulcX^{(2,r)} \tto \ulcX \tto \ulcX^{(1)}_R \tto 0\, .\leqno(\refn{TorsPts}.1)
$$
By induction we may assume that $\ulcX{}^{(2,r)}$, which is a lifting of $\ul{X}{}^{(2,r)}$, is isogenous to the 
canonical lifting 
$$
\ul{X}^{(2,r),\can} \otimes_{W(K)} R = \ulcX^{(2)}_R \times \cdots \times \ulcX^{(r)}_R\, .
$$
But if the class of the extension (\refn{TorsPts}.1) is torsion then $\ulcX$ is isogenous to $\ulcX{}^{(2,r)} 
\times \ulcX^{(1)}_R$. Hence (a) implies~(b). 

That (b) implies~(c) is immediate. For (c) $\Rightarrow$~(a) we may assume that $K = k$ is algebraically closed. 
Let $L$ be the fraction field of~$\cO$. Suppose that $\ulcX$ is of CM-type, i.e., there is a commutative 
semi-simple $L$-subalgebra $\cE \subset \End^0(\ulcX)$ with $\dim_L(\cE) = d$. Again we are going to use 
induction on~$r$. For $r=1$ there is nothing to prove. For $r\geq 2$ it suffices to show that that the extension 
class of~(\refn{TorsPts}.1) is torsion and that $\ulcX{}^{(2,r)}$ is of CM-type, too.

As $\End_R(\ulcX) \otimes \Qp$ maps injectively to 
$$
\End_k(\ul{X}) \otimes \Qp = M_{d^1}(L) \times \cdots \times M_{d^r}(L)\, ,\leqno(\refn{TorsPts}.2)
$$
the algebra~$\cE$ is a product $\cE_1 \times \cdots \times \cE_r$ with $\cE_j$ a field extension of degree~$d^j$ 
of~$L$. If $O_j$ is the ring of integers in~$\cE_j$ then $\cE \cap \End_R(\ulcX)$ is a subring of finite index 
in $O_1 \times \cdots \times O_r$. Further, every $\alpha \in \End_R(\ulcX)$ maps $\ulcX{}^{(2,r)} \subset 
\ulcX$ into itself; indeed, the composition $\ulcX{}^{(2,r)} \hookrightarrow \ulcX \mapright{\alpha} \ulcX 
\twoheadrightarrow \ulcX^{(1)}_R$ is zero, as it is zero on the special fibre for slope reasons. Of course, 
under the decomposition~(\refn{TorsPts}.2) the resulting homomorphism $h\colon \End_R(\ulcX) \to 
\End_R\big(\ulcX{}^{(2,r)}\big)$ is given by the projection onto the last $(r-1)$ factors. The kernel of~$h$ 
maps injectively to $\End_R\big(\ulcX^{(1)}_R\big) = M_{d^1}(\cO)$. Combining these remarks we readily find 
that $\ulcX{}^{(2,r)}$ is of CM-type. Finally, because $\cE \cap \End_R(\ulcX)$ is of finite index in $O_1 
\times \cdots \times O_r$, there are non-zero integers~$n_j$ such that $(n_1,0,\ldots,0)$ and 
$(0,n_2,\ldots,n_r)$ are both in~$\End_R(\cX)$. This implies that the class of~(\refn{TorsPts}.1) is 
torsion.\QED

\section{Ordinary polarized Barsotti-Tate $\cO$-modules}{PolBT}

\ssection{Generalities on BT$_n$ with $(\cO,\ast,\epsilon)$-structure}{BTn+O*e} 

\sssection{BTO*eDef} 
{}From now on we assume that $p > 2$. If $X$ is a commutative finite locally free group scheme over some basis~$S$ then we write $X^D$ for its Cartier dual. If $X$ is a BT over~$S$ then we write $X^D$ for its Serre dual. In both cases there is a canonical isomorphism $\kappa_X\colon X \isomarrow X^{DD}$.

Let $n \in \mN \cup \{\infty\}$ and $\epsilon \in \{\pm 1\}$. If $X$ is a BT$_n$ over a base scheme~$S$ then by an {\it $\epsilon$-duality\/} of~$X$ we mean an isomorphism $\lambda\colon X \isomarrow X^D$ such that $\lambda = \epsilon \cdot \lambda^D \circ \kappa_X$. Such an $\epsilon$-duality induces an involution $f \mapsto f^\dagger$ on the ring $\End_S(X)$. We also refer to an $\epsilon$-duality as a polarization.  

Let $(\cO,\ast)$ be a $\Zp$-algebra equipped with a $\Zp$-linear involution $b \mapsto b^\ast$. Let $\epsilon \in \{\pm 
1\}$. By a BT$_n$ with $(\cO,\ast,\epsilon)$-structure over~$S$ we mean a triple $\ul{X} = (X,\iota,\lambda)$ where $(X,\iota)$ is a BT$_n$ with $\cO$-structure and $\lambda\colon X \to X^D$ is an $\epsilon$-duality, such that $\iota(b^\ast) = \iota(b)^\dagger$ for all $b \in \cO$.

Let $K$ be a perfect field, $\charact(K) = p$. Let~$\sigma$ be the Frobenius automorphism of~$W_n(K)$. Then a BT$_n$ with 
$(\cO,\ast,\epsilon)$-structure over~$K$ corresponds to a $5$-tuple $(M,F,V,\phi,\iota)$, where
\item{---}$M$ is a free $W_n(K)$-module of finite rank,
\item{---}$F\colon M \to M$ is a $\sigma$-linear endomorphism, 
\item{---}$V \colon M \to M$ is a $\sigma^{-1}$-linear endomorphism,
\item{---}$\phi\colon M \times M \rightarrow W_n(K)$ is a perfect, $\epsilon$-symmetric bilinear form, and
\item{---}$\iota\colon \cO \to \End(M,F,V)$ is a $\Zp$-linear homomorphism. 

\noindent
In addition to the relation $F \circ V = p \cdot \id_M = V \circ F$ we should have 
$$
\leqalignno{
\phi(Fm_1,m_2) = \sigma\big(\phi(m_1,Vm_2)\big)\qquad
&\hbox{{\rm for all~$m_1$,~$m_2 \in M$;}}&(\refn{BTO*eDef}.1)\cr
\phi(bm_1,m_2) = \phi(m_1,b^\ast m_2)\qquad
&\hbox{{\rm for all $b\in \cO$ and $m_1$, $m_2 \in M$.}}&\cr}
$$
We shall mainly use this in the cases $n=1$ and $n=\infty$.

\sssection{BasicCas}
We call a $\Qp$-algebra {\it unramified\/} if it is isomorphic to a product of matrix algebras over finite unramified field extensions of~$\Qp$. We are interested in BT$_n$ with $(\cO,\ast,\epsilon)$-structure, where $\cO$ is a maximal order in an unramified $\Qp$-algebra. By Morita equivalence (see e.g.\ [\Knus], Chap.~\Romno 1, \S~9), the study of such objects reduces to the following four special cases.

\typesitem{Case C:} $\cO \cong W(\kappa)$, with $\kappa$ a finite field, $\ast=\id$ and $\epsilon = -1$.

\typesitem{Case D:} $\cO \cong W(\kappa)$, with $\kappa$ a finite field, $\ast=\id$ and $\epsilon = +1$.

\typesitem{Case AU:} $\cO \cong W(\tilde\kappa)$, with $\tilde\kappa \cong \mF_{p^{2m}}$ a finite field of even degree over~$\Fp$, with $\ast = \sigma^m$ the unique non-trivial involution, and $\epsilon = +1$.

\typesitem{Case AL:} $\cO \cong W(\kappa) \times W(\kappa)$, with $\kappa$ a finite field, with $\ast$ given by $(x,y)^\ast = (y,x)$, and $\epsilon = +1$.

In case~AL every BT$_n$ with $(\cO,\ast,+1)$-structure is of the form $\ul{X} \cong \ul{X}_1 \times \ul{X}_1^D$, where $\ul{X}_1$ is a BT$_n$ with $W(\kappa)$-structure, and where the $+1$-duality of~$\ul{X}$ is given by switching the factors $\ul{X}_1$ and $\ul{X}_1^D$. This reduces case~AL to the study of BT with $\cO$-structure.

\sssection{GPENot}
Let us now briefly review the second classification theorem proved in~[\GSAS]; this concerns a variant of Thm.~\refn{GEThm} for polarized objects. We shall state the result in its general form, not only for the basic cases C, D and~A.

Let $(B,\ast)$ be a finite dimensional semi-simple $\Fp$-algebra equipped with an involution $b \mapsto b^\ast$. Let $\epsilon \in \{\pm 1\}$. Let $\tilde\kappa$ be the center of~$B$ and $\kappa := \{z \in \tilde\kappa \mid z^\ast = z\}$. We can decompose $(B,\ast)$ as a product of simple factors, say $(B,\ast) = \prod_{n=1}^l (B_n,\ast_n)$. Accordingly we have decompositions $\tilde\kappa = \prod \tilde\kappa_n$ and $\kappa = \prod \kappa_n$. The $\kappa_n$ are finite fields. We have $B_n \cong M_{r_n}(\tilde\kappa_n)$ for some $r_n \geq 1$.

If $\ast_n$ is an involution of the second kind then either $\tilde\kappa_n \cong \kappa_n \times \kappa_n$ or $\tilde\kappa_n$ is a quadratic field extension of~$\kappa_n$. We say in this case that $(B_n,\ast_n)$ is of type~A. Next suppose $\ast_n$ is of the first kind; in this case $\tilde\kappa_n = \kappa_n$. Set $\epsilon_n = +1$ if $\ast_n$ is orthogonal, $\epsilon_n = -1$ if $\ast_n$ is symplectic. We say that $(B_n,\ast_n)$ is of type~C if $\epsilon \cdot \epsilon_n = -1$ and that it is of type~D if $\epsilon \cdot \epsilon_n = +1$.

Let $\cI = \cI_1 \cup \cdots \cup \cI_l$ be the set of homomorphisms $\kappa \rightarrow k$. For $X \in \{{\rm C},{\rm D},{\rm A}\}$, let $\cI^X \subset \cI$ be the union of all subsets $\cI_n \subset \cI$ for which $(B_n,\ast_n)$ is of type~$X$. Let~$\cItil = \cItil_1 \cup \cdots \cup \cItil_l$ be the set of homomorphisms $\tilde\kappa \rightarrow k$. We have a restriction map $\res\colon \cItil \rightarrow \cI$. For $\tau \in\cItil$ define $\bar\tau := \tau \circ \ast$. If $i\in\cI^{{\rm C}} \cup \cI^{{\rm D}}$ there is a unique $\tau \in \cItil$ with $\res(\tau) = i$, and $\tau = \bar\tau$; if $i \in \cI^{{\rm A}}$ there are precisely two elements $\tau$, $\bar\tau \in \cItil$ that restrict to the embedding $i$ on~$\kappa$. 

\sssection{pairs2}
Let $k$ be an algebraically closed field, $\charact(k) = p>2$. Consider triples $(N,L,\psi)$ consisting of a finitely generated $B \otimes_\Fp k$-module~$N$, an $\epsilon$-$\ast$-hermitian form $\psi\colon N \times N \to B \otimes_{\Fp} k$, and a maximal isotropic submodule $L \subset N$. With a similar construction as in~\refn{pairs}, such a triple is classified, up to isomorphism, by a pair $(d,\gf)$ consisting of functions $d\colon \cI \to \mZ_{\geq 0}$ and $\gf\colon \cItil \to \mZ_{\geq 0}$ such that $\gf(\tau) + \gf(\bar\tau) = d(i)$ for all $\tau\in\cItil$ and $i = \res(\tau) \in \cI$.

\sssection{dfdelta}
Notation as above. Let $\ul{Y}$ be a BT$_1$ with $(B,\ast,\epsilon)$-structure over~$k$. To~$\ul{Y}$ we associate a triple of invariants $(d,\gf,\delta)$, referred to as its type.

Let $\ul{N} = (N,F,V,\iota,\phi)$ be the Dieudonn\'e module of~$\ul{Y}$. There is a unique $\epsilon$-$\ast$-hermitian form $\psi\colon N \times N \to B \otimes_{\Fp} k$ such that $\phi = \Trd \circ \psi$, where $\Trd\colon B \otimes_\Fp k \to k$ is the reduced trace. Set $L := \Ker(F)$. Let $(d,\gf)$ be the pair of functions corresponding to $(N,L,\psi)$. It can be shown ([\GSAS], 4.3 and~6.5) that the function~$d$ is constant on each of the subsets $\cI_n \subset \cI$. Note that for $i \in \cI^{{\rm C}} \cup \cI^{{\rm D}}$ there is a unique $\tau = \bar\tau$ with $\res(\tau) = i$, hence $d(i) = 2\cdot \gf(\tau)$.

Finally we define a function $\delta \colon \cI^{{\rm D}} \to \mZ/2\mZ$. Given $i \in \cI^{{\rm D}}$, let $\tau \in \cItil$ be the unique element with $\res(\tau) = i$, and write $N_i := N_\tau \subset N$. Then let
$$
\delta(i) = \len_{B \otimes_\Fp k}\Bigl(\Ker\big(F_{|N_i}\big)\Bigm/ \Ker\big(F_{|N_i}\big) \cap \Ker\big(V_{N_i}\big)\Bigr) \bmod{2}\, .
$$
If there are no factors of type~D then $\cI^{{\rm D}} = \emptyset$ and the invariant~$\delta$ is void.

\sssection{wdefGPE}
Fix a triple $(N_0,L_0,\psi_0)$ as in~\refn{pairs2}, corresponding to a pair~$(d,\gf)$ with $d$ constant on each subset $\cI_n \subset \cI$. Define $G := \Sp_{B \otimes_\Fp k}(N_0,\psi_0)$, the algebraic group (over~$k$) of $B \otimes_\Fp k$-linear automorphisms of~$N_0$ that preserve the form~$\psi_0$. We have $G = \prod_{i\in\cI} G_i$, with $G_i$ isomorphic to $\Sp_{d(i),k}$ if $i\in \cI^{{\rm C}}$, to $\OO_{d(i),k}$ if $i\in \cI^{{\rm D}}$ and to $\GL_{d(i),k}$ if $i\in \cI^{{\rm A}}$.

Let $G^0 \subset G$ be the identity component. Define $\mX^0$ to be the conjugacy class of parabolic subgroups of~$G^0$ containing $\Stab(L_0)$. Write $W_{G^0}$ for the Weyl group of~$G^0$, and let $W_{\mX^0} \subset W_{G^0}$ be the subgroup corresponding to~$\mX^0$.

Let $\ul{Y}$ be a BT$_1$ with $(B,\ast,\epsilon)$-structure over~$k$, of type $(d,\gf,\delta)$. To~$\ul{Y}$ we associate an 
element $w(\ul{Y}) \in W_{\mX^0} \backslash W_{G^0}$. This works essentially the same as in~\refn{wdefGE}: Choose an isometry $\xi\colon (N,\psi) \isomarrow (N_0,\psi_0)$ that restricts to $L \isomarrow L_0$. Then we choose a Borel subgroup $Q \subset G^0$ that stabilizes the canonical filtration $\calC_\gdot$ (viewed as a flag in~$N_0$ via~$\xi$), and we define $w(\ul{Y})$ to be the Weyl group coset measuring the relative position of~$\Stab(L_0)$ and~$Q$. This is independent of the choices of~$\xi$ and~$Q$.

With these notations the second main theorem of~[\GSAS] is the following.

\sssection{GPEThm}
{\it Theorem. --- Let $k$ be an algebraically closed field, $\charact(k) > 2$. Sending a BT$_1$ with $(B,\ast,\epsilon)$-structure $\ul{Y}$ to the element $w(\ul{Y})$ gives a bijection
$$
\left\{
\vcenter{
\setbox0=\hbox{{\rm isomorphism classes of}}
\copy0
\hbox to \wd0{{\rm \hfil $\ul{Y}$ of type $(d,\gf,\delta)$\hfil}}}
\right\}
\longisomarrow
W_{\mX^0} \backslash W_{G^0}\, .
$$\vskip-\lastskip\smallskip}

\sssection{GG0Rem}
{\it Remark.\/} --- In~[\GSAS] we have given two versions of the above theorem: the result as stated here, and a version working with the possibly non-connected group~$G$. (This is only relevant if there are factors of type~D.) In this paper we shall exclusively work with the connected group~$G^0$. The notation $G^0$ and $\mX^0$ should remind us of this.  

\ssection{Ordinary BT with $(\cO,\ast,\epsilon)$-structure}{OrdBT+O*e}

\sssection{OrdPolSit}
{\it Situation.} --- We assume $p > 2$. Let $\cB$ be an unramified semi-simple $\Qp$-algebra, equipped with an involution~$\ast$. Let $\cO \subset \cB$ be a maximal order that is stable under~$\ast$. Write $B := \cO/p\cO$, which is a finite dimensional semi-simple $\Fp$-algebra on which we have an induced involution~$\ast$. Let $k = \kbar$, $\charact(k) = p$. Let $\epsilon \in \{\pm 1\}$.

Let $\ul{X} = (X,\iota,\lambda)$ be a BT with $(\cO,\ast,\epsilon)$-structure over~$k$. Write $\ul{Y} := \ul{X}[p]$, which is a BT$_1$ with $(B,\ast,\epsilon)$-structure. Let $(d,\gf,\delta)$ be its type. Let $(G^0,\mX^0)$ be as in~\refn{wdefGPE}.
 
We should like to have a notion of ordinariness for the polarized object~$\ul{X}$. We shall take the same approach as in the non-polarized case. Thus, we define a notion of $[p]$-ordinariness, depending only on the structure of the $p$-kernel, and a notion of $\mu$-ordinariness, depending only on the isogeny class of~$\ul{X}$. Our main goal is then to show that the two notions are equivalent, and that, working over $k=\kbar$ and fixing $(d,\gf,\delta)$, there is a unique ordinary object, up to isomorphism. For factors of type~C or~A, most of this is a rather straightforward extension of the results in the non-polarized case. The factors of type~D require some extra work.

\sssection{[p]OrdPol}
{\it Definition.} --- Situation as in~\refn{OrdPolSit}. Let $w^\ord \in W_{\mX^0}\backslash W_{G^0}$ be the class of the longest element of~$W_{G^0}$. We say that $\ul{X}$, as a BT with $(\cO,\ast,\epsilon)$-structure, is {\it $[p]$-ordinary\/} if $w(\ul{Y}) = w^\ord$.

\sssection{XordPol}
We define a $[p]$-ordinary object $\ul{X}^\ord = \ul{X}^\ord(d,\gf,\delta)$ over~$k$. We shall only do this in the basic cases C, D and~AU. As explained in~\refn{BasicCas}, case~AL reduces to the study of BT with $\cO$-structure (without polarization), and $\ul{X}^\ord$ corresponds to the standard ordinary object described in~\refn{XordNPol}. In the general case we can define $\ul{X}^\ord$ by ``reversing'' the reduction step discussed in~\refn{BasicCas}, based on Morita equivalence; we leave the details of this to the reader. In the cases C and~AU the invariant~$\delta$ plays no role, and we simply omit it in the discussion.

{\it Case~C.\/} In this case the pair $(d,\gf)$ has a very simple form: there is a natural number~$q$ such that $d=2q$ and $\gf(i) = q$ for all $i \in \cI$. Let $\ul{X}_\et$ and $\ul{X}_\mult$ be as in~\refn{EtMult}. Then $\ul{X}_\et \times \ul{X}_\mult$ has a natural $-1$-duality, and $\ul{X}^\ord = (\ul{X}_\et \times \ul{X}_\mult)^q$.

{\it Case~D.\/} The pair $(d,\gf)$ is as in case~C: $d=2q$ and $\gf(i) = q$ for all $i \in \cI$. Further, $\delta$ is an arbitrary function $\cI \to \mZ/2\mZ$. 

First we do the case $q=1$. Up to isomorphism there is a unique BT$_1$ with $(\kappa,\id,+1)$-structure of type $(2,1,\delta)$; we call it~$\ul{Y}(\delta)$. The Dieudonn\'e module of the corresponding standard ordinary object $\ul{X}(\delta) = \ul{X}^\ord(2,1,\delta)$ is given as follows. Let $M$ be the free $W(k)$-module with basis $\{e_{i,j}\}$ for $i \in \cI$ and $j \in \{1,2\}$. Let $b \in \cO$ act on $e_{i,j}$ as multiplication by~$i(b) \in W(k)$. Frobenius is given on base vectors by
$$
\left\{\matrix{
F(e_{i,1}) = e_{i+1,1}\hfill\cr
F(e_{i,2}) = p\cdot e_{i+1,2}\cr}
\right.\quad \hbox{if $\delta(i) = \bar 1$;}
\qquad\quad
\left\{\matrix{
F(e_{i,1}) = p\cdot e_{i+1,2}\cr
F(e_{i,2}) = e_{i+1,1}\hfill\cr}
\right.\quad \hbox{if $\delta(i) = \bar 0$.}
$$
Verschiebung is determined by the rule that $FV = p = VF$. The form~$\phi$ is an orthogonal sum of the forms~$\phi_i$ on $M_i = \Span(e_{i,1},e_{i,2})$ given by the matrix $0\, 1\choose1\, 0$.

For $q>1$ we have $\ul{X}^\ord = (\ul{X}_\et \times \ul{X}_\mult)^{q-1} \times \ul{X}(\delta)$, where this time we equip $\ul{X}_\et \times \ul{X}_\mult$ with its natural $+1$-duality. Note that if $\delta$ is the constant function~$\bar 1$ then $\ul{X}(\delta) = \ul{X}_\et \times \ul{X}_\mult$, so in this case we have $\ul{X}^\ord = (\ul{X}_\et \times \ul{X}_\mult)^q$.

{\it Case~AU.\/} We have $d \in \mN$, and $\gf\colon \cItil \to \mZ_{\geq 0}$ is a function with $\gf(\tau) + \gf(\bar\tau) = d$ for all $\tau\in\cItil$. Since $\tilde\kappa$ is a finite field, $\cItil$ is a finite set with cyclic ordering. Let $M$ be the free $W(k)$-module with basis $\{e_{\tau,j}\}$ for $\tau \in \cItil$ and $j \in \{1,\ldots,d\}$. Define $F$ and~$V$ by
$$
F(e_{\tau,j}) = \cases{e_{\tau+1,j}&if $j \leq d-\gf(\tau)$;\cr
p \cdot e_{\tau+1,j}&if $j>d-\gf(\tau)$;}
\qquad
V(e_{\tau+1,j}) = \cases{p \cdot e_{\tau,j}&if $j \leq d-\gf(\tau)$;\cr
e_{\tau,j}&if $j>d-\gf(\tau)$.}
$$
The pairing~$\phi$ can be chosen in such a way that $\phi(e_{\tau,j},e_{\tau^\prime,j^\prime}) \neq 0$ only if $\tau^\prime = \bar\tau$ and $j=j^\prime$, and such that $\phi(e_{\tau,j},e_{\bar\tau,j}) =: c_\tau$ only depends on~$\tau$. In order for this pairing to satisfy~(\refn{BTO*eDef}.1) we should then choose the function $\tau \mapsto c_\tau$ such that $c_{\tau + 1} = \sigma(c_\tau)$ for all~$\tau$. In particular, if $E \subset k$ is the subfield with $p^{2m} = \#\tilde \kappa$ elements then $c_\tau \in W(E)^\times$ for all~$\tau$. The choice of the constants~$c_\tau$ is not unique, but it can be shown that, up to isomorphism, the resulting Dieudonn\'e module~$\ul{M}$ is independent of this choice.

\sssection{muOrdPol}
Our next objective is to define the notion of $\mu$-ordinariness for BT with $(\cO,\ast,\epsilon)$-structure, analogous to the definition in~\refn{muordDef}. In the polarized case we cannot give the definition in terms of a single Newton polygon; instead we have to go deeper into the theory developed in [\KottIsoc] and~[\RR]. We closely follow Wedhorn~[\WedhOrd], to which the reader is referred for more details. For simplicity of exposition we shall assume that we are in one of the basic cases C, D or~AU.

Let $\cV := \cB^d$, with its natural structure of a left $\cB$-module. Let $\psi\colon \cV \times \cV \to \cB$ be an $\epsilon$-$\ast$-hermitian form. Write $\gamma \mapsto \bar \gamma$ for the associated involution of the $\Qp$-algebra $\End_{\cB}(\cV)$, and let $\cG = \CSp_\cB(\cV,\psi)$ be the algebraic group over~$\Qp$ given, as a functor on $\Qp$-algebras, by
$$
\cG(A) = \big\{\gamma \in \End_{\cB}(\cV) \otimes_\Qp A \bigm| \gamma \bar \gamma \in A^\times\big\}\, .
$$
Let $(X^\ast,R^\ast,X_\ast,R_\ast,\Delta)$ be the based root datum of~$\cG$. We have a natural action of $\Gamma := \Gal(\Qpbar/\Qp)$ on~$X_\ast$. Let $W_{\cG^0}$ be the Weyl group of~$\cG^0$ (= the Weyl group of the root datum). The closed Weyl chamber $\overline{C} \subset (X_\ast \otimes \mQ)$ corresponding to the root base~$\Delta$ is stable under the action of~$\Gamma$ and is a fundamental domain for the action of~$W_{\cG^0}$.

We define a subset $\Ord(d,\gf) \subset (X_\ast \otimes \mQ)/W_{\cG^0}$ of ordinary points. Choose a decomposition $(\cV \otimes \Qpbar) = \cW_0 \oplus \cW_1$, where $\cW_0$ and~$\cW_1$ are $(\cB \otimes \Qpbar)$-submodules, totally isotropic with respect to~$\psi$, with $\cW_1$ of type~$\gf$. Define a cocharacter $\mu\colon \mG_{m,\Qpbar} \tto \cG_{\Qpbar}$ by the requirement that $\mu(z)$ acts as multiplication by $z^j$ on $\cW_j$. The set~$\gc$ of all cocharacters obtained in this way is a union of $\cG^0(\Qpbar)$-conjugacy classes, say $\gc = \gc_1 \cup \cdots \cup \gc_s$. (Of course, $s > 1$ occurs only if $(\cO,\ast,\epsilon)$ is of type~D.) Let $\tilde\mu_j$ be the unique representative of~$\gc_j$ in~$\overline{C}$. If $\Gamma^\prime \subset \Gamma$ is the stabilizer of~$\tilde\mu_j$ in the Galois group then we define $\bar\mu_j \in (X_\ast \otimes \mQ)/W_{\cG^0}$ to be the image of the ``averaged'' element
$$
{1\over[\Gamma:\Gamma^\prime]}\, \sum_{\gamma \in \Gamma/\Gamma^\prime} \gamma \cdot \tilde\mu_j\, .
$$
Finally, define 
$$
\Ord(d,\gf) := \{\bar\mu_1,\ldots,\bar\mu_s\} \subset (X_\ast \otimes \mQ)/W_{\cG^0}
$$
to be the set of classes $\bar\mu_j$ thus obtained. As the notation suggests, $\Ord(d,\gf)$ takes the role of the Newton polygon $\Ord(d,\gf)$ defined in~\refn{ajDef}. If $(\cO,\ast,\epsilon)$ is of type~C or type~A then $\Ord(d,\gf)$ consists of a single element; in case~D we may get a set of more than~$1$ element. See Wedhorn~[\WedhOrd], section~2.3, for an explicit calculation of the set~$\Ord(d,\gf)$.

\sssection{NewtPt}
To $\ul{X}$ we can associate a {\it Newton point\/} $\bar\nu(\ul{X}) \in (X_\ast \otimes \mQ)/W_{\cG^0}$. For the definition we refer to~[\RR]. (One needs to combine loc.\ cit.\ (1.8), (3.4) and~(3.5).) The Newton point takes the role of the Newton polygon in the classical theory. Note however, that in general $\bar\nu(\ul{X})$ does not determine $\ul{X}$ up to isogeny, as the Newton map need not be injective.   

\sssection{muordPolDef}
{\it Definition.} --- Situation as in~\refn{OrdPolSit}. We say that $\ul{X}$, as a BT with $(\cO,\ast,\epsilon)$-structure, is {\it $\mu$-ordinary\/} if $\bar\nu(\ul{X}) \in \Ord(d,\gf)$.

\sssection{OrdPolThm}
{\it Theorem. --- Situation as in\/ {\rm \refn{OrdPolSit}}. Then the following are equivalent:
\item{\rm (a)} $\ul{X}$ is $\mu$-ordinary;
\item{\rm (b)} $\ul{X}$ is $[p]$-ordinary;
\item{\rm (c)} $\ul{X} \cong \ul{X}^\ord(d,\gf,\delta)$.

\noindent
If there are no factors of type~D or if the function $\delta\colon \cI^{{\rm D}} \to \mZ/2\mZ$ is the constant function~$\bar 1$ then\/ {\rm (a)--(c)} are equivalent to the condition that $(X,\iota)$, the underlying BT with $\cO$-structure, is ordinary in the sense of section\/ {\rm \refn{OrdBT+O}}.}
\Cskip

We divide the proof into a couple of steps.

\sssection{OPTStep1}
As usual we can reduce to the basic cases C, D, AU and~AL. In case~AL there is a further reduction to a statement about non-polarized BT, and the result follows from Thm.~\refn{ord=ord}.

We sketch the argument if we are in one of the cases C or~A, or if $\delta\colon \cI^{{\rm D}} \to \mZ/2\mZ$ is the constant function~$\bar 1$. Write $\ul{X}^\prime$ for the underlying BT with $\cO$-structure, without polarization. Similar notation for $\ul{Y} := \ul{X}[p]$. If $\ul{X}$ is $[p]$-ordinary then by inspection of~\refn{XordNPol} and~\refn{XordPol} we see that $\ul{X}^\prime$ is $[p]$-ordinary too. Conversely, suppose $\ul{X}^\prime$ is $[p]$-ordinary. Up to isomorphism there is a unique polarization on~$\ul{Y}^\prime$ that makes it into a BT$_1$ with $(B,\ast,\epsilon)$-structure; see~[\GSAS], 5.5 and~6.7. Hence $\ul{X}$ is $[p]$-ordinary. In particular this proves the last assertion of the theorem.

Suppose $\ul{X}$ is $\mu$-ordinary. It can be shown that then also $\ul{X}^\prime$ is $\mu$-ordinary. By Thm.~\refn{ord=ord} this implies that $\ul{X}^\prime$ is $[p]$-ordinary, and by the above it follows that $\ul{X}$ is $[p]$-ordinary. 

If $\ul{X}$ is $[p]$-ordinary then by Thm.~\refn{ord=ord} and the above we know that $\ul{X}^\prime \cong \ul{X}^{\prime,\ord}(d,\gf)$. We claim that up to isomorphism there is a unique polarization form on~$\ul{X}^{\prime,\ord}(d,\gf)$ making it into a BT with $(\cO,\ast,\epsilon)$-structure. In the cases~C and~D (still assuming that $\delta$ is the constant function~$\bar 1$) we have $\ul{X}^{\prime,\ord}(d,\gf) \cong (\ul{X}_\et \times \ul{X}_\mult)^q$ for some~$q$, and the claim follows without difficulty. In case~AU we may assume that $\ul{X}^{\prime,\ord}(d,\gf)$ is isoclinic (one slope), which means that it is isomorphic to the $d$-fold product of a height~$1$ object. The polarization forms then correspond to the isometry classes of rank~$d$ hermitian forms over~$W(\tilde\kappa)$. But there is only one such class, by [\Knus], Chap.~\Romno 2, (4.6.5) and the fact that the norm map $W(\tilde\kappa)^\times \to W(\kappa)^\times$ is surjective. Our claim follows.  

The implication (c) $\Rightarrow$~(a) follows by direct computation of the Newton point of~$\ul{X}^\ord(d,\gf,\delta)$.

\sssection{OPTStep2}
Let us now assume that we are in case~D and $\delta$ is not the constant function~$\bar 1$. Recall that $\cO = W(\kappa)$ for some finite field~$\kappa \cong \mF_{p^m}$ and that the type $(d,\gf)$ is given by $d=2q$ and $\gf(i) = q$ for all~$i$. The implication (c) $\Rightarrow$~(a) is again done by direct computation of the Newton point. Next suppose $\ul{X}$ is $\mu$-ordinary. To prove that $\ul{X}$ is $[p]$-ordinary it suffices to show that~$X$, the underlying~BT, has $p$-rank $\geq m\cdot (q-1)$; the point is that $\ul{Y}^\ord(d,\gf)$ is the {\it only\/} BT$_1$ with $(\kappa,\id,+1)$-structure which is of type $(2q,q,\delta)$ and for which the $p$-rank is $\geq m\cdot(q-1)$. We use the notation of~\refn{muOrdPol}, applied to case~D. (In particular, $\cB$ is the fraction field of~$W(\kappa)$.) Let $\cG^\prime := \GL_\cB(\cV)$, write $X_\ast^\prime$ for its coroot lattice and $W_{\cG^\prime}$ for its Weyl group. The inclusion $\cG \hookrightarrow \cG^\prime$ gives rise to a map $\beta\colon (X_\ast \otimes \mQ)/W_{\cG^0} \to (X_\ast^\prime \otimes \mQ)/W_{\cG^\prime}$. If $\ul{X} = (X,\iota,\lambda)$ is a BT with $\big(W(\kappa),\id,+1\big)$-structure and $\ul{X}^\prime = (X,\iota)$ is the underlying non-polarized BT with $W(\kappa)$-structure then $\beta$ maps the Newton point of~$\ul{X}$ to that of~$\ul{X}^\prime$. But the Newton point of $\ul{X}^\prime$ can be represented by a single Newton polygon (cf.\ the proof of~\refn{ulHIsog}), and its $p$-rank is simply $m$~times the multiplicity of the slope~$0$ in that polygon. Hence everything boils down to verification that in each of the polygons $\beta(\bar\mu)$, for $\bar\mu \in \Ord(d,\gf) \subset (X_\ast \otimes \mQ)/W_{\cG^0}$, the slope~$0$ has multiplicity~$\geq q-1$. This easily follows from the computations by Wedhorn in [\WedhOrd],~(2.3.4).

Finally, assume $\ul{X}$ is $[p]$-ordinary. With the notation of~\refn{XordPol}, Case~D, we have $\ul{Y} \cong (\ul{Y}_\et \times \ul{Y}_\mult)^{q-1} \times \ul{Y}(\delta)$. By $p$-rank considerations we have a similar decomposition for~$\ul{X}$, say $\ul{X} \cong (\ul{X}_\et \times \ul{X}_\mult)^{q-1} \times \ul{X}^{(2)}$, where $\ul{X}^{(2)}$ is BT with $\big(\cO,\id,+1\big)$-structure of type $(2,1,\delta)$. Hence to prove that $\ul{X} \cong \ul{X}^\ord(d,\gf,\delta)$ we may assume that~$q=1$. As usual we write $\ul{M}$ for the Dieudonn\'e module of~$\ul{X}$. We have a natural decomposition $M = \oplus_{i\in\cI} M_i$; write $\phi_i$ for the restriction of~$\phi$ to~$M_i$. Let $\ul{N} = \ul{M}/p\ul{M}$ and write $\ol\phi_i = \phi_i \bmod p$. As $\ul{Y} \cong \ul{Y}(\delta)$ we can choose a basis $\{e_{i,j}\}$ for~$N$ (with $i\in\cI$ and $j \in \{1,2\}$) such that $b \in \kappa$ acts on $N_i = k \cdot e_{i,1} + k \cdot e_{i,2}$ as multiplication by $i(b) \in k$, such that
$$
\left\{\matrix{
F(e_{i,1}) = e_{i+1,1}\cr
V(e_{i+1,2}) = e_{i,2}\cr}
\right.\quad \hbox{if $\delta(i) = \bar 1$,}
\qquad
\left\{\matrix{
F(e_{i,2}) = e_{i+1,1}\cr
V(e_{i+1,2}) = e_{i,1}\cr}
\right.\quad \hbox{if $\delta(i) = \bar 0$,}
$$
and such that the form $\ol\phi_i$ on~$N_i$ is given by the matrix~$0\, 1\choose1\, 0$.

We claim that for every $i\in\cI$ there exists an orthonormal basis $\{\tilde e_{i,1},\tilde e_{i,2}\}$ for~$M_i$ such that $\tilde e_{i,j}$ reduces to $e_{i,j}$ modulo~$p$. Further, this lifted basis is unique up to a scalar: any other such basis is of the form $\{c\tilde e_{i,1}, c^{-1}\tilde e_{i,2}\}$ with $c \in 1 + pW(k) \subset W(k)^\times$. To prove the claim, let $e_{i,j} \in N_i$ be the vector generating the Frobenius kernel (i.e, $j=1$ if $\delta(i) = \bar 0$ and $j=2$ if $\delta(i) = \bar 1$). Clearly it suffices to show that $e_{i,j}$ can be lifted to a vector $\tilde e_{i,j} \in M_i$ such that $\phi_i\big(\tilde e_{i,j},\tilde e_{i,j}\big) = 0$, and that this lifting is uniquely determined up to a scalar in $1 + pW(k)$. Start with an arbitrary $u \in M_i$ reducing to $e_{i,j}$ modulo~$p$. Choose any~$v$ such that $\{u,v\}$ is a $W(k)$-basis for~$M_i$. As $\phi_i(u,u) \equiv 0 \bmod p$ we have $\phi_i(u,v) \in W(k)^\times$, so after rescaling the vector~$v$ we can assume that $\phi_i(u,v) = 1$. Let $\gamma = -\phi_i(u,u)/2$ and set $u^\prime := u + \gamma v$. Note that $\gamma \equiv 0 \bmod p$, as $p\neq 2$; hence $u^\prime$ lifts $e_{i,j}$ and $\{u^\prime,v\}$ is again a basis for~$M_i$. Finally, $\phi_i(u^\prime,u^\prime) = \phi_i(u,u)^2 \phi_i(v,v)/4$, which is $p$-adically closer to~$0$ than~$\phi_i(u,u)$. As $M_i$ is $p$-adically complete, the existence of the desired lifting $\tilde e_{i,j}$ follows by approximation. That this lifting is unique up to a scalar is straightforward to check, again using that~$p\neq 2$.

The rest of the argument is easy. Choose a starting point $i_0 \in \cI$. As just shown we can choose an orthonormal basis $\{\tilde e_{i_0,1},\tilde e_{i_0,2}\}$ for~$M_{i_0}$. Let $j_0 \in \{1,2\}$ be the index such that $F(\tilde e_{i_0,j_0}) \equiv 0 \bmod p$; let $l_0$ be the other index. Note that there is a unique vector in $M_{i+1}$ which maps to $\tilde e_{i_0,j_0}$ under~$V$; hence we can define $\tilde e_{i_1,1} := F(\tilde e_{i_0,l_0})$ and $\tilde e_{i_1,2} := V^{-1}(\tilde e_{i_0,j_0})$. It readily follows from (\refn{BTO*eDef}.1) that $\{\tilde e_{i_1,1},\tilde e_{i_1,2}\}$ is an orthonormal basis for~$M_{i_1}$. Iterating this construction we arrive, after $m$ steps, at a second orthonormal basis $\{\tilde e_{i_m,1},\tilde e_{i_m,2}\}$ for $M_{i_m} = M_{i_0}$. As shown above, this second basis differs from the first one by a scalar $c \in 1+pW(k)$. But if we rescale $\{\tilde e_{i_0,1},\tilde e_{i_0,2}\}$ by a factor $\gamma \in 1+pW(k)$ then this affects the resulting basis $\{\tilde e_{i_m,1},\tilde e_{i_m,2}\}$ by a factor $\sigma^m(\gamma)$. Choosing~$\gamma$ such that $\gamma = \sigma^m(\gamma) \cdot c$ (such a~$\gamma$ exists!) we have brought the Dieudonn\'e module $\ul{M}$ in standard form. This completes the proof of~\refn{OrdPolThm}. \QED

\sssection{PolOrdDef}
{\it Definition.\/} --- Let $K$ be a field of characteristic~$p$. Let $(\cO,\ast,\epsilon)$ be as in~\refn{OrdPolSit}. If $\ul{X}$ is a BT with $(\cO,\ast,\epsilon)$-structure over~$K$ then we say that $\ul{X}$ is {\it ordinary\/} if $\ul{X} \otimes_K k$ satisfies the equivalent conditions of~\refn{OrdPolThm} for some (equivalently: every) algebraically closed field~$k$ containing~$K$.

\ssection{Deformation theory of ordinary polarized objects}{DefoOrdPol}

\sssection{PolDefo}
To finish this section we describe the deformation theory in the polarized case. Let $\ul{X} = (X,\iota,\lambda)$ be an ordinary BT with $(\cO,\ast,\epsilon)$-structure over a perfect field~$K$ of characteristic~$p$. Let us first assume that we are in one of the basic cases C or~AU; case~D shall be discussed in~\refn{PolDefD} below. (As always, case~AL reduces to a study of non-polarized BT with given endomorphisms and requires no further explanation.) Write $\ul{X}^\prime = (X,\iota)$ for the underlying non-polarized BT with$\cO$-structure, and let $\mD := \Def\big(\ul{X}^\prime\big)$. The given duality $\lambda\colon \ul{X}^\prime \isomarrow \ul{X}^{\prime,D}$ induces an isomorphism of cascades
$$
\gamma\colon \mD \isomarrow \Def\big(\ul{X}^{\prime,D}\big) \quad{\buildrel\hbox{\eightpoint 
\refn{CascHDual}}\over\cong}\quad \mD^\vee\, .
$$
Note that $\mD^\vee$ has the same underlying space as~$\mD$; the duality ``~${}^\vee$~'' only involves the cascade structure. Hence we can define a formal subscheme $\mD^\lambda \subset \mD$ by $\mD^\lambda(R) := \{x \in \mD(R) \mid \gamma(x) = x\}$. Then we find that we have a natural isomorphism
$$
\Def\big(\ul{X}\big) \isomarrow \mD^\lambda\, ,
$$
and that $\mD^\lambda$ is stable under the Frobenius lifting $\Phi^\can$ defined in~\refn{PhiCan}.

\sssection{PolDefD}
Assume now we are in case~D. Let $(2q,q,\delta)$ be the type of~$\ul{X}$. If $\delta$ is the constant function~$\bar 1$ then~$X$, the underlying BT, is ordinary in the classical sense and the structure of~$\Def(\ul{X})$ is fully explained by the classical theory. {}From now on we therefore restrict our attention to the case that $\delta(i) = \bar 0$ for some~$i$. Note that in this case the underlying object~$\ul{X}^\prime$ is {\it not\/} ordinary in the sense of~\refn{BTordDef}, so we cannot directly use the theory developed in~\S~\refn{DefoTh}. With some easy modifications we have a similar theory in this case, though. We outline the main features. 

(a) To begin with, recall that the ordinary object $\ul{X}$ of type $(2q,q,\delta)$ canonically decomposes as $\ul{X} = \ul{X}^{(1)} \times \ul{X}^{(2)} \times \ul{X}^{(3)}$ with $\ul{X}^{(1)} \cong \ul{X}_\et^{q-1}$ and $\ul{X}^{(3)} \cong \ul{X}_\mult^{q-1}$, and with $\ul{X}^{(2)}$ ordinary of type $(2,1,\delta)$. It is important to note that $\ul{X}^{(2)}$ is again a polarized object, whereas the factors $\ul{X}^{(1)}$ and~$\ul{X}^{(3)}$ are BT with $\cO$-structure, dual under the given polarization on~$\ul{X}$.

The first fact we need is that the ``middle'' factor $\ul{X}^{(2)}$ is rigid, as a polarized object. Thus, for $R \in \catC_{W(K)}$ there is a unique lifting $\ulcX^{(2)}_R$ of $\ul{X}^{(2)}$ over~$R$. In particular there is a canonical lifting $\ul{X}^{(2),\can}$ over~$W(K)$. The factors $\ul{X}^{(1)}$ and $\ul{X}^{(3)}$, as BT with $\cO$-structure, are rigid too, and we use a similar notation for their liftings.

(b) Consider the functor $\mE$ from $\catC_{W(K)}$ to $\cO$-modules given by $R \mapsto \Ext_R\big(\ulcX^{(1)}_R,\ulcX^{(2)}_R\big)$. Here we view Barsotti-Tate groups as sheaves for the flat topology, and the $\Ext$ is taken in the category of sheaves of $\cO$-modules on~$\Spec(R)$. We claim that $\mE$ is represented by a BT with $\cO$-structure over~$W(K)$ which is geometrically isomorphic to the product of $q-1$ copies of~$\ul{X}^{(2),\can}$. To see this, we first observe that we have a morphism of functors $\Def\big(\ul{X}^{(1)} \times \ul{X}^{(2),\prime}\big) \to \Def\big(\ul{X}^{(2),\prime}\big)$, where the prime indicates that we now view $\ul{X}^{(2)}$ {\it without\/} its polarization. Using Grothendieck-Messing deformation theory one can show that this morphism is formally smooth. Now the object $\ul{X}^{(2),\can}$ gives us a section $\Spf\big(W(K)\big) \to \Def\big(\ul{X}^{(2)}\big)$, and the functor~$\mE$ represents the pull-back of $\Def\big(\ul{X}^{(1)} \times \ul{X}^{(2)}\big)$ via this section. In this way we see that~$\mE$ is pro-representable and formally smooth.

{}From now on let us assume that $K = k$ is algebraically closed. We identify $\cO = W(\kappa)$ and write $L$ for its fraction field. Note that $L/\cO = \ul{X}_\et$. We have a short exact sequence of sheaves of $\cO$-modules $0 \to \cO \to L \to L/\cO \to 0$. As $\Hom_R\big(L,\ulcX^{(2)}_R\big) = 0$, this gives rise to injective maps
$$
\ul{X}^{(2),\can}(R) = \Hom_R\big(\cO,\ulcX^{(2)}_R\big) \longhookrightarrow \Ext_R\big(L/\cO,\ulcX^{(2)}_R\big)\, ,
$$
functorial in~$R$. Put differently, we have an injective map $j\colon \big(\ul{X}^{(2),\can}\big)^{q-1} \hookrightarrow \mE$. One easily verifies that this map is an isomorphism on tangent spaces. Hence $j$ is an isomorphism and $\mE_k$ is isomorphic to the product of $q-1$ factors~$\ul{X}^{(2),\can}$.

As Serre duality gives an isomorphism of~$\mE$ with the functor $R \mapsto \Ext_R\big(\ulcX^{(2)}_R,\ulcX^{(3)}_R\big)$, we have the same conclusions for the latter.

(c) Let $\ul{X}^\prime$ be the pair $(X,\iota)$, without polarization form. Any deformation of $\ul{X}^\prime$ over~$R$ admits a slope decomposition, with graded pieces $\ulcX^{(3)}_R$, $\ul{Z}$ and $\ulcX^{(1)}_R$, where $\ul{Z}$ is a deformation of~$\ul{X}^{(2),\prime}$. Consider the closed formal subscheme $\mD \subset \Def(\ul{X}^\prime)$ given by the condition that $Z = \ulcX^{(2)}_R$; this is equivalent to the condition that the polarization form on~$\ul{X}^{(2)}$ lifts to~$Z$. The slope filtration gives rise to a morphism $\mD \to \mE \times \mE$, where the first (resp.\ second) factor~$\mE$ controls the extension of $\ulcX^{(2)}_R$ by $\ulcX^{(1)}_R$ (resp.\ the extension of $\ulcX^{(3)}_R$ by $\ulcX^{(2)}_R$).

Similar to our construction in~\refn{CascConstr}, we can give $\mD$ the structure of a biextension (= $3$-cascade) over~$\mE \times \mE$. The structure group is of course $\Ext\big(\ul{X}^{(1)},\ul{X}^{(3)}\big)$, which is a formal torus of rank~$(q-1)^2$.

(d) Finally, the deformations of~$\ul{X}$ are parametrized by a closed formal subscheme $\mD^\lambda \subset \mD$, defined as the fixed point locus in~$\mD$ of an involution $\mD \isomarrow \mD^\vee$. This fixed point locus lives over the diagonal in $\mE \times \mE$; its fibres are principal homogeneous under a formal torus of rank~$q(q-1)/2$.

\section{Moduli spaces of PEL type, and congruence relations}{ModulPEL}

\ssection{The Ekedahl-Oort stratification on moduli spaces of PEL type}{EOonPEL}

\sssection{PELdata}
We consider a moduli problem of PEL type with good reduction at a prime $p>2$. The data involved are the following.

\item{---} $(\cB,\ast)$ is a finite dimensional semi-simple $\mQ$-algebra with a positive involution;
\item{---} $\cV$ is a finitely generated faithful left $\cB$-module;
\item{---} $\phi\colon \cV \times \cV \rightarrow \mQ$ is a symplectic form ($\mQ$-bilinear, alternating and perfect) with 
the property that $\phi(bv_1,v_2) = \phi(v_1,b^\ast v_2)$ for all $b \in \cB$ and $v_1$, $v_2 \in \cV$;
\item{---} $p$ is a prime number $>2$ such that $\cB \otimes \Qp$ is unramified (see~\refn{BasicCas});
\item{---} $O_\cB$ is a $\mZ_{(p)}$-order in~$\cB$, stable under~$\ast$, such that $O_\cB \otimes \Zp$ is a maximal order in 
$\cB \otimes \Qp$;
\item{---} $\Lambda \subset \cV \otimes \Qp$ is a $\Zp$-lattice which is also an $O_\cB$-submodule, such that $\phi$ induces 
a perfect pairing $\Lambda \times \Lambda \to \Zp$;
\item{---} $\cG := \CSp(\Lambda,\phi) \cap \GL_{O_\cB \otimes \Zp}(\Lambda)$ is the (not necessarily connected) reductive 
group over~$\Zp$ given by the symplectic similitudes of $(\Lambda,\phi)$ that commute with the action of~$O_\cB$;
\item{---} $\cX$ is a $\cG(\mR)$-conjugacy class of homomorphisms $\mS \to \cG_\mR$ (with $\mS := \Res_{\mC/\mR} \mG_m$) that 
define a Hodge structure of type $(-1,0) + (0,-1)$ on~$\cV_\mR$ for which either $2\pi i\cdot \phi$ or $-2\pi i\cdot \phi$ 
is a polarization form;
\item{---} $\gc$ is the $\cG(\mC)$-conjugacy class of cocharacters of~$\cG_\mC$ associated to~$\cX$; concretely, if 
$h\in\cX$ then we have a cocharacter $\mu=\mu_h$ through which $z\in\mC^\times$ acts on $\cV^{-1,0}$ (resp.\ $\cV^{0,-1}$) 
as multiplication by~$z$ (resp.\ by~$1$);
\item{---} $E = E(\cG,\cX)$ is the reflex field, i.e., the field of definition of the conjugacy class~$\gc$.

\sssection{ADdef}
Fix data $\cD = (\cB,\ast,\cV,\phi,O_\cB,\Lambda,\cX)$ as in \refn{PELdata}. Let $\Qbar$ be the algebraic closure of~$\mQ$ inside~$\mC$. We fix an embedding $\Qbar \to \Qbar_p$. Let $v$ be the corresponding place of~$E$ above~$(p)$. We write $O_{E,v}$ for the localization of~$O_E$ at~$v$.

Let $C_p := \cG(\Zp)$. Let $C^p$ be a compact open subgroup of $\cG(\mA^p_f)$, and put $C := C_p \times C^p$. We consider 
the moduli problem $\cA_{\cD,C}$ over $\Spec(O_{E,v})$ defined by Kottwitz in [\Kott], \S~5. If $T$ is a locally noetherian 
$O_{E,v}$-scheme then the $T$-valued points of $\cA_{\cD,C}$ are the isomorphism classes of four-tuples 
$\ul{A} = (A,\bar\lambda,\iota,\bar\eta)$ with
\item{---} $A$ an abelian scheme up to prime-to-$p$ isogeny over~$T$;
\item{---} $\bar\lambda \in \big(\NS(A) \otimes \Zp\big)/\mZ_p^\times$ the class of a prime-to-$p$ polarization;
\item{---} $\iota \colon O_\cB \to \End_T(A) \otimes \mZ_{(p)}$ a homomorphism of $\mZ_{(p)}$-algebras with $\iota(b^\ast) = 
\iota(b)^\dagger$; here $\dagger$ is the Rosati involution associated to~$\bar\lambda$;
\item{---} $\bar\eta$ a level structure of type $C^p$ on~$A$;

\noindent
such that a certain determinant condition is satisfied. For precise details we refer to Kottwitz~[\Kott], \S~5. If $C^p$ is 
sufficiently small, which we from now on assume, then $\cA_{\cD,C}$ is representable by a smooth quasi-projective 
$O_{E,v}$-scheme.

\sssection{ADBasCas}
In the rest of this section we assume that $(\cB,\ast)$ is simple as an algebra with involution. Then $(\cB,\ast)$ is of one of the four types \Romno 1--\Romno 4 in Albert's classification; see Mumford~[\MAV], \S~21. Let $Z = Z_\cB$ be the center of~$\cB$ and $Z_0 \subset Z$ the subfield of $\ast$-symmetric elements. Define $\cO := (O_\cB \otimes \Zp) \cap (Z \otimes \Qp)$. We again write~$\ast$ for the involution of~$\cO$ induced by the given involution of~$\cB$. If $(\cB,\ast)$ is of Albert type \Romno 1 or \Romno 2, set $\epsilon = -1$; otherwise set $\epsilon = +1$. The triple $(\cO,\ast,\epsilon)$ thus obtained is a product of triples of type C, D, AU or~AL (cf.\ section~\refn{BasicCas}); more precisely:

\item{---} if $(\cB,\ast)$ is of type \Romno1 or \Romno2 then $(\cO,\ast,\epsilon)$ is a product of triples of type C, where 
the factors are indexed by the primes of~$Z$ above~$p$;
\item{---} if $(\cB,\ast)$ is of type \Romno3 then $(\cO,\ast,\epsilon)$ is a product of triples of type D, where the 
factors are indexed by the primes of~$Z$ above~$p$;
\item{---} if $(\cB,\ast)$ is of type \Romno4 then $(\cO,\ast,\epsilon)$ is a product of triples of type AU and AL; the 
factors of type AU correspond to the primes of~$Z_0$ above~$p$ that are inert in the extension $Z_0 \subset Z$; the factors 
of type AL correspond to the primes of~$Z_0$ above~$p$ that split in~$Z$.
\smallskip

Let $T$ be a scheme over~$O_{E,v}$. Let $s \in \cA_{\cD,C}(T)$ be a $T$-valued point, corresponding to a four-tuple 
$(A,\bar\lambda,\iota,\bar\eta)$. The Barsotti-Tate group $A[p^\infty]$ has a $(O_\cB \otimes \Zp,\ast,-1)$-structure. By 
our assumptions, $O_\cB \otimes \Zp$ is isomorphic to a matrix algebra over~$\cO$. Therefore Morita equivalence 
applies, to the effect that $A[p^\infty]$ comes from a BT $\ul{X} = \ul{X}_s$ with $(\cO,\ast,\epsilon)$-structure.

Morita equivalence also applies to $(\Lambda,\phi)$. For the rest of this section we fix a pair $(\Lambda_0,\phi_0)$ consisting of an $\cO$-module and an $\epsilon$-$\ast$-hermitian pairing such that $(\Lambda_0,\phi_0)$ is Morita equivalent to the original pair $(\Lambda,\phi)$. Note that $\cG = \CSp_\cO(\Lambda_0,\phi_0)$.

\sssection{Sh(G0,X0)}
Let $\cX^0 \subset \cX$ be a $\cG^0(\mR)$-orbit. The pair $(\cG^0,\cX^0)$ is a Shimura datum. Define $\gc^0$ to be the 
$\cG^0(\mC)$-conjugacy class of cocharacters of~$\cG_\mC$ with $\mu_h \in \gc^0$ for all $h\in\cX^0$.    The reflex field 
$E^0 := E(\cG^0,\cX^0)$ is a finite extension of $E = E(\cG,\cX)$. Let $v^0$ be the place of~$E^0$ determined by the chosen embedding $\Qbar \to \Qbar_p$.

Let $S_C = S_C(\cG^0,\cX^0)$ denote the canonical model (over~$E^0$) of the Shimura variety associated to $(\cG^0,\cX^0)$ at 
level $C \cap \cG^0(\mA_f)$. Then $S_C$ can be identified with an open and closed subscheme of the generic fibre 
of~$\cA_{\cD,C} \otimes O_{E^0,v^0}$. In fact, we have a decomposition of $\cA_{\cD,C} \otimes \mC$ as a union of open and 
closed subschemes, say
$$
\cA_{\cD,C} \otimes \mC = \cA^{(1)} \amalg \cdots \amalg \cA^{(s)}\, ,
$$
such that each $\cA^{(j)}$ is a Shimura variety. Here $s$ is the order of $\Ker\big(H^1(\mQ,\cG) \to \prod_p 
H^1(\Qp,\cG)\big)$. In general, the Shimura varieties that constitute the generic fibre are not all associated to the same 
$\mQ$-group. (Note that different PEL data $\cD$ may give rise to the same moduli problem $\cA_{\cD,C}$, since this problem 
only involves all local information.)

For some results that we want to discuss it is more natural to work with the individual Shimura varieties. We define $\cS_C$ 
to be the open and closed subscheme of~$\cA_{\cD,C} \otimes O_{E^0,v^0}$ whose generic fibre is~$S_C$. If there is no risk 
of confusion we simply write~$\cS$ instead of~$\cS_C$. Write $\cS_0 = \cS_{C,0}$ for the special fibre.

\sssection{TypeCst/S}
{\it Proposition. --- Let $k$ be an algebraically closed field containing~$\kappa(v^0)$. In the cases C, AU or AL (Albert 
types \Romno1, \Romno2 or \Romno4), the type-function $s \mapsto (d_s,\gf_s)$ is constant on~$\cS_0(k)$. In case~D (Albert 
type~\Romno3), the type-function $s \mapsto (d_s,\gf_s,\delta_s)$ is constant on~$\cS_0(k)$.}
\Cskip

In fact, it is clear that $d$ and~$\gf$ are constant. We postpone the proof of the assertion in case~D to the end of this subsection.

\sssection{EOStrat}
Let $k$ be an algebraically closed field containing~$\kappa(v^0)$. A $k$-valued point of~$\cS_0$ gives rise to a BT $\ul{X} 
= \ul{X}_s$ with $(\cO,\ast,\epsilon)$-structure, of some fixed type $(d,\gf)$ or, in case~D, $(d,\gf,\delta)$. Let 
$(G^0,\mX^0)$ be the corresponding pair consisting of an algebraic group and a conjugacy class of parabolic subgroups; see~\refn{wdefGPE}. Our classification results of the $p$-kernel group schemes $\ul{Y} := \ul{X}[p]$ give rise to an {\it Ekedahl-Oort stratification\/}
$$
\cS_0 = \disunion_{w \in W_{\mX^0}\backslash W_{G^0}} \cS_0(w)\, ,\leqno{(\refn{EOStrat}.1)}
$$
where $s \in \cS_0(k)$ lies in $\cS_0(w)$ if and only if $\ul{Y}_s$ is of type~$w$. By the statement that (\refn{EOStrat}.1) is a stratification we mean that it gives a decomposition of~$\cS_0$ as a disjoint union of locally closed subspaces and that the closure of each stratum~$\cS_0(w)$ is a union of strata. For proofs of these facts we refer to Wedhorn's paper~[\Wedh].

In our paper~[\DFEO] we have proved the following result.

\sssection{EODimForm}
{\it Theorem. --- If $\cS_0(w) \neq \emptyset$ then $\cS_0(w)$ is equi-dimensional, of dimension $\ell(w)$.}
\Cskip

Combining this with Thm.~\refn{OrdPolThm} we obtain a new proof of the main result of Wedhorn~[\WedhOrd]:

\sssection{OrdDense}
{\it Corollary. {\rm (Wedhorn)} --- The ordinary locus in $\cS_0$ is Zariski dense.}

\sssection{TpLam0}
Let $Q$ be a field containing $E^0$. Let $\ul{X}$ be the BT with $(\cO,\ast,\epsilon)$-structure associated to a $Q$-valued point of~$\cS$. Write $T_p = T_p(\ul{X})$ for its Tate-$p$-module, which is to be viewed as a free $\Zp$-module of finite rank with $\cO$-action and with an $\epsilon$-symmetric perfect bilinear form $\psi_p \colon T_p \times T_p \to \Zp(1)$ satisfying $\psi_p(bx,y) = \psi_p(x,b^\ast y)$ for all $b \in \cO$ and $x$, $y \in T_p$. 

The interpretation of the generic fibre of~$\cS$ as the Shimura variety associated to the datum $(\cG^0,\cX^0)$ gives us an isomorphism of $\cO$-modules $\alpha \colon \Lambda_0 \isomarrow T_p$ such that $\alpha^\ast \psi_p = c \cdot \phi_0$ for some $c \in \Zp(1)^\times$. (See~\refn{ADBasCas} for the definition of $(\Lambda_0,\phi_0)$.) This isomorphism~$\alpha$ is canonical up to the action of an element of~$\cG^0(\Zp)$.

\sssection{TypeCstPf}
{\it Proof of Proposition\/ {\rm \refn{TypeCst/S}}.}\enspace The only non-trivial part of the proposition is the statement that, in case~D, the function $\delta$ is constant. 

Suppose we are in case~D. As shown in~[\DFEO], Lemma.~3.1.4, the function~$s \mapsto \delta_s$ is locally constant in families. On the other hand, the ordinary locus of $\cS_0$ is Zariski dense; see~\refn{OrdDense}. Hence it suffices to show that any two ordinary points of~$\cS_0$ have the same~$\delta$.

Let $\ul{A}$ be an ordinary $k$-valued point of~$\cS$. Write $\ul{X}$ for the associated BT with $(\cO,\id,+1)$-structure. For simplicity of exposition, let us assume that $\cO$ is a domain, i.e., there is only one prime of~$Z = Z_\cB$ above~$p$. In the general case the argument is the same, but we first have to decompose~$\ul{X}$ according to the decomposition of~$\cO$ as a product of domains. As usual we write $\kappa = \cO/p\cO$, we let $\cI = \Hom(\kappa,k) = \Hom\big(\cO,W(k)\big)$, and we put $m=\#\cI = [\cO:\Zp]$. Recall that we have an integer~$q$ such that $d=2q$ and $\gf$ is the constant function~$q$.

Let $\delta$ be the invariant of~$\ul{X}$ as in~\refn{dfdelta}. We say we are in the {\it split case\/} if $\sum_{i\in\cI} \delta(i) \equiv m$ modulo~$2$, in the {\it non-split case\/} if not. As we shall see, this is independent of the choice of the ordinary point~$\ul{A}$.

In the split case, let $\cOtil = \cO \times \cO$ with involution~$\ast$ given by $(y_1,y_2) \mapsto (y_2,y_1)$. In the non-split case, let $\cOtil$ be the unramified quadratic extension of~$\cO$ and $\ast$ the non-trivial automorphism of $\cOtil$ over~$\cO$. Set $\cItil := \Hom\big(\cOtil,W(k)\big)$. We use the letter~$\tau$ for elements of~$\cItil$. Similar to the notation introduced in~\refn{GPENot} we have a natural $2:1$ map $\cItil \to \cI$ and we set $\bar\tau := \tau \circ \ast$.

Let $R = \cO^{q-1} \times \cOtil \times \cO^{q-1}$ with involution~$\ast$ given by $(x,y,z)^\ast = (z,y^\ast,x)$. Let $\cT$ be the torus over~$\Zp$, of rank $mq+1$, given on points by
$$
\cT(A) = \big\{\xi \in (R \otimes_\Zp A)^\times \bigm| \xi \xi^* \in A^\times\big\}\, .
$$
The cocharacter group of~$\cT$ is given by
$$
X_\ast(\cT) = \big\{(a_{i,j},b_\tau,c_{i,j}) \in \big(\mZ^\cI\big)^{q-1} \times \mZ^{\cItil} \times \big(\mZ^\cI\big)^{q-1} 
\bigm| a_{i,j} + c_{i,q-j} = \hbox{const} = b_\tau + b_{\bar\tau}\big\}\, .
$$
The constant appearing here is called the weight of the cocharacter. The fundamental group of $\Spec(\Zp)$ acts on~$X_\ast(\cT)$ through its natural action on the sets~$\cI$ and~$\cItil$. 

Let $\nu\colon \mG_m \to \cT$ be a cocharacter of weight~$1$ over~$W(k)$. To~$\nu$ we associate a Dieudonn\'e module with $(R,\ast,+1)$-structure: Take $M_\nu = R \otimes_\Zp W(k)$, with $F$ and~$V$ given by
$$
F(r \otimes w) = (1 \otimes \sigma)\big(\nu(p) \cdot (r \otimes w)\big)
\quad\hbox{and}\quad
V(r \otimes w) = (1 \otimes \sigma^{-1})\big(\nu(p)^\ast \cdot (r \otimes w)\big)\, ,
$$
and with $+1$-duality given by $\psi(r_1 \otimes w_1, r_2 \otimes w_2) = \trace_{R/\Zp}(r_1 r_2^\ast) w_1 w_2$. If instead of the full $R$-action we only remember the action of~$\cO$ (embedded diagonally into~$R$) then we obtain a BT with $(\cO,\id,+1)$-structure, denoted~$\ul{X}_\nu$. For later use let us remark that $X_\nu$ is ordinary when viewed as a BT with $R$-structure; the point is that it has height~$1$ (``relative to its $R$-structure''), and height~$1$ objects are always ordinary.

The point of all this is that~$\ul{X}$, our ordinary BT with $(\cO,\id,+1)$-structure, is of the form $\ul{X}=\ul{X}_\nu$ for some cocharacter~$\nu$. In the given description of the cocharacter group, we can choose~$\nu$ in such a way that $a_{i,j} = 0$ and $c_{i,j} = 1$ for all~$i\in\cI$ and $j\in \{1,\ldots,q-1\}$; with $\ul{X} = \ul{X}^{(1)} \times \ul{X}^{(2)} \times \ul{X}^{(3)}$ as in~\refn{PolDefD} this is equivalent to the requirement that the $n$th factor ($n=1,2,3$) of $R = \cO^{q-1} \times \cOtil \times \cO^{q-1}$ acts by endomorphisms of~$\ul{X}^{(n)}$. Once we fix the $R$-action on~$\ul{X}$, the corresponding $\nu$ is uniquely determined. 

We can compute~$\delta$ from~$\nu$, as follows. Let $i \in \cI$. Choose an element $\tau \in \cItil$ that maps to~$i$ under the natural $2:1$ map $\cItil \to \cI$. Write $\nu = (a_{i,j},b_\tau,c_{i,j})$. Then
$$
\delta(i) = \cases{0 \bmod 2&if $b_\tau \neq b_{\tau-1}$;\cr
1 \bmod 2&if $b_\tau = b_{\tau-1}$.\cr}
$$

Let $T_p = T_p(\ul{X}^\can)$ be the Tate-$p$-module of the canonical lifting of~$\ul{X}$. As discussed in~\refn{TpLam0} we have an  isomorphism $\alpha\colon \Lambda_0 \isomarrow T_p$, canonical up to an element of~$\cG^0(\Zp)$. As remarked above, $X$ is ordinary when viewed as a BT with $R$-structure, so by~\refn{TorsPts} the full $R$-action on~$X$ lifts to~$X^\can$. Via~$\alpha$ this gives rise to an embedding $j\colon \cT \hookrightarrow \cG^0$ (over~$\Zp$), realizing~$\cT$ as a maximal torus of~$\cG^0$. The cocharacter~$\nu$ is defined over a finite unramified extension~$V$ of $W\big(k(v^0)\big)$ inside~$W(k)$. Choose an embedding $V \hookrightarrow \Qpbar$. We obtain a cocharacter $j \circ \nu$ of~$\cG^0$ over~$\Qpbar$. On the other hand, writing $\gC(L)$ for the $\cG^0(L)$-conjugacy classes of cocharacters $\mG_m \to \cG^0$ over a field~$L$, we have natural bijections $\gC(\mC) \isomworra \gC(\Qbar) \isomarrow \gC(\Qpbar)$, via which we can view~$\gc^0$ (as in~\refn{Sh(G0,X0)}) as an element of~$\gC(\Qpbar)$. By Reimann-Zink~[\ReimZink], Thm.~(1.6), we have $j \circ \nu \in \gc^0$. (This may be off by a normalization factor, due to the fact that we use a different version of Dieudonn\'e theory, and due to various sign conventions in Hodge theory. Such a normalization does not affect our argument, though, and we save ourselves the trouble of getting it exactly right.)

As claimed earlier, whether we are in the split or in the non-split case is independent of the choice of the ordinary point~$\ul{A} \in \cS(k)$. Let us now prove this. Let $\cG_1 = \Ker(c) \subset \cG^0$, where $c\colon \cG \to \mG_m$ is the multiplier character. We have $\cG_1 = \Res_{\cO/\Zp} \cG_1^\prime$ with $\cG_1^\prime$ an algebraic group over~$\cO$. We claim that we are in the split case if and only if $\cG_1^\prime$ is split (over~$\cO$). If we are in the split case then $\cT_1 := j(\cT) \cap \cG_1$ is of the form $\cT_1 = \Res_{\cO/\Zp} \cT_1^\prime$ and $\cT_1^\prime \subset \cG_1^\prime$ is a split maximal torus. Conversely, if $\cG_1^\prime$ is split then every element of $\gC(\Qpbar)$ is defined over the fraction field of~$\cO$. But if $\sum_{i\in\cI} \delta(i) \not\equiv m$ modulo~$2$ then we find that the conjugacy class of~$j \circ \nu$ is defined only over a quadratic extension of~$\cO$. This proves our claim.

To complete the proof, let us now show that $\delta$ is determined by the conjugacy class~$\gc^0$. If $W_{\cG^0}$ is the Weyl group of~$\cG^0$ then there is a natural bijection $\gC(\Qpbar) \isomarrow X_\ast(\cT)/W_{\cG^0}$. We observe that $\nu \in X_\ast(\cT)$ is the unique representative of the class $\gc^0 \in X_\ast(\cT)/W_{\cG^0}$ with the property 
that for all $\tau \in \cItil$ and $i = \res(\tau) \in\cI$,
$$
a_{i,1} \leq \cdots \leq a_{i,q-1} \leq \min(b_{\tau},b_{\bar\tau}) < \max(b_\tau,b_{\bar\tau}) \leq c_{i,1} \leq \cdots \leq c_{i,q-1}\, .
$$
Since we can compute~$\delta$ from~$\nu$, it follows that $\gc^0$ determines~$\delta$. \QED

\sssection{kappa(v)}
{\it Remark.\/} --- There are two key points in the above proof. Firstly, we have a direct relation between the conjugacy class~$\gc^0$ and the conjugacy class of the cocharacter~$\nu$. Secondly, the cocharacter~$\nu$ is directly related to the explicit description of the Dieudonn\'e module of $\ul{X}^\ord(d,\gf,\delta)$. We can further exploit these relations to obtain information on the residue field $\kappa(v^0)$ of the reflex field~$E^0$ at the place~$v^0$. In the cases~A and~C (where similar ideas apply), we find that $\kappa(v) = \kappa(v^0)$ equals the field $E(d,\gf)$ defined in~\refn{modpRefl}. 

In case~D something similar can be done. Given a type $(2q,q,\delta)$, define $E(\delta) \subset k$ as the fixed field of $\{\alpha \in \Aut(k) \mid {}^\alpha \delta = \delta\}$. Then we find that in case~D(split) we have $\kappa(v^0) = E(\delta)$ and in case D(non-split) $\kappa(v^0)$ is the quadratic extension of~$E(\delta)$ in~$k$.

\sssection{deltaval}
In case~D, let $\cG_1 = \Res_{\cO/\Zp} \cG_1^\prime$ as in the above proof. Set $n=0$ if $\cG_1^\prime$ is a split group, $n=1$ otherwise. Then every value for~$\delta$ such that $\sum_{i\in\cI} \delta(i) \equiv m+n$ modulo~$2$ occurs on the special fibre of~$\cA_{\cD,C}$. To see this we have to analyse what happens if we replace $\cX^0 \subset \cX$ by another $\cG^0(\mR)$-orbit. This amounts to changing the conjugacy class~$\gc^0$ by an element of~$\cG$. Computing~$\delta$ from~$\gc^0$ as in the above proof, we find that all~$\delta$ with $\sum_{i\in\cI} \delta(i) \equiv m+n$ are obtained.

\ssection{Congruence relations}{CongRel}

\sssection{pIsog}
We retain the notation of section~\refn{EOonPEL}. If there is no risk of confusion we simply write $\cA$ for~$\cA_{\cD,C}$. Write $\cA^\ord \subset \cA$ for the open subscheme obtained by removing the non-ordinary locus on the special fibre. (Note: if we decompose $\cA$ as a union of Shimura varieties, as in~\refn{Sh(G0,X0)}, then on {\it each\/} of the ``Shimura components'' we remove the non-ordinary locus.)

Suppose given two four-tuples $\ul{A}_i = (A_i,\bar\lambda_i,\iota_i,\bar\eta_i)$, for $i=1$, $2$, corresponding to $T$-valued points of~$\cA$. By a {\it $p$-isogeny\/} $f\colon \ul{A}_1 \to \ul{A}_2$ we mean an $O_\cB$-linear isogeny such that $f^\ast \bar\lambda_2 = p^c \cdot \bar\lambda_1$ for some $c \geq 0$. Note that $f$ necessarily has $p$-power degree.

Write $\pIsog = \pIsog_{\cD,C}$ for the $O_{E,v}$-scheme of such $p$-isogenies; it comes equipped with two morphisms $s$, $t\colon \pIsog \to \cA$, sending an isogeny to its source and target, respectively. Fixing $p^d = \deg(f)$ gives an open and closed subscheme $\pIsog^{(d)} \subset \pIsog$ which is locally of finite type over~$O_{E,v}$. Using [\FaCh], Chap.~\Romno 1, Prop.~2.7 and the valuative criterion we see that the morphisms~$s$ and $t\colon \pIsog^{(d)} \to \cA$ are proper. 

For every $n \geq 0$ we have a morphism $\cA \to \pIsog$ that sends $\ul{A}$ to the isogeny ``multiplication by~$p^n$ on~$\ul{A}$''. As this is clearly a section of the (separable) morphism~$s$, its image is a reduced closed subscheme $\Mult(p^n) \subset \pIsog$.

Composition of isogenies defines a morphism
$$
c\colon \pIsog \times_{t,\cA,s} \pIsog \tto \pIsog\, .
$$
We claim that this morphism is proper. To see this, work over a d.v.r.\ $R$ with fraction field~$K$. Suppose given an isogeny $f\colon \ul{A}_1 \to \ul{A}_2$, and suppose that $f_K = \psi_K \circ \phi_K$. Since $R$ is a d.v.r., the flat closure of $\Ker(\phi_K)$ inside~$A_1$ is a finite flat subgroup scheme, and we get (in a unique way) a factorization $f = \psi \circ \phi$ over~$R$. By the valuative criterion it follows that $c$ is proper.

Suppose given a homomorphism $O_{E,v} \to L$ with $L$ a field. Using the composition morphism~$c$ we can define an ``algebra of isogenies'' over~$L$. Let $Z_\mQ(\pIsog \otimes L)$ be the group of algebraic cycles on $\pIsog \otimes L$, taken with $\mQ$-coefficients. If $Y_1$ and $Y_2$ are two cycles, let
$$
Y_1 \cdot Y_2 := c_\ast(Y_1 \times_{t,s} Y_2)\, .
$$
Note that the push-forward is defined on the level of cycles, since $c$ is proper. Extending this product bilinearly, we obtain the structure of a $\mQ$-algebra on~$Z_\mQ(\pIsog \otimes L)$. The identity element is the cycle $\Mult(1) = \pIsog^{(0)}$. Finally define $\mQ[\pIsog \otimes L]$ to be the subalgebra of $Z_\mQ(\pIsog \otimes L)$ generated by the irreducible components.

The previous constructions also work on $\pIsog^\ord$, which we define as the inverse image under~$s$ of $\cA^\ord \subset \cA$. As our notion of ordinariness is invariant under isogenies, $\pIsog^\ord$ is also the inverse image of $\cA^\ord \times \cA^\ord$ under $(s,t)$. Let $\pIsog^{\ord,(d)}$ be the open and closed subscheme of $p$-isogenies~$f$ with $\deg(f) = p^d$. 

\sssection{AlgIsogLem}
{\it Lemma. --- Suppose given a homomorphism $O_{E,v} \to L$ with $L$ a field.

{\rm (\romno1)} If $\charact(L) = 0$ then $\mQ[\pIsog \otimes L] \subset Z_\mQ(\pIsog \otimes L)$ is the $\mQ$-subspace spanned by the irreducible components of $\pIsog \otimes L$. In other words, if $Y_1$ and $Y_2$ are irreducible components of $\pIsog \otimes L$ then $Y_1 \cdot Y_2$ is a $\mQ$-linear combination of irreducible components.

{\rm (\romno2)} If $\charact(L) = p$ an analogous statement holds over the ordinary locus: $\mQ[\pIsog^\ord \otimes L] \subset Z_\mQ(\pIsog^\ord \otimes L)$ is the $\mQ$-subspace spanned by the irreducible components of $\pIsog^\ord \otimes L$.}
\Dskip

\Proof~If $\charact(L)=0$ we reduce to the case $L=\mC$; then we use the complex uniformization of the components of $\cA \otimes \mC$ by hermitian symmetric domains. We omit the details. Next suppose $\charact(L) = p$. For the purpose of this proof let us abbreviate $\pIsog^\ord \otimes L$ to~$\cI$.

{\it Step 1:\/} It suffices to show that, fixing~$d$, the morphisms $s$, $t \colon \cI^{(d)} \to \cA^\ord \otimes L$ are finite and flat. Indeed, suppose this is true. Let $Y_1$ and~$Y_2$ be irreducible components of~$\cI$, say with $Y_i$ contained in $\cI^{(d_i)}$. Then $W := \cI^{(d_1)} \times_{t,s} \cI^{(d_2)}$ is finite flat over $\cA^\ord \otimes L$, and $Y_1 \times_{t,s} Y_2$ is a union of irreducible components of~$W$. It follows that $\cI$ and $Y_1 \times_{t,s} Y_2$ are both equidimensional of dimension equal to $\dim(\cA)$. This readily implies that $Y_1 \cdot Y_2$ is a linear combination of irreducible components of~$\cI$.

{\it Step 2:\/} Let $\cS_0 \subset \cA \otimes \kappa(v^0)$ be the special fibre of the Shimura variety~$\cS$; see the discussion in~\refn{Sh(G0,X0)}. Let $\cS_0^\ord$ be the ordinary locus, and let $\cJ_0^\ord \subset \cI$ be the inverse image of $\cS_0^\ord$ under the source morphism~$s$. Let $k$ be an algebraically closed field containing~$L$. A $k$-valued point $a \in \cS_0^\ord(k)$ gives rise to a four-tuple $\ul{A}$, which in turn gives rise to a BT $A[p^\infty]$ with $(O_\cB \otimes \Zp,\ast,-1)$-structure (cf.~\refn{ADBasCas}). By Thm.~\refn{OrdPolThm} and Prop.~\refn{TypeCst/S}, this object is independent of the choice of~$a$, up to isomorphism.

If $R$ is a $k$-algebra and if $f\colon \ul{A}_1 \to \ul{A}_2$ corresponds to an $R$-valued point of~$\cJ_0^{\ord,(d)}$ then up to isomorphism $f$ only depends on its kernel $\Ker(f) \subset \ul{A}_1$. It follows that the fibres $s^{-1}(a)$, for $a \in \cS_0^\ord(k)$, are all isomorphic as schemes.

{\it Step 3:\/} We claim that, fixing $d = \log_p\big(\deg(f)\big)$, the fibres $s^{-1}(a)$ are finite. First note that any $p$-isogeny~$f$ factors as $f = f_r \circ f_{r-1} \circ \cdots \circ f_1$ in such a way that $\Ker(f_n \circ f_{n-1} \circ \cdots \circ f_1) = \Ker(f)[p^n]$ for all~$n$. Moreover, up to isomorphism this factorization is unique. It therefore suffices to prove that in the fibre over a given point $a \in \cS_0^\ord(k)$ there are only finitely many $p$-isogenies~$f$ with the property that $\Ker(f)$ is killed by~$p$. Such an isogeny is completely determined by the induced homomorphism $f[p]\colon \ul{A}_1[p] \to \ul{A}_2[p]$. Moreover, sending $(a_1,a_2)$ to $\big(a_1,a_2 + f[p](a_1)\big)$ gives an automorphism of $\ul{A}_1[p] \times \ul{A}_2[p]$ as a BT$_1$ with $(O_\cB/pO_\cB,\ast,-1)$-structure. But $\ul{A}_1[p] \times \ul{A}_2[p]$ is ordinary, so by Thm.~2.1.2 of~[\DFEO] its automorphism group scheme is finite.  

{\it Step 4:\/} We shall use the following general fact: If $\phi\colon X \to Y$ is a finite morphism of schemes such that $Y$ is reduced and such that the function $y \mapsto \dim_{\kappa(y)}\big(\phi_\ast O_X \otimes_{O_Y} \kappa(y)\big)$ is constant on~$Y$ then $\phi$ is flat. By what was explained in Steps 2 and~3 we can apply this to the morphism $s\colon \cJ_0^\ord \to \cS_0^\ord$. As finite flatness is a local notion on the target scheme, this shows that $s \colon \cI^{(d)} \to \cA^\ord \otimes L$ is finite and flat.

Finally, that the target morphism $t$ is also finite and flat follows from duality. Namely, sending a $p$-isogeny~$f\colon \ul{A}_1 \to \ul{A}_2$ to the dual isogeny $f^t \colon \ul{A}_2^t \to \ul{A}_1^t$ gives an isomorphism between $\pIsog$ and another scheme of $p$-isogenies, which interchanges the roles of~$s$ and~$t$. (In general the ``other'' scheme of $p$-isogenies lives over another moduli scheme of PEL type, as the dual abelian schemes $\ul{A}_i^t$ with the inherited $O_\cB$-action may have a different CM-type.) We leave the details of this to the reader. \QED

\sssection{ScJandJ}
Let $\cS = \cS_C \hookrightarrow \cA \otimes O_{E^0,v^0}$ be as in~\refn{Sh(G0,X0)}. Define $\cJ \hookrightarrow \pIsog \otimes O_{E^0,v^0}$ to be the inverse image of $\cS \times \cS$ under the morphism $(s,t) \colon \pIsog \to \cA \times \cA$. We write~$J$ for the generic fibre of~$\cJ$ and $\cJ_0$ for its special fibre. 

Consider a homomorphism $O_{E^0,v^0} \to L$ with $L$ a field. If $\charact(L)=0$ we define $\mQ[J \otimes L]$ to be the subalgebra of $\mQ[\pIsog \otimes L]$ generated by the irreducible components of~$J_L$. Similarly, if $\charact(L)=p$ then we define $\mQ[\cJ^\ord \otimes L]$ to be the subalgebra of $\mQ[\pIsog^\ord \otimes L]$ generated by the irreducible components of~$\cJ^\ord \otimes L$.

\sssection{PhionS}
Let $q=p^m$ be the cardinality of the residue field~$\kappa(v^0)$. We have a section $\phi\colon \cS_0 \to \cJ_0$ of the source morphism, sending a four-tuple~$\ul{A}$ to the $m$th power Frobenius isogeny $\phi_{\ul{A}}\colon \ul{A} \to \ul{A}^{(q)}$. The image of this section is a closed reduced subscheme $\Phi \subset \cJ_0$. As the source morphism~$s$ is finite and flat over the ordinary locus, and the ordinary locus in~$\cS_0$ is dense, $\Phi$~is a union of irreducible components of~$\cJ_0$. We shall henceforth view~$\Phi$ as an element of the algebra $\mQ[\cJ_0 \otimes \kappa(v^0)]$, or as an element of $\mQ[\cJ_0^\ord \otimes \kappa(v^0)]$. We refer to this element~$\phi$ as the Frobenius correspondence. The main theme of this section is that $\Phi$~satisfies a polynomial equation with coefficients in the Hecke algebra of~$\cG$.   

\sssection{HeckeAlg}
Recall that we have a conjugacy class $\gc^0$ of cocharacters of~$\cG^0$, and that $E^0$ is the field of definition of~$\gc^0$. Let $\cE$ be the $v^0$-adic completion of~$E^0$; its ring of integers is $O_\cE := \hat O_{E^0,v^0}$. Note that $\cE$ is an unramified extension of~$\Qp$; see~[\Milne], Cor.~4.7. By the same arguments as in~[\WedhCong], Lemma~5.1, there exists a cocharacter $\mu\in\gc^0$ that is defined over~$\cE$. Similar to what we did in~\refn{muOrdPol} we can consider the quasi-cocharacter obtained from~$\mu$ by averaging its Galois conjugates. More precisely, let $\Gamma := \Gal(\Qpbar/\Qp) \supseteq \Gamma^\prime := \Gal(\Qpbar/\cE)$, and consider
$$
N(\mu) := {1\over [\Gamma:\Gamma^\prime]}\, \sum_{\gamma \in \Gamma/\Gamma^\prime} \gamma \cdot \mu\, ,
$$
which is a quasi-cocharacter of~$\cG^0$ defined over~$\Qp$. This $N(\mu)$ extends to a quasi-cocharacter over~$\Zp$, and we define $\cM \subset \cG^0$ to be the centralizer of~$N(\mu)$.

We denote by $\cH(\cG^0,\mQ)$ the Hecke algebra of $\cG^0_\Qp$ with respect to its hyperspecial subgroup~$\cG^0(\Zp)$, with $\mQ$ as coefficient field. Let $\cH_0(\cG^0,\mQ) \subset \cH(\cG^0,\mQ)$ be the subalgebra of $\mQ$-valued functions that have support contained in $\cG^0(\Qp) \cap \End(\Lambda_0)$. The Hecke algebras $\cH_0(\cM,\Qp) \subset \cH(\cM,\mQ)$ are defined in a similar manner; see Wedhorn's paper~[\WedhCong], \S~1, for more details. We write
$$
\dot S^{\cG^0}_\cM \colon \cH(\cG^0,\mQ) \tto \cH(\cM,\mQ)
$$ 
for the twisted Satake homomorphism. It restricts to a map $\cH_0(\cG^0,\mQ) \tto \cH_0(\cM,\mQ)$, which we again call~$\dot S^{\cG^0}_\cM$.

\sssection{NmuCaseD}
{\it Remark.} --- Suppose $(\cO,\ast,\epsilon)$ is of type~D. Let $\ul{A}$ be a $K$-valued point of~$\cS_0^\ord$, where $K$ is a perfect field containing~$\kappa(v^0)$. Let $\ul{X}$ be the corresponding BT with $(\cO,\ast,\epsilon)$-structure. In~\refn{muOrdPol} we have defined conjugacy classes of quasi-cocharacters~$\bar\mu_j$ of~$\cG^0$ over~$\Qpbar$. One of these conjugacy classes contains~$N(\mu)$. By definition, saying that $\ul{A}$ is ordinary means that the associated Newton quasi-cocharacter~$\nu(\ul{X})$ lies in one of the conjugacy classes~$\bar\mu_j$. Under the assumption that $\ul{A}$ is a point of~$\cS_0$ (not just a point of~$\cA$) we can sharpen this: $\ul{A}$ is ordinary if and only if $\nu(\ul{X})$ is conjugate to~$N(\mu)$. This can be shown using arguments as in~\refn{TypeCstPf}.

\sssection{Typef}
As usual we write $T_p(?)$ for Tate-$p$-modules and $V_p(?) := T_p(?) \otimes_{\Zp} \Qp$. Let $L$ be a field containing~$E^0$. Let $f\colon \ul{A}_1 \to \ul{A}_2$ be an isogeny corresponding to an $L$-valued point of~$J$. Choose identifications $\alpha_i\colon \Lambda_0 \isomarrow T_p X_i$ as in~\refn{TpLam0}. The linear map $V_p f\colon V_p X_1 \to V_p X_2$ is an isomorphism and $\alpha_2^{-1} \circ V_p f \circ \alpha_1\colon \Lambda_0 \otimes \Qp \isomarrow \Lambda_0 \otimes \Qp$ is an element of~$\cG^0(\Qp)$. Its class
$$
\tau(f) := \left[\alpha_2^{-1} \circ V_p f \circ \alpha_1\right] \in \cG^0(\Zp)\backslash \cG^0(\Qp)/\cG^0(\Zp)
$$
is independent of the choices of the~$\alpha_i$. We refer to~$\tau(f)$ as the {\it type\/} of the $p$-isogeny~$f$.

The type of an isogeny is constant on irreducible components of the scheme~$J$. This allows us to define a map
$$
h \colon \cH_0(\cG^0,\mQ) \tto \mQ[J]
$$
sending the characteristic function of a class~$[\gamma]$ with $\gamma \in \cG^0(\Qp) \cap \End(\Lambda_0)$ to the sum of all irreducible components of~$J$ on which the type is equal to~$[\gamma]$. By extending scalars to~$\mC$ and using the complex uniformization of~$S_\mC$, it can be checked that $h$~is a homomorphism.

\sssection{ptypefA}
Our next goal is to define the {\it $p$-type\/} of an ordinary isogeny $f\colon \ul{A}_1 \to \ul{A}_2$ corresponding to a point of~$\cJ^\ord_0$. In the Siegel modular case this notion is defined by Chai and Faltings in~[\FaCh], Chap.~\Romno 7,~\S~4. The $p$-type of an isogeny~$f$ will be an element in $\cM(\Zp)\backslash \cM(\Qp)/\cM(\Zp)$.

Let $K$ be a perfect field containing the residue field~$\kappa(v^0)$. Let $\ul{A}$ be a $K$-valued point of~$\cS_0^\ord$. Let $\ul{X}$ be the corresponding BT with $(\cO,\ast,\epsilon)$-structure. Write $T_p = T_p(\ul{X}^\can)$. As before we have an identification $\alpha\colon \Lambda_0 \isomarrow T_p$, canonical up to an element of~$\cG^0(\Zp)$. 

Write $\ul{X}^\prime$ for the BT with $\cO$-structure underlying~$\ul{X}$ (forgetting the polarization). We have a slope decomposition
$$
\ul{X}^\prime = \prod\nolimits_{\nu \in \mQ} \ul{X}^{\prime,(\nu)}\, ,
$$
in such a way that the BT underlying $\ul{X}^{\prime,(\nu)}$, is isotypic of slope~$\nu$. The canonical lifting~$\ul{X}^\can$ has the property that as a BT with $\cO$-structure it is the product of factors $\ul{X}^{\can,\prime,(\nu)}$ where the factor indexed by~$\nu$ lifts~$\ul{X}^{\prime,(\nu)}$. In particular this gives a decomposition
$$
T_p(X^\can) = \dirsum_{\nu \in \mQ} T_p^{(\nu)}\, .\leqno(\refn{ptypefA}.1)
$$\vskip-\lastskip\smallskip

\sssection{ConjugLem}
{\it Lemma. --- Possibly after changing $\alpha \colon \Lambda_0 \isomarrow T_p$ by an element of~$\cG^0(\Zp)$, the decomposition\/~{\rm (\refn{ptypefA}.1)} agrees with the decomposition of~$\Lambda_0$ into eigenspaces with respect to the quasi-cocharacter~$N(\mu)$.}
\Dskip

\Proof Let $Q$ be the fraction field of~$W(K)$. Consider the category $\MF^\adm_{Q}(\phi)$ of admissible filtered modules over~$Q$. It is a neutral Tannakian category over~$\Qp$. Let $M$ be the Dieudonn\'e module of~$\ul{X}^\can$. Then $M \otimes \Qp$, equipped with its Frobenius automorphism and Hodge filtration is an object of $\MF^\adm_{Q}(\phi)$. Write $\langle M_\Qp\rangle^\otimes \subset \MF^\adm_{Q}(\phi)$ for the tensor subcategory that it generates. We have a diagram as follows.
$$
\matrix{
\Rep_\Qp(\cG^0)&\mapright{u}&\langle M_\Qp\rangle^\otimes&\mapright{\omega_\Gal}&\Vec_{\Qp}\cr
\noalign{\smallskip}
&&\mapdown{\omega}&&\cr
\noalign{\smallskip}
&&\Vec_{Q}&&\cr
}
$$
Explanation: (a)~We can view $M_\Qp$ as an object with $\cG^0$-structure in~$\MF^\adm_{Q}(\phi)$; this means precisely that we have a tensor functor~$u$ as in the diagram. This functor~$u$ sends the tautological representation of~$\cG^0$ on $\Lambda_0 \otimes \Qp$ to the object~$M \otimes \Qp$. (b)~The functor~$\omega$ is the fibre functor that sends an object of~$\MF^\adm_{Q}(\phi)$ to the underlying $Q$-vector space. (c)~The functor~$\omega_\Gal$ is the functor that sends an object~$(L,\phi,\Fil^\gdot)$ to the $\Qp$-vector space $\Fil^0(L \otimes_Q B_\cris)^{\phi=1}$. We refer to the paper~[\ColFon] of Colmez and Fontaine for further details.

We have a $\mQ$-grading on the object~$M_\Qp$, coming from the decomposition of~$\ul{X}^\can$ as a product of isotypical factors. This gives us a $\mQ$-grading on the functor~$u$. The induced $\mQ$-gradings on the functors $\omega \circ u$ and $\omega_\Gal \circ u$ give rise to quasi-cocharacters
$$
\gamma\colon \tilde\mG_m \tto \Aut^\otimes(\omega \circ u)
\qquad\hbox{and}\qquad
\gamma_\Gal\colon \tilde\mG_m \tto \Aut^\otimes(\omega_\Gal \circ u)\, .
$$
(See Saavedra Rivano~[\Saav], \Romno 4, \S~1.) Here $\tilde\mG_m$ is the pro-algebraic torus with character group~$\mQ$. 

The choice of an identification $\alpha \colon \Lambda_0 \isomarrow T_p$ gives an isomorphism of the functor $\omega_\Gal \circ u$ with the forgetful functor $\Rep_\Qp(\cG^0) \to \Vec_\Qp$, and this induces an isomorphism $\cG^0 \cong \Aut^\otimes(\omega_\Gal \circ u)$. The conjugacy class~$[\gamma_\Gal]$ of the quasi-cocharacter~$\gamma_\Gal$ that we get does not depend on the choice of~$\alpha$ (within its $\cG^0$-conjugacy class). By construction, if we let $\cG^0$ act on~$T_p$ via~$\alpha$ then the decomposition $T_p = \oplus T_p^{(\nu)}$ is precisely the eigenspace decomposition with respect to~$\gamma_\Gal$. 

It now suffices to prove that the quasi-cocharacter~$N(\mu)$ is geometrically conjugate to~$\gamma_\Gal$. Indeed, if this holds then $\gamma_\Gal$ and~$N(\mu)$, being both defined over~$\Qp$, are already in the same $\cG^0(\Qp)$-conjugacy class (use [\KottTOI], Lemma~1.1.3), and this gives our claim.

The automorphism group $\Aut^\otimes(\omega \circ u)$ is an inner form of~$\cG^0_Q$, so over an algebraic closure~$\overline{Q}$ we have an isomorphism $\Aut^\otimes(\omega \circ u) \otimes_Q \overline{Q} \cong \cG^0 \otimes_\Qp {\overline{Q}}$, canonical up to inner automorphisms. In particular, $\gamma$ gives a conjugacy class~$[\gamma]$ of quasi-cocharacters of~$\cG^0$ over~$\overline{Q}$. The two conjugacy classes $[\gamma]$ and~$[\gamma_\Gal]$ are (geometrically) the same; this results from the fact that they are both obtained from the same $\mQ$-grading of the functor~$u$. But $[\gamma]$~represents the Newton quasi-cocharacter associated to~$\ul{X}$. Because~$\ul{X}$ is ordinary, $N(\mu)$~is in this conjugacy class and we are done. \QED

\sssection{ptypefB}
Let $K$ be a field containing $\kappa(v^0)$. Let $f\colon \ul{A}_1 \to \ul{A}_2$ be a $K$-valued point of~$\cJ^\ord$. Choose isomorphisms $\alpha_i \colon \Lambda_0 \isomarrow T_p(X_i^\can)$ as in the lemma. The canonical lifting $f^\can \colon \ul{X}_1^\can \to \ul{X}_2^\can$ respects slope decompositions. Because $f$ is an isogeny, the induced map $V_pf^\can\colon V_pX_1^\can \to V_pX_2^\can$ is an isomorphism. Hence $\alpha_2^{-1} \circ V_p f^\can \circ \alpha_1$ is an element of~$\cM(\Qp)$. Its class in $\cM(\Zp)\backslash \cM(\Qp)/\cM(\Zp)$ is independent of choices. We call
$$
\tau_p(f) := \left[\alpha_2^{-1} \circ V_p f^\can \circ \alpha_1\right] \in \cM(\Zp)\backslash \cM(\Qp)/\cM(\Zp)
$$
the {\it $p$-type\/} of~$f$.

\sssection{ptypeLC}
{\it Lemma. --- Let $T$ be a reduced scheme over $\kappa(v^0)$. If $f \colon T \to \cJ^\ord_0$ is a $T$-valued point of~$\cJ^\ord_0$ then the map $t \mapsto \tau_p(f_t)$ that associates to~$t \in T$ the $p$-type of~$f_t$ is locally constant.}
\Dskip

\SketchProof~First one shows that if $f\colon \ul{A}_1 \to \ul{A}_2$ is a $p$-isogeny over a field then its $p$-type is completely determined by the structure of~$\Ker(f)$. Next remark that it suffices to test local constancy on schemes $T = \Spec(R)$ with $R$ a discrete valuation ring of equal characteristic~$p$. Let $k \subset R$ be a coefficient field. Choose an integer~$N$ large enough such that $\Ker(f) \subseteq X_1[p^N]$. Possibly after passing to a finite extension of~$R$ we can assume that $X_1[p^N]$ is trivialized, meaning that we have an isomorphism $\ul{X}_1[p^N] \cong \ul{Z}[p^N] \otimes_k R$ with $\ul{Z}$ a standard ordinary BT with $(\cO,\ast,\epsilon)$-structure over~$k$. Suppose $m \in \cM(\Qp) \cap \End(\Lambda_0)$ represents the $p$-type of the isogeny~$f$ over the generic point $\eta \in \Spec(R)$. Then $\Ker(f)_\eta$, viewed as a subgroup scheme of $\ul{Z}[p^N] \otimes_k k(\eta)$, can be expressed directly in terms of~$m$, and we find that $\Ker(f)_\eta = H(m) \otimes_k k(\eta)$ for some subgroup scheme $H(m) \subset \ul{Z}[p^N]$ depending only on~$m$. But then $H(m) \otimes_k R$ is the unique flat subgroup scheme of $\ul{Z}[p^n] \otimes_k R$ extending~$\Ker(f)_\eta$. Hence via the chosen trivialization of~$X_1[p^N]$ we have $\Ker(f) = H(m) \otimes_k R$, from which it follows that the $p$-type of~$f$ over the special fibre is also given by the double coset~$[m]$. \QED     

\sssection{hbarDef}
The $p$-type of isogenies allows us to define a map
$$
\bar h\colon \cH_0(\cM,\mQ) \to \mQ[\cJ^\ord_0]\, ,
$$
sending the characteristic function of a double coset~$[m]$ with $m \in \cM(\Qp) \cap \End(\Lambda_0)$ to the sum of all irreducible components of~$\cJ^\ord_0$ on which the $p$-type is equal to~$[m]$.

It should be noted that in the Siegel modular case this definition agrees with the one given by Chai and Faltings in~[\FaCh], Chap.~\Romno 7, \S~4; the normalization factor $1/\# \Sym^2\big(\mZ_p^g/d\mZ_p^g\big)$ used in loc.\ cit.\ (page~261) arises only to compensate for the difference between ``irreducible components'' and ``connected components''. Although it is presumably true that $\bar h$ is a ring homomorphism, we did not check this. Fortunately we only need this property on the image of the twisted Satake homomorphism, where it follows from the commutativity of the diagram below. 

\sssection{CommSqThm}
{\it Theorem. --- Let $\sigma\colon \mQ[J] \tto \mQ[\cJ^\ord_0]$ be the homomorphism given by specialization of cycles. Then we have a commutative diagram of $\mQ$-algebra homomorphisms
$$
\matrix{
\cH_0(\cG^0,\mQ)&\sizedmapright{hhhh}{h}&\mQ[J]\cr
\noalign{\smallskip}
\mapdownl{\dot S^{\cG^0}_\cM}&&\mapdownr{\sigma}\cr
\noalign{\smallskip}
\cH_0(\cM,\mQ)&\sizedmapright{hhhh}{\bar h}&\mQ[\cJ^\ord_0]\cr}
$$}
\vskip-\lastskip\Dskip

The proof of this result is essentially the same as in the Siegel modular case; see Chai and Faltings~[\FaCh], page~263. 

\sssection{H(F)=0}
{\it Corollary. --- Let $\Phi$ be the Frobenius correspondence on~$\cS_0$, as in\/~{\rm\refn{PhionS}}. Let $H_{(\cG^0,\cX^0)} \in \cH_0(\cG^0,\mQ)\big[t\big]$ be the Hecke polynomial associated to the datum $(\cG^0,\cX^0)$, as defined in\/~{\rm [\WedhCong]}, Section~2. Regarding $\mQ[\cJ_0^\ord]$ as an algebra over $\cH(\cG^0,\mQ)$ via $\sigma \circ h$, we have the relation $H_{(\cG^0,\cX^0)}(\Phi)=0$.}
\Dskip

\Proof~As in Wedhorn's paper~[\WedhCong], this is a direct consequence of the theorem together with the purely group-theoretic result loc.\ cit., Prop.~(2.9), due to B\"ultel. \QED

\sssection{IfRelDense}
{\it Corollary. --- If $\cJ^\ord_0$ is Zariski dense in~$\cJ_0$ then the relation $H_{(\cG^0,\cX^0)}(\Phi)=0$ holds in the algebra $\mQ[\cJ_0]$, viewed as an algebra over $\cH_0(\cG^0,\mQ)$ via~$\sigma \circ h$.}

\sssection{CohAppl}
For cohomological applications it is the latter result that is most interesting. The condition that $\cJ^\ord_0$ is Zariski dense in~$\cJ_0$ is referred to as the {\it relative density condition\/} in~[\WedhCong]. Note that even though $\cS_0^\ord$ is dense in~$\cS_0$, the analogous property for~$\cJ_0$ may fail; see for instance Stamm~[\Stamm]. In the cases where the relative density condition is satisfied, Cor.~\refn{IfRelDense} proves the conjecture formulated by Blasius and Rogawski in~[\BlaRog], Section~6. 

It is not clear to the author whether it is reasonable to expect that $H_{(\cG^0,\cX^0)}(\Phi)=0$ if the relative density condition fails. If in the polynomial $H_{(\cG^0,\cX^0)}$ we replace all coefficients (viewed as elements of the algebra $\mQ[\cJ_0]$) by the Zariski closures of their ordinary part then we obtain a polynomial $H^\prime_{(\cG^0,\cX^0)}$ for which $H^\prime_{(\cG^0,\cX^0)}(\Phi) = 0$ by~\refn{H(F)=0}. Therefore the question is whether for the difference $H^\pprime_{(\cG^0,\cX^0)} := H_{(\cG^0,\cX^0)} - H^\prime_{(\cG^0,\cX^0)}$ we again have $H^\pprime_{(\cG^0,\cX^0)}(\Phi)=0$. It is not so clear to us why this should hold. Note that the coefficients of $H^\pprime_{(\cG^0,\cX^0)}$ are linear combinations of irreducible components of~$\cJ_0$ that are {\it not\/} in the closure of the ordinary locus. 

None the less, if we test this in the Hilbert modular case at inert primes (where the relative density condition fails) then it is still true that $H_{(\cG^0,\cX^0)}(\Phi)=0$. The point is that the minimum polynomial is in this case only a factor (of degree~$2$) of the Hecke polynomial (cf.\ Example~(2.13) in~[\WedhCong]), and the coefficients of this mimimum polynomial do not have truly ``non-ordinary'' terms. It would be interesting to investigate whether this is a general phenomenon.

\Askip\noindent
{\sectitlefont References}\writetocentry{nnsect}{References}
\nobreak\bigskip\nobreak
{\eightpoint
\item{SGA 3}
M.~Demazure et al., {\it Sch\'emas en groupes, \Romno1, \Romno2, \Romno3\/}, 
Lect.\ Notes in Math.\ {\bf 151}, {\bf 152}, {\bf 153}, 
Springer-Verlag, Berlin, 1970. 

\item{SGA 7}
A.~Grothendieck et al., {\it Groupes de monodromie en g\'eom\'etrie 
alg\'ebrique\/}, Lect.\ Notes in Math.\ {\bf 288}, {\bf 340},
Springer-Verlag, Berlin, 1972.
\smallskip

\item{[\BlaRog]}
D.~Blasius and J.D.~Rogawski, {\it Zeta functions of Shimura varieties\/}, in: 
Motives (U.~Jannsen, S.~Kleiman, J-P.~Serre, eds.), Proc.\ Symp.\ Pure Math.\ 
{\bf 55}, part~2, AMS, Providence, 1994, pp.~525--571.  

\item{[\ColFon]}
P.~Colmez, J.-M.~Fontaine, {\it Construction des repr\'esentations $p$-adiques
semi-stables\/}, Invent.\ math.\ {\bf 140} (2000), 1--43.

\item{[\Conr]}
B.~Conrad, {\it Background notes on $p$-divisible groups over local fields\/}, 
unpublished manuscript, available at 
{\tt www-math.mit.edu/\kern-1pt\lower 2pt\hbox{\~{\ }}\kern-4pt dejong}~.

\item{[\DelCoord]}
P.~Deligne, {\it Cristaux ordinaires et coordonn\'ees canoniques\/}; with the 
collaboration of L.~Illusie; with an appendix by N.M.~Katz; in: Surfaces 
Alg\'ebriques (J.~Giraud, L.~Illusie, M.~Raynaud, eds.), Lect.\ Notes in Math.\ 
{\bf 868}, Springer-Verlag, Berlin, 
1981, pp.~80--137.

\item{[\Falt]}
G.~Faltings, {\it Integral crystalline cohomology over very ramified valuation 
rings\/}, J.A.M.S.\ {\bf 12} (1999), 117--144.

\item{[\FaCh]}
G.~Faltings and C.-L.~Chai, {\it Degeneration of abelian varieties\/}, Ergebnisse 
der Math., 3.~Folge, {\bf 22}, Springer-Verlag, Berlin, 1990.

\item{[\Font]}
J.-M.~Fontaine, {\it Groupes $p$-divisibles sur les corps locaux\/}, 
Ast\'erisque {\bf 47--48} (1977).  

\item{[\GorOo]}
E.Z.~Goren and F.~Oort, {\it Stratifications of Hilbert modular varieties\/}, J.\ 
Alg.\ Geom.\ {\bf 9} (2000), 111--154.

\item{[\Illus]}
L.~Illusie, {\it D\'eformations de groupes de Barsotti-Tate (d'apr\`es 
A.~Grothendieck)\/}, in: S\'emi\-naire sur les pinceaux arithm\'etiques: la 
conjecture de Mordell (L.~Szpiro, ed.), Ast\'erisque {\bf 127} (1985), 151--198.

\item{[\Katz]}
N.M.~Katz, {\it Slope filtration of $F$-crystals\/}, Ast\'erisque {\bf 63}
(1979), 113--163.

\item{[\KatzST]}
N.M.~Katz, {\it Serre-Tate local moduli\/}, in: Surfaces Alg\'ebriques 
(J.~Giraud, L.~Illusie, M.~Raynaud, eds.), Lect.\ Notes in Math.\ {\bf 868}, 
Springer-Verlag, Berlin, 1981, pp.~138--202.

\item{[\Knus]}
M.-A.~Knus, {\it Quadratic and hermitian forms over rings\/}, Grundlehren der 
math.\ Wiss.\ {\bf 294}, Springer-Verlag, Berlin, 1991.

\item{[\KottTOI]}
R.E.~Kottwitz, {\it Shimura varieties and twisted orbital integrals\/}, Math.\ 
Ann.\ {\bf 269} (1984), 287--300. 

\item{[\KottIsoc]}
R.E.~Kottwitz, {\it Isocrystals with additional structure\/}, Compos.\ Math.\ 
{\bf 56} (1985), 201--220; {\it ---~\Romno 2\/}, Compos.\ Math.\ {\bf 109} 
(1997), 
255--339.

\item{[\Kott]}
R.E.~Kottwitz, {\it Points on some Shimura varieties over finite fields\/}, 
J.A.M.S.\ {\bf 5} (1992), 373--444.

\item{[\Kraft]}
H.~Kraft, {\it Kommutative algebraische $p$-Gruppen (mit Anwendungen auf 
$p$-divisible Gruppen und abelsche Variet\"aten)\/}, manuscript, Univ.\ Bonn, 
Sept.\ 1975, 86 pp. (Unpublished)

\item{[\LST]}
J.~Lubin, J-P.~Serre and J.~Tate, {\it Elliptic curves and formal groups\/}, 
Woods Hole Summer Institute, 1964. (Mimeographed notes.) Available at 
{\tt www.ma.utexas.edu/users/voloch/lst.html}~.

\item{[\Manin]}
Yu.I.~Manin, {\it The theory of commutative formal groups over fields of finite 
characteristic\/}, Uspehi Mat.\ Nauk.\ {\bf 18} (1963), 3--90; English 
translation: Russian Math.\ Surv.\ {\bf 18} (1963), 1--83.

\item{[\Mess]}
W.~Messing, {\it The crystals associated to Barsotti-Tate groups: with 
applications to abe\-lian sche\-mes\/}, Lect.\ Notes in Math.\ {\bf 264}, 
Springer-Verlag, Berlin, 1972.

\item{[\Milne]}
J.S.~Milne, {\it Shimura varieties and motives\/}, in: Motives (U.~Jannsen, 
S.~Kleiman, J-P.~Serre, eds.), Proc.\ Symp.\ Pure Math.\ {\bf 55}, part~2, AMS, 
Providence, 1994, pp.~447--523. 

\item{[\Durh]}
B.J.J.~Moonen, {\it Models of Shimura varieties in mixed characteristics\/}, in: 
Galois Representations in Arithmetic Algebraic Geometry (A.J.~Scholl, R.~Taylor, 
eds.), London Math. Soc., Lecture Notes Series {\bf 254}, Cambridge Univ. Press, 
Cambridge, 1998, pp.\ 271--354.

\item{[\GSAS]}
B.J.J.~Moonen, {\it Group schemes with additional structures and Weyl group 
cosets\/}, in: Moduli of Abelian Varieties (C.~Faber, G.~van der Geer, F.~Oort, 
eds.), Progr.\ Math.\ {\bf 195}, Birkh\"auser, Basel, 2001, pp.\ 255--298.

\item{[\DFEO]}
B.J.J.~Moonen, {\it A dimension formula for Ekedahl-Oort strata\/}, math.AG/0208161.

\item{[\MAV]}
D.~Mumford, {\it Abelian varieties\/}, Oxford Univ.\ Press, Oxford, 1970.

\item{[\Noot]}
R.~Noot, {\it Models of Shimura varieties in mixed characteristic\/}, J.\ Alg.\ 
Geom.\ {\bf 5} (1996), 187--207.

\item{[\FOTexel]}
F.~Oort, {\it A stratificiation of a moduli space of abelian varieties\/}, in:
Moduli of Abelian Varieties (C.~Faber, G.~van der Geer, F.~Oort, eds.), Progr.\ 
Math.\ {\bf 195}, Birkh\"auser, Basel, 2001, pp.\ 345--416.

\item{[\FOTexelB]}
F.~Oort, {\it Newton polygon strata in the moduli space of abelian varieties\/}, 
in:
Moduli of Abelian Varieties (C.~Faber, G.~van der Geer, F.~Oort, eds.), Progr.\ 
Math.\ {\bf 195}, Birkh\"auser, Basel, 2001, pp.\ 417--440.

\item{[\RapBourb]}
M.~Rapoport, {\it On the Newton stratification\/}, S\'em.\ Bourbaki, Mars 2002.

\item{[\RR]}
M.~Rapoport and M.~Richartz, {\it On the classification and specialization of 
$F$-isocrystals with additional structure\/}, Compos.\ Math.\ {\bf 103} (1996), 
153--181.

\item{[\Rayn]}
M.~Raynaud, {\it Sch\'emas en groupes de type $(p,\ldots,p)$\/}, Bull.\ Soc.\ 
math.\ France {\bf 102} (1974), 241--280.

\item{[\ReimZink]}
H.~Reimann und Th.~Zink, {\it Der Dieudonn\'emodul einer polarisierten abelschen
Mannigfaltigkeit vom CM-Typ\/}, Ann.\ Math.\ {\bf 128} (1988), 461--482.

\item{[\Saav]}
N.~Saavedra Rivano, {\it Cat\'egories Tannakiennes\/}, Lect.\ Notes in Math.\ 
{\bf 265}, Springer-Verlag, Berlin, 1972.

\item{[\Stamm]}
H.~Stamm, {\it On the reduction of the Hilbert-Blumenthal-moduli scheme with 
$\Gamma_0(p)$-level structure\/}, Forum Math.\ {\bf 9} (1997), 405--455. 

\item{[\WedhOrd]}
T.~Wedhorn, {\it Ordinariness in good reductions of Shimura varieties of 
PEL-type\/}, Ann.\ scient.\ \'Ec.\ Norm.\ Sup.\ (4) {\bf 32} (1999), 575--618.

\item{[\WedhCong]}
T.~Wedhorn, {\it Congruence relations on some Shimura varieties\/}, J.\ reine 
angew.\ Math.\ {\bf 524} (2000), 43--71. 

\item{[\Wedh]}
T.~Wedhorn, {\it The dimension of Oort strata of Shimura varieties of PEL-type\/}, 
in: Moduli of Abelian Varieties (C.~Faber, G.~van der Geer, F.~Oort, eds.), 
Progr.\ Math.\ {\bf 195}, Birkh\"auser, Basel, 2001, pp.\ 441--471.
\Bskip 

\noindent
Ben Moonen, University of Amsterdam, Korteweg-de Vries Institute for 
Mathematics, Plantage Muidergracht~24, 1018~TV Amsterdam, The Netherlands. 
Email: {\tt bmoonen@science.uva.nl}\par}
\bye